\documentclass{ws-m3as}

\usepackage{todonotes}
\usepackage{amsfonts}
\usepackage{graphicx}
\usepackage{epstopdf}
\usepackage{algorithmic}
\usepackage{diagbox}
\usepackage{paralist}
\usepackage{overpic}
\usepackage{mathtools}
\usepackage{booktabs}
\usepackage{enumitem}
\usepackage{url}

\newcommand{\inner}[2]{\langle #1, #2 \rangle}

\newcommand{\R}{\mathbb{R}}

\newcommand{\dx}{\,\mathrm{d}x}

\newcommand{\K}{\text{K}}
\newcommand{\Na}{\text{Na}}
\newcommand{\Cl}{\text{Cl}}

\newcommand{\NaP}{\text{NaP}}
\newcommand{\GHK}{\text{GHK}}

\newcommand{\KIR}{\text{KIR}}
\newcommand{\KDR}{\text{KDR}}
\newcommand{\KA}{\text{KA}}
\newcommand{\ATP}{\text{ATP}}
\newcommand{\norm}[1]{\lVert #1 \rVert}
\newcommand{\T}[1]{\text{#1}}

\DeclareMathOperator{\Grad}{\nabla}

\begin{document}

\markboth{A. J. Ellingsrud, N. Boull\'e, P. E. Farrell, and M. E. Rognes}{Numerics for brain electrochemomechanics}

%
\catchline{}{}{}{}{}
%

\title{Accurate numerical simulation of electrodiffusion and water
movement in brain tissue}
    

\author{Ada Johanne Ellingsrud}

\address{Simula Research Laboratory, Oslo, Norway. \\ ada@simula.no}

\author{Nicolas Boull\'e}
\address{Mathematical Institute, University of Oxford, Oxford, UK. \\ boulle@maths.ox.ac.uk}

\author{Patrick E.~Farrell}
\address{Mathematical Institute, University of Oxford, Oxford, UK. \\ patrick.farrell@maths.ox.ac.uk}

\author{Marie E.~Rognes}
\address{Simula Research Laboratory, Oslo, Norway. \\ meg@simula.no}

\maketitle

\begin{history}
\end{history}

\begin{abstract}
    Mathematical modelling of ionic electrodiffusion and water movement
    is emerging as a powerful avenue of investigation to provide new
    physiological insight into brain homeostasis. However, in order to
    provide solid answers and resolve controversies, the accuracy of the
    predictions is essential. Ionic electrodiffusion models typically
    comprise non-trivial systems of non-linear and highly coupled
    partial and ordinary differential equations that govern phenomena on
    disparate time scales. Here, we study numerical challenges related
    to approximating these systems. We consider a homogenized model for
    electrodiffusion and osmosis in brain tissue and present and
    evaluate different associated finite element-based splitting schemes
    in terms of their numerical properties, including accuracy,
    convergence, and computational efficiency for both idealized scenarios
    and for the physiologically relevant setting of cortical
    spreading depression (CSD). We find that the schemes display optimal
    convergence rates in space for problems with smooth manufactured
    solutions. However, the physiological CSD setting is challenging: we
    find that the accurate computation of CSD wave characteristics (wave
    speed and wave width) requires a very fine spatial and fine temporal
    resolution.
\end{abstract}

\keywords{electrodiffusion, osmosis, brain electrophysiology and mechanics,
finite element method, splitting scheme, numerical convergence}

\ccode{AMS Subject Classification: 22E46, 53C35, 57S20}

\section{Introduction}
\label{sec:intro}
The movement of ions and molecules in and between cellular compartments is
fundamental for brain function, and importantly for neuronal excitability and
activity\cite{aitken1986sources,nicholson1978calcium,utzschneider1992mutual}.
Vital processes such as action potential firing, transmitter release, and
synaptic transmission are all driven by ionic gradients across the neuronal
membrane. Regulation of the extracellular volume by cellular swelling is
closely related to ionic dynamics, including potassium
buffering\cite{hertz1965possible,kuffler1966physiological}. Concurrently,
several pathologies are associated with disruption to ionic homeostasis in the
brain, e.g.~Huntington's disease\cite{tong2014astrocyte}, multiple
sclerosis\cite{srivastava2012potassium},
migraine\cite{staehr2019involvement}, epilepsy\cite{kohling2016potassium},
Alzheimer's disease\cite{noh2019transient}, and cortical spreading
depression\cite{pietrobon2014chaos}. Recent research efforts indicate that
ion concentrations in the extracellular space are not static, but vary across states such as locomotion\cite{brocard2013activity} and the sleep cycle\cite{ding2016changes}.

In spite of these aspects, mathematical and numerical models for describing dynamics in
brain tissue traditionally assume that the ion concentrations are constant in
time and space. Although such models have provided valuable insight into the mechanisms underlying excitable cells, they fail to represent essential dynamics related to altered
ion concentrations in brain tissue. Recently, several mathematical models
also including electrodiffusive effects have been
presented\cite{cressman2009influence,ellingsrud2020finite,halnes2016effect,hubel2014dynamics,kager2000simulated,mori2015multidomain,oyehaug2012dependence,pods2013electrodiffusion,saetra2020electrodiffusive,somjen2008computer,ullah2009influence}.
However, to date little attention has been paid to the numerical solution of
these models.

In this paper, we consider a mathematical framework proposed by Mori\cite{mori2015multidomain},
consisting of a system of partial differential equations (PDEs)
governing ionic electrodiffusion and water flow in biological tissue,
coupled to a system of ordinary differential equations (ODEs) describing the
temporal evolution of ionic membrane mechanisms. The
system predicts the dynamics of volume fractions, ion concentrations,
electrical potentials and mechanical pressure in an arbitrary number of
cellular compartments and in the extracellular space (ECS). The cellular
compartments can communicate with the ECS via transmembrane ion and water
fluxes. This mathematical model extends on the celebrated bidomain
model\cite{eisenberg1979electrical,eisenberg1970three,tung1978bi}, and both
represent the tissue in a homogenized manner. Homogenized models are
coarse-grained, and hence well suited for simulating phenomena on the tissue
scale (mm). Importantly, the two models differ in that the classical bidomain
model only predicts electrical potentials, whereas the model for ionic
electrodiffusion and water flow takes into account how the movement of ions
affect the excitable tissue, both in terms of electrochemical and mechanical
effects.

Previously, the electrodiffusive model (in its zero flow limit) has been used
to study dynamics in brain tissue, and in particular cortical spreading
depression (CSD)\cite{o2016effects,tuttle2019computational}.  CSD is a slowly
propagating wave of depolarization of brain cells, characterized by elevated
levels of extracellular potassium, calcium, and glutamate, cellular swelling
and pronounced ECS shrinkage\cite{charles2013cortical,pietrobon2014chaos}.
Importantly, CSD is a fundamental pattern of brain signalling that challenges
ionic homeostasis mechanisms in the brain.  As such, a better understanding of
the sequence of events in CSD has the potential to provide new insight into
underlying processes both in cerebral physiology and
pathology\cite{enger2017deletion}. The aforementioned computational studies have focused
on providing new insight into the role of glial cells in
CSD\cite{o2016effects}, and the role of glutamate dynamics in
CSD\cite{tuttle2019computational}. However, in order to provide true
physiological insight, the accuracy of the numerical and computational
predictions is key; the difference between conflicting experimental
observations may very well be within the numerical error of underresolved
models.

\begin{figure}[h]
        \textbf{A} 
        \vspace{0.25cm}
        \begin{center}
        \begin{overpic}[width=0.13\linewidth]{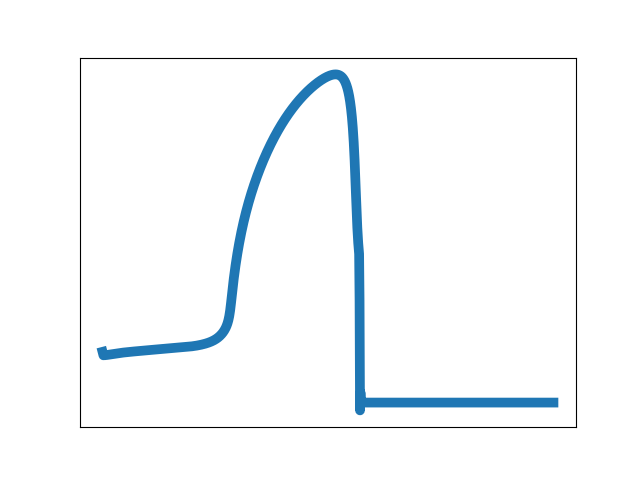}
            \put(-3,82){\color{black}\vector(1,0){50}}
            \put(55,81){Temporal refinement}
            \put(-18,-140){\rotatebox{90}{Spatial refinement}}
            \put(-3,82){\color{black}\vector(0,-1){50}}
        \end{overpic}
        \includegraphics[width=0.13\linewidth]{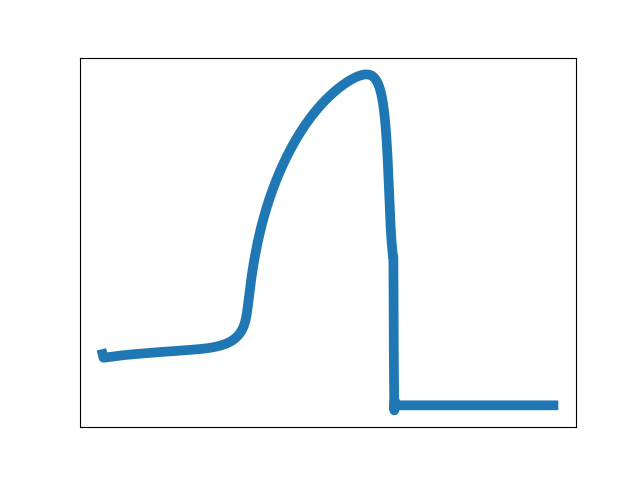}
        \includegraphics[width=0.13\linewidth]{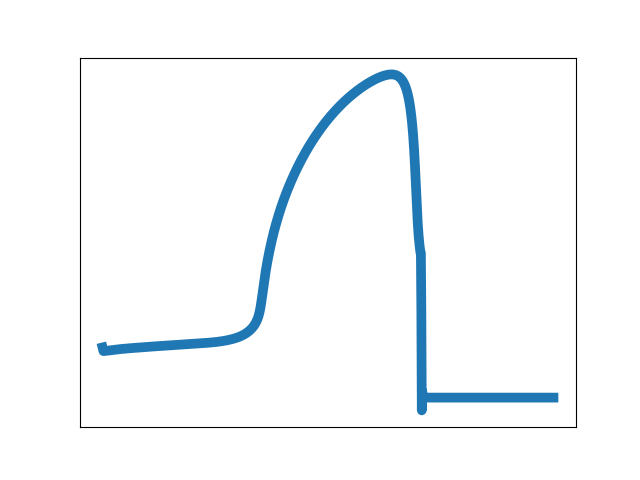}
        \includegraphics[width=0.13\linewidth]{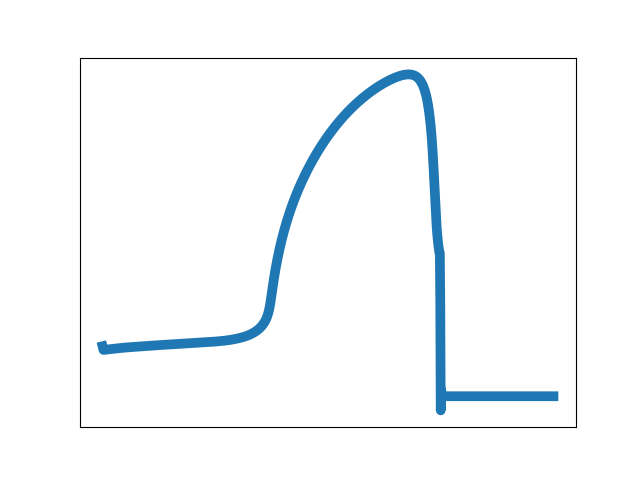}
        \includegraphics[width=0.13\linewidth]{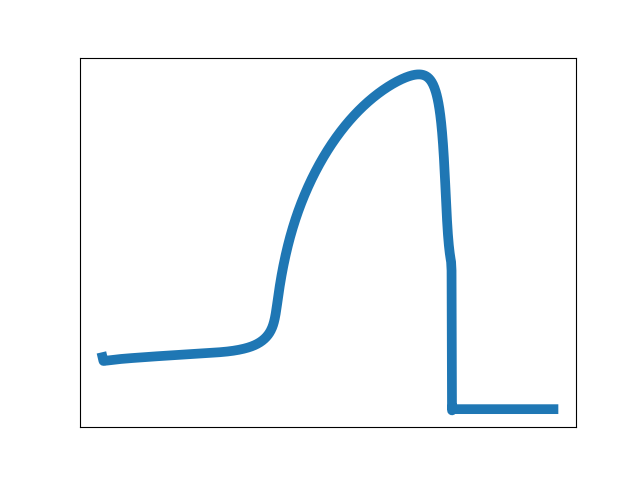}
        \includegraphics[width=0.13\linewidth]{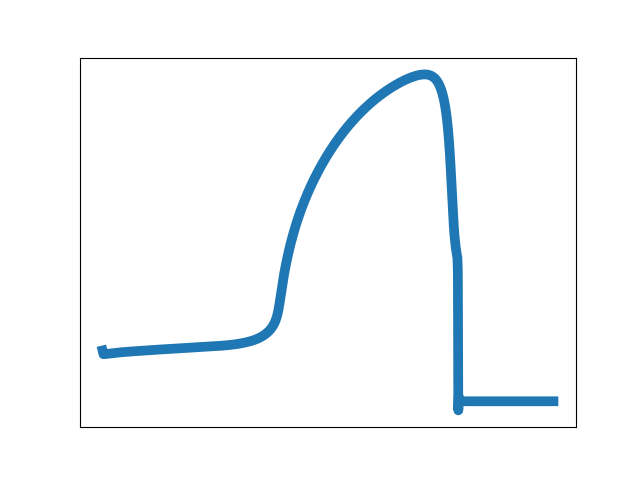}
        \includegraphics[width=0.13\linewidth]{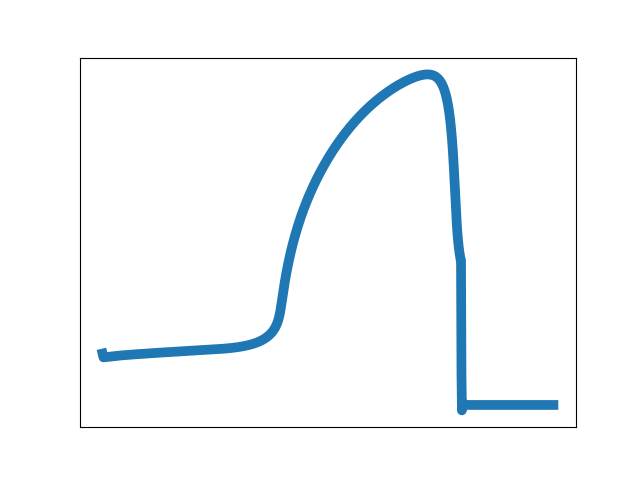}  
         \\

        \includegraphics[width=0.13\linewidth]{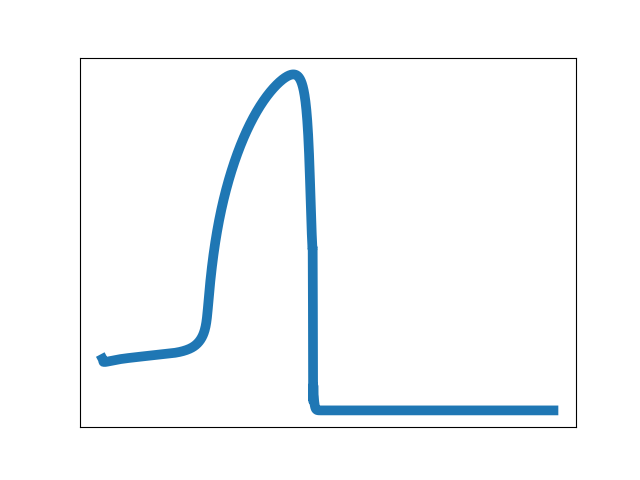}
        \includegraphics[width=0.13\linewidth]{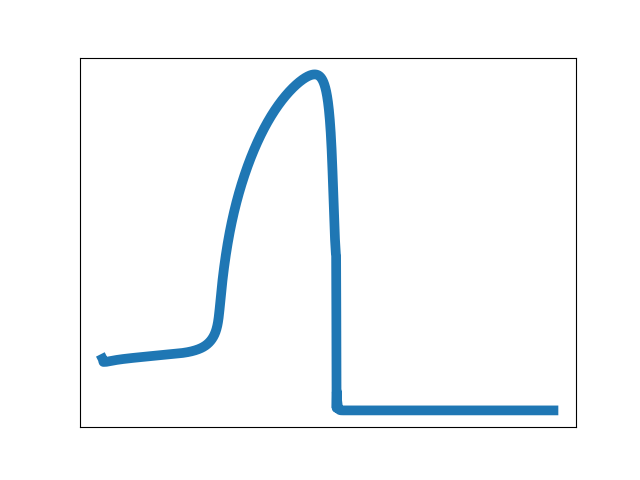}
        \includegraphics[width=0.13\linewidth]{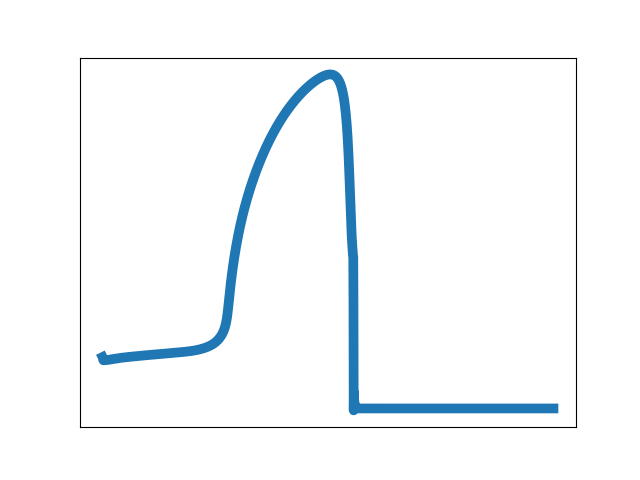}
        \includegraphics[width=0.13\linewidth]{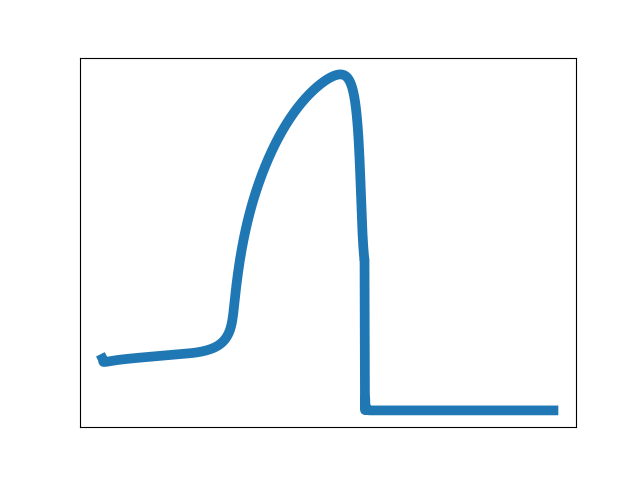}
        \includegraphics[width=0.13\linewidth]{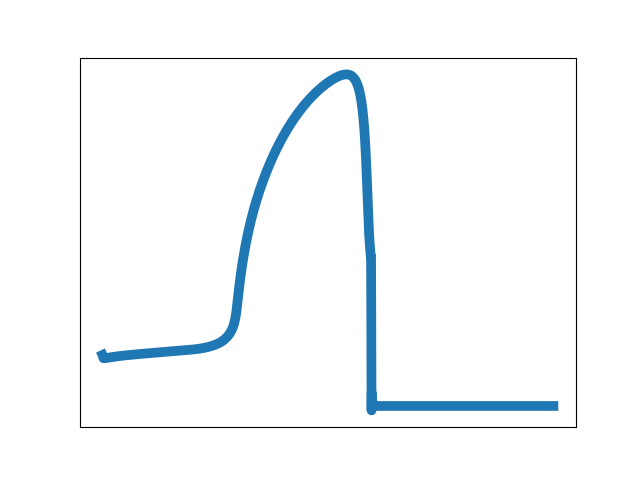}
        \includegraphics[width=0.13\linewidth]{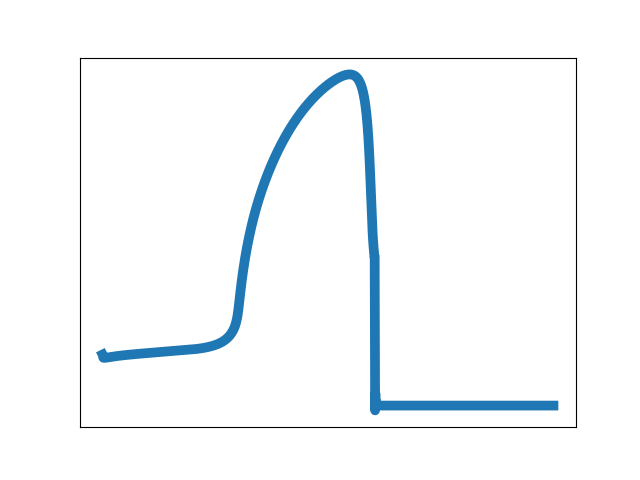}
        \includegraphics[width=0.13\linewidth]{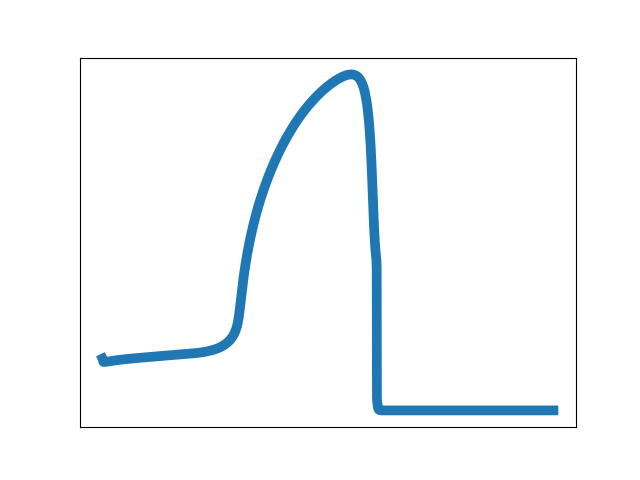} \\

        \includegraphics[width=0.13\linewidth]{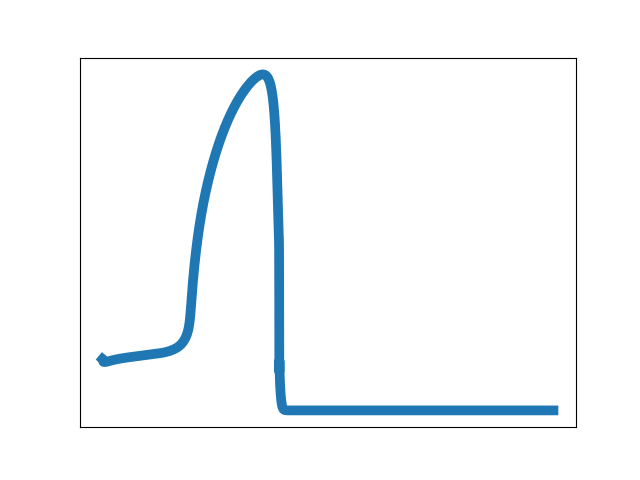}
        \includegraphics[width=0.13\linewidth]{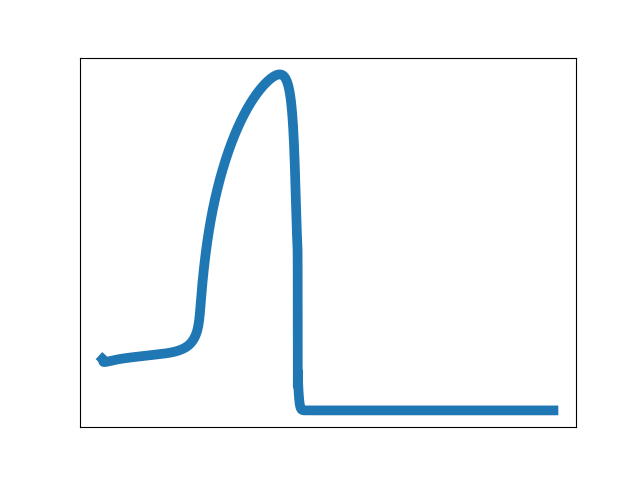}
        \includegraphics[width=0.13\linewidth]{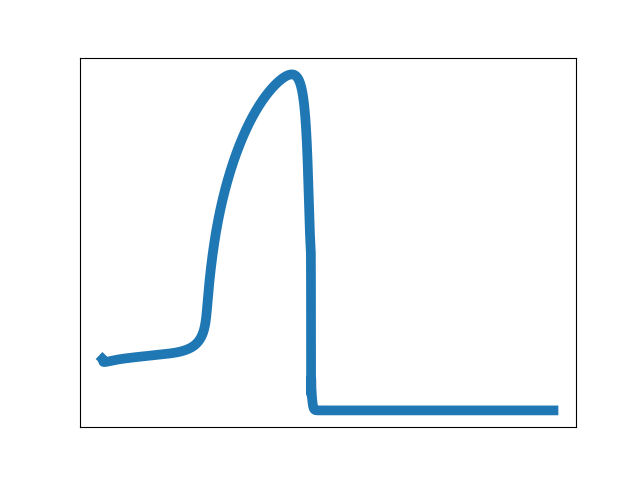}
        \includegraphics[width=0.13\linewidth]{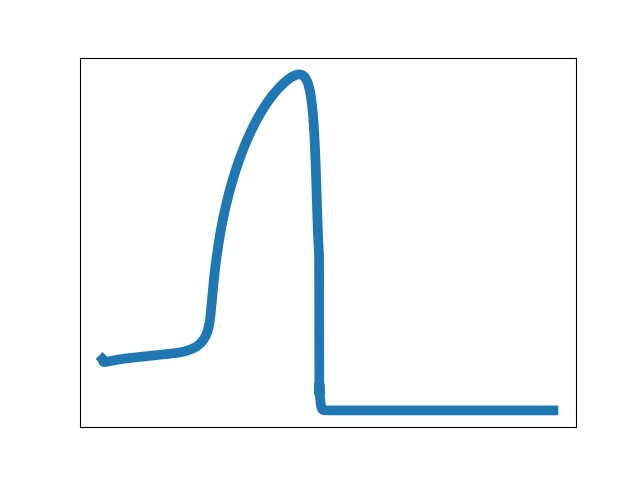}
        \includegraphics[width=0.13\linewidth]{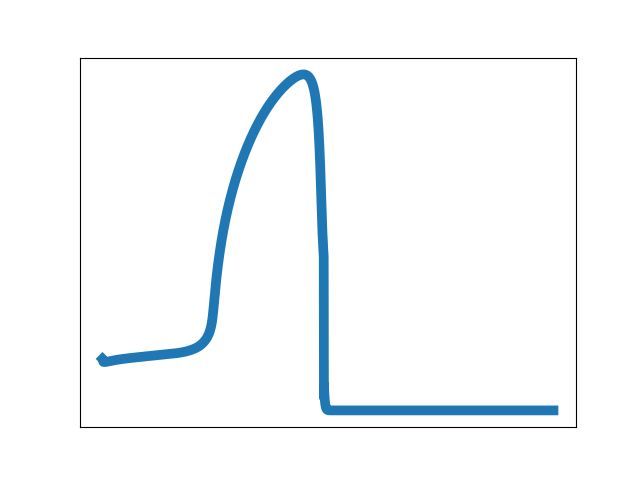}
        \includegraphics[width=0.13\linewidth]{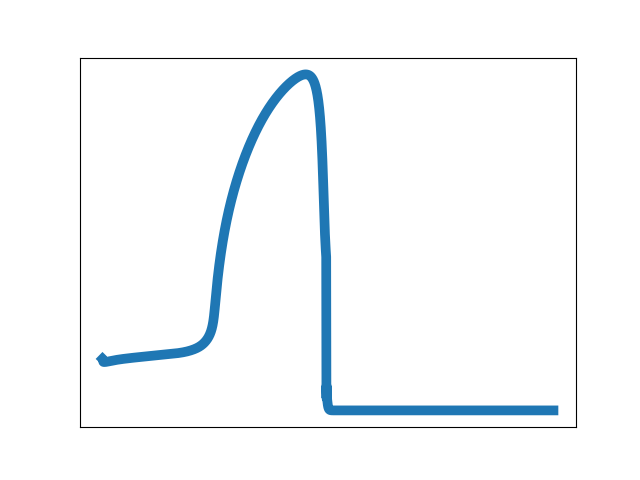}
        \includegraphics[width=0.13\linewidth]{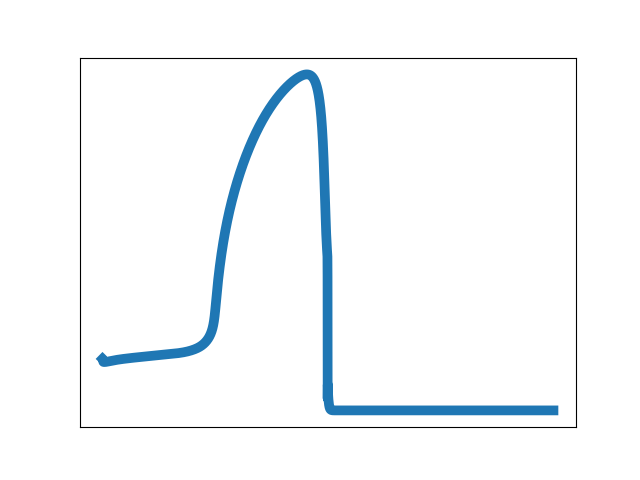} \\

        \includegraphics[width=0.13\linewidth]{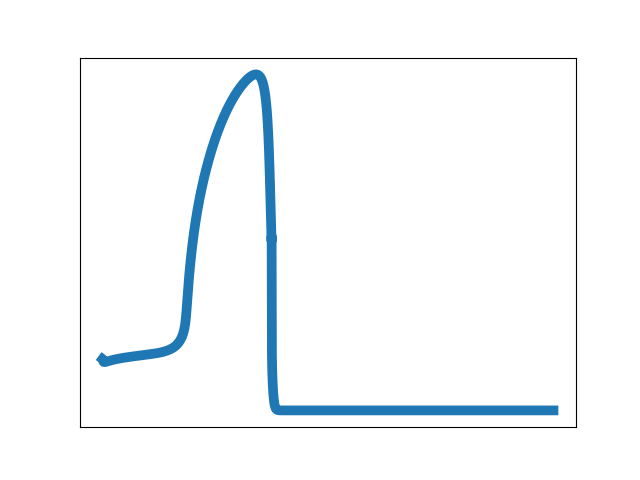}
        \includegraphics[width=0.13\linewidth]{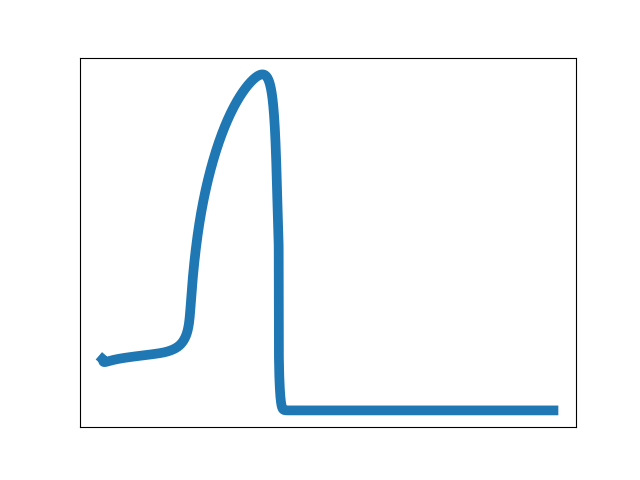}
        \includegraphics[width=0.13\linewidth]{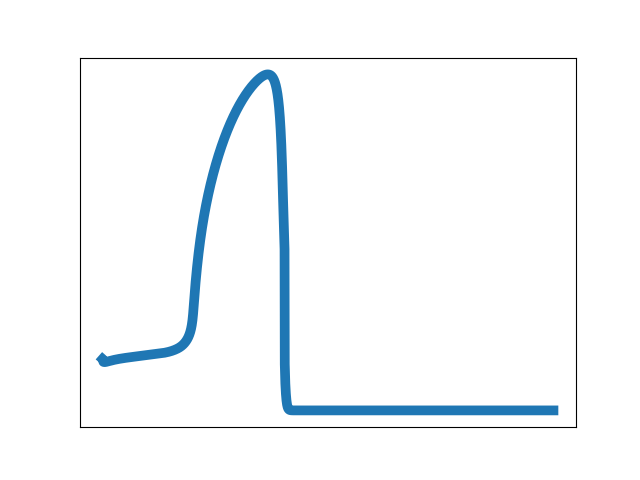}
        \includegraphics[width=0.13\linewidth]{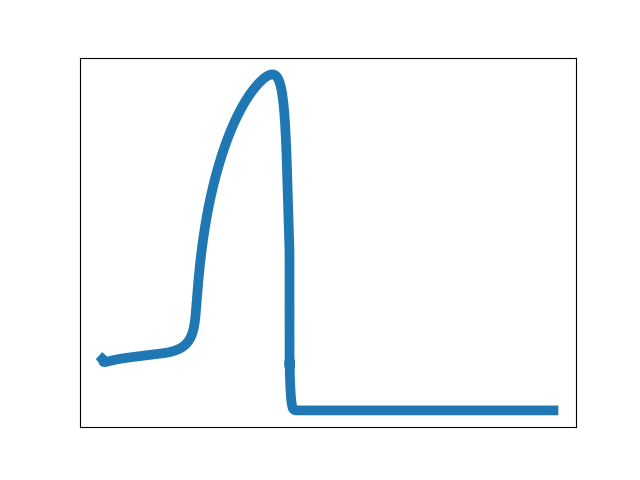}
        \includegraphics[width=0.13\linewidth]{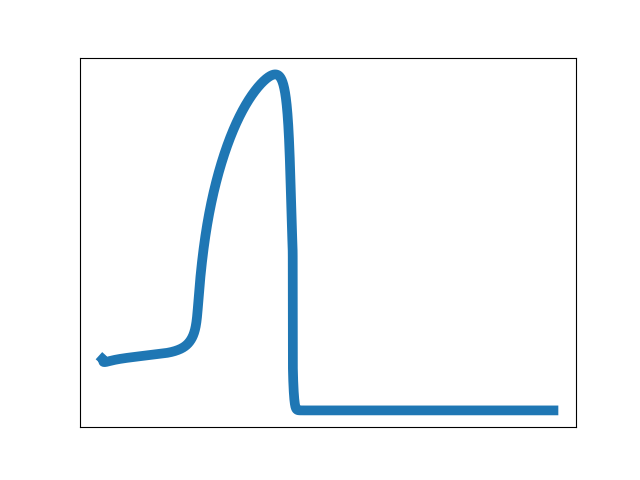}
        \includegraphics[width=0.13\linewidth]{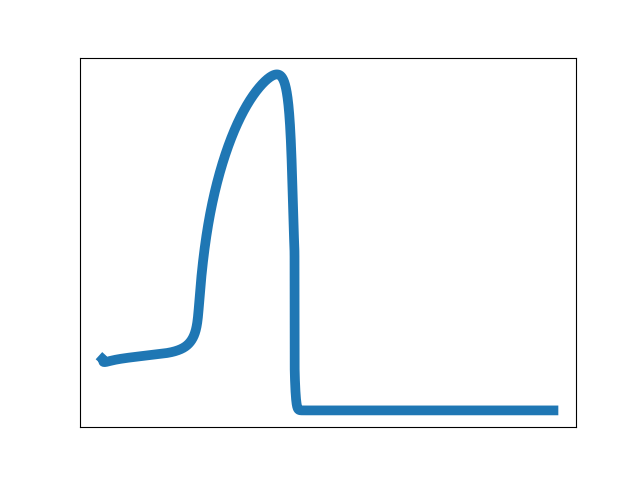}
        \includegraphics[width=0.13\linewidth]{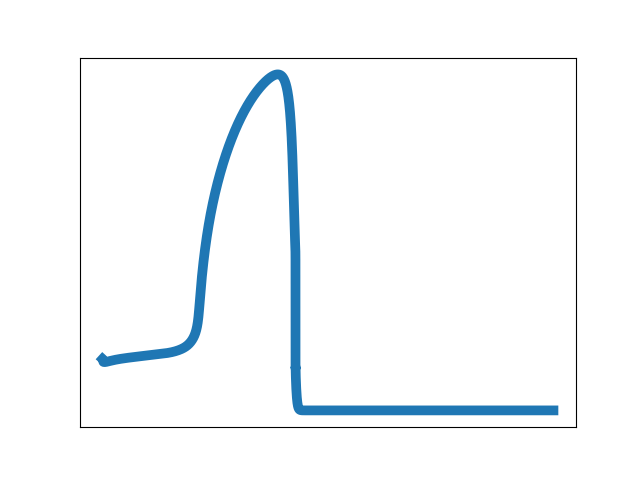}
    \end{center}
    \textbf{B} 
    \begin{center}
    \begin{tabular}{ c | c c c c c c c | c}
        \toprule
        \diagbox{$N$}{$\Delta t$} & 12.5 & 6.25 & 3.125 & 1.563 & 0.781 & 0.391 & 0.195 & $\Delta \bar{v}_{\rm CSD}$ \\
        \midrule
        1000  &7.015 &7.923 &8.677 &9.185 &9.477 &9.662 &9.738 & -- \\
        2000  &5.763 &6.385 &6.846 &7.146 &7.331 &7.423 &7.469 &2.269 \\
        4000  &4.867 &5.361 &5.716 &5.931 &6.054 &6.123 &6.158 &1.312 \\
        8000  &4.688 &4.865 &5.019 &5.147 &5.232 &5.282 &5.305 &0.852 \\
        \midrule
        $\Delta \bar{v}_{\rm CSD}$ & -- &0.178 &0.154 &0.128 &0.085 &0.049 &0.024 & \\
        \bottomrule
    \end{tabular}
    \end{center}
    \caption{Wave properties during refinement in space (N) and time
    ($\Delta$t, ms) in a 1D domain of length $10$ mm at $t=50$ s. The PDEs are
    discretized in time by the previously presented first order scheme
    from Mori, O'Connell, and Tuttle et al.\protect\cite{mori2015multidomain,o2016effects,tuttle2019computational} and
    the ODEs are solved using backward Euler. \textbf{A}: Neuron potential
    $\phi_n(x, 50)$ (mV) versus $x \in \Omega$ (mm). \textbf{B}: CSD mean wave
    speed $\bar{v}_{\rm CSD}$ (mm/min) and difference $\Delta \bar{v}_{\rm
    CSD}$ between consecutive refinements.}
    \label{fig:intro}
\end{figure}

Previously considered numerical schemes for the electrodiffusive model
are based on finite difference\cite{mori2015multidomain,o2016effects}
or finite volume\cite{tuttle2019computational} discretizations in
space and a backward Euler scheme (with explicit treatment of the
active membrane flux) in time. The spatio-temporal discretization
sizes of these schemes are reported to be on the order of $\Delta x
\approx 0.02-0.2$ mm and $\Delta t \approx 10$ ms.  By applying a
(comparable) finite element scheme in space and a similar
discretization scheme in time, we find that the CSD wave properties
change substantially during spatial and temporal refinement. In
particular, we observe that the wave speed and the width of the wave
\emph{increase} with decreasing time resolution, and \emph{decrease}
with decreasing spatial resolution (Figure~\ref{fig:intro}).
As the Mori framework comprises a system of non-linear and highly
coupled partial and ordinary differential equations, governing
phenomena on disparate time scales (both fast electrotonic effects and
much slower effects mediated by diffusion), theoretical analysis of
the full equations is challenging. In this paper, we take a natural
first step towards accurate and efficient approximation schemes for
the system by comparing different schemes numerically, with an
emphasis on convergence and performance for simulating CSD waves.

In the first instance, we consider a low-order finite element scheme in space for
the electrodiffusive model in the zero flow limit, and numerically study the
effect of: (i) different operator splitting schemes, (ii) choice of time
stepping schemes for the PDEs, (iii) choice of time stepping schemes for the
ODEs, and (iv) a higher order spatial discretization scheme. All schemes display
optimal convergence rates during refinement in space and time for problems with
smooth manufactured solutions. However, we find that the accurate computation of
CSD wave characteristics (such as wave speed and wave width) requires a very
fine spatial and fine temporal resolution for all schemes tested. In
particular, the different splitting schemes and PDE time stepping give
comparable results in terms of accuracy. We observe that higher order
ODE time stepping schemes (ESDIRK4, RK4) yield slightly faster convergence
than the lower order backward Euler scheme. Finally, we find that applying a
higher order finite element scheme with a coarser spatial resolution gives
comparable results to the lower order finite element scheme in terms of
accuracy, but at a higher cost, justifying the use of a low order scheme with fine spatial
resolution.

We then turn to consider the full mathematical model, where the
electrochemical and mechanical dynamics described in the zero flow
limit are coupled with microscopic fluid dynamics. To the best of our
knowledge, the only previous study involving this model was presented
by O'Connell\cite{o2016computational}, studying the effects of
mechanical pressures and compartmental fluid flow on CSD wave
characteristics in a setting with two compartments (neurons and
ECS). Here, we present a low order finite element scheme and
simulation results for the full model in three compartments (neurons,
glial and ECS) and study the convergence of the scheme.  The scheme
again displays optimal convergence rates during refinement in space
and time for a problem with smooth manufactured solutions, and meets
similar challenges as the zero flow limit model for the CSD case.

The paper is organized as follows. The Mori framework is summarized
in Section~\ref{sec:mathmodel}. In Section~\ref{sec:scheme}, we present new numerical schemes
for the zero flow limit version of the model, along with numerical convergence
and performance studies for the suggested schemes in Section~\ref{sec:MMS}
and Section~\ref{sec:convergence}. The next sections consider the full model. We present a
numerical scheme in Section~\ref{sec:full:scheme} along with numerical convergence
studies in Section~\ref{sec:full:MMS}, while in Section~\ref{sec:full:CSD} we present results
from a full model simulation of CSD including microscopic fluid dynamics.
Finally, Section~\ref{sec:conclusions} contains a discussion and concluding remarks.
The code used to obtain the simulation results presented within this work is
based on FEniCS\cite{logg2012automated}, and is publicly
available\cite{ellingsrud2021supplementary}.

\section{Mathematical model}
\label{sec:mathmodel}

The tissue of interest is represented as a domain $\Omega \subset
\mathbb{R}^d$, with $d \in \{1,2,3\}$. Moreover, the tissue is composed of $R$
compartments indexed by $r = 1, \dots, R$. We assume that $r=R$
always denotes the ECS, and that other compartments communicate with the ECS
only. Let time $t \in (0, T]$. For a self-contained exposition, we
summarize the Mori framework\cite{mori2015multidomain} below.

\subsection{Governing equations} We will consider the following system of
coupled, time-dependent, nonlinear partial differential equations. Find the
\emph{volume fractions} $\alpha_r: \Omega \times (0, T] \rightarrow [0,1)$ such
that for each $t \in (0, T]$:
\begin{subequations} \label{eq_alpha_rR}
\begin{align}
    \label{eq:alpha:r}
    \frac{\partial \alpha_r}{\partial t} + \nabla \cdot (\alpha_r u_r) &=
    -\gamma_{rR} w_{rR}, \quad \textup{for }  r = 1, \ldots, R-1, \\
    \label{eq:alpha:R}
    \frac{\partial \alpha_R}{\partial t} + \nabla \cdot (\alpha_R u_R) &=
    \sum_{r=1}^{R-1} \gamma_{rR} w_{rR},
\end{align}
\end{subequations}
where $u_r: \Omega \times (0, T] \rightarrow \R^d$ is the \emph{fluid
velocity field} of compartment $r$ (m/s). The \emph{transmembrane water flux} $w_{rR}$ between
compartment $r$ and the ECS is driven by osmotic and oncotic pressure, and will
be discussed further in Section~\ref{sec:sub:memmech}. The coefficient
$\gamma_{rR}$ (1/m) represents the area of cell membrane between compartment
$r$ and the ECS per unit volume of tissue. 
By definition, the volume fractions
sum to $1$, and hence:
\begin{equation}
    \label{eq:alpha:R:sum}
    \alpha_R = 1 - \sum_{r=1}^{R-1} \alpha_r.
\end{equation}
where $\alpha_r \ge 0$. 
Further, for each ion species $k \in K$, find
the \emph{ion concentrations} $[k]_r: \Omega \times (0, T] \rightarrow \R$
(mol/m$^3$) and the \emph{electrical potentials} $\phi_r: \Omega \times (0, T]
\rightarrow \R$ (V) such that for each $t \in (0, T]$:
\begin{subequations} \label{eq_conservation_rR}
\begin{align}
    \label{eq:concentration:r}
    \frac{\partial (\alpha_r [k]_r)}{\partial t} + \nabla \cdot J_r^k &= -
    \gamma_{rR} {J}_{rR}^k, \quad r = 1, \ldots, R-1, \\
    \label{eq:concentration:R}
     \frac{\partial (\alpha_R [k]_R)}{\partial t} +  \nabla \cdot J_R^k&= 
     \sum_{r=1}^{R-1} \gamma_{rR} {J}_{rR}^k,
\end{align}
\end{subequations}
where $J_r^k: \Omega \times (0, T] \rightarrow \R^d$ is the \emph{ion flux
density} (mol/(m$^2$s)) for each ion species $k$. Modelling the
\emph{transmembrane ion flux density} (mol/(m$^2$s)) $J_{rR}^k$ for each ion
species $k \in K$ will be discussed further in Section~\ref{sec:sub:memmech}.
Note that~\eqref{eq_conservation_rR} follow from first principles and express
conservation of ion concentrations for the bulk of each region.  Moreover, we
assume that the ion flux densities (and \emph{a fortiori} the ion
concentrations and electrical potentials) satisfy:
\begin{subequations} \label{eq_phi_rR}
\begin{align}
    \label{eq:phi:r}
    \gamma_{rR} C_{rR} \phi_{rR} &= z^0 F a_r + \alpha_r F \sum_k z^k [k]_r,
    \quad r = 1,\ldots, R-1, \\
    \label{eq:phi:R}
    - \sum_r^{R-1} \gamma_{rR} C_{rR} \phi_{rR} &= z^0 F a_R + \alpha_R F
    \sum_k z^k [k]_R,
\end{align}
\end{subequations}
where $z^k$ (dimensionless) is the \emph{valence} of ion species $k$, $a_r$
(mol/m$^3$) is the \emph{amount of immobile ions}, $F$ (C/mol) denotes
\emph{Faraday's constant}, and $C_{rR}$ (F/m$^2$) is the \emph{membrane
capacitance} of the membrane between compartment $r$ and the ECS.

In this paper, we assume that the ion flux densities $J_r^k$ can be expressed in
terms of the ion concentrations, the electrical potentials, the volume
fractions and the fluid velocity field as:
\begin{equation}
  \label{eq:fluxJ:r}
  J_r^k  = - D_r^k \nabla[k]_r
    - \frac{D_r^k z^k}{ \psi}[k]_r \nabla \phi_r + \alpha_r u_r[k]_r, 
    \quad r = 1,\ldots, R.
\end{equation}
Here, $D_r^k$ (m$^2$/s) denotes the \emph{effective diffusion coefficient} of
ion species $k$ in the region $r$ and may be a given constant or e.g.~a
spatially varying scalar field dependent of $\alpha_r$. The constant $\psi = R
T F^{-1}$ combines Faraday's constant $F$, the \emph{absolute temperature} $T$
(K), and the \emph{gas constant} $R$ (J/(K mol)). We assume that the
system is isothermal, i.e.~$T$ is constant. The ion flux density, i.e.~the
flow rate of ions per unit area, is thus modelled as the sum of three terms:
(i) the diffusive movement of ions due to ionic gradients $- D_r^k \Grad
[k]_r$, (ii) the ion concentrations that are transported via electrical
potential gradients, i.e.~the ion migration $- D_r^k z^k \psi^{-1} [k]_r \Grad
\phi_r$ where $D_r^k \psi^{-1}$ is the \emph{electrochemical mobility} and
(iii) the convective movement $\alpha_r u_r[k]_r$ of ions.

We now turn to the dynamics of the fluid velocity field $u_r$ and the
\emph{mechanical pressure} $p_r:\Omega \times (0, T] \rightarrow \R$ (Pa).
Let the compartmental fluid velocity $u_r$ (m/s) be expressed as:
\begin{equation}  \label{eq:u:r}
    u_r = -\kappa_r\left(\nabla \tilde{p}_r
    + F \sum_{k} z^k[k]_r \nabla \phi_r\right),
    \quad \tilde{p}_r = p_r - RT\frac{a_r}{\alpha_r}, \quad  r = 1,\ldots, R,
\end{equation}
where $\kappa_r$ (m$^4$/(N s)) is the water permeability in compartment
$r$. The fluid velocity field is thus modelled as the sum of three terms: (i)
the mechanical pressure gradient $- \kappa_r\nabla p_r$, (ii) the oncotic
pressure gradient $\kappa_r\nabla (RT \frac{a_r}{\alpha_r})$, and (iii) the
electrostatic forces $- \kappa_r F \sum_{k} z^k[k]_r \nabla \phi_r$. The
mechanical pressure $p_r$ in compartment $r$ and in the ECS $p_R$ may be
related as:
\begin{equation}
    \label{eq:p:r}
    p_r - p_R = \tau_r(\alpha_r),  \qquad r = 1, \ldots, R-1,
\end{equation}
where the \emph{mechanical tension per unit area of the membrane} $\tau_r =
\tau_r(\alpha_r)$ is to be modelled. A simple relation could be $\tau_r
= S_r(\alpha_r - \alpha_r^0)$, were $S_r$ (Pa/m$^3$) denotes the stiffness of
the membrane between compartment $r$ and the ECS, and $\alpha_r^0$ is the
volume fraction at which the membrane has no mechanical tension.  Furthermore,
we assume that the volume--fraction weighted velocity is divergence free, that
is:
\begin{equation} \label{eq:p:R}
\nabla \cdot \left(\sum_{r=1}^{R} \alpha_r u_r\right) = 0.
\end{equation}
Upon inserting~\eqref{eq:u:r}--\eqref{eq:p:r} into~\eqref{eq:p:R}, we obtain
the following equation for the {extracellular mechanical pressure} $p_R:
\Omega \times (0, T] \rightarrow \R$ (Pa),
\begin{equation} \label{eq:p:N_2}
    \nabla \cdot \left(\sum_{r=1}^{R} \kappa_r \alpha_r \left(- \nabla p_R + \nabla RT\frac{a_r}{\alpha_r}
    - F \sum_{k} z^k[k]_r \nabla \phi_r\right)
    - \sum_{r=1}^{R-1} \kappa_r \alpha_r \nabla \tau_r\right)= 0.
\end{equation}

The combination of~\eqref{eq:alpha:r} and~\eqref{eq:alpha:R:sum}--\eqref{eq_phi_rR},
together with the insertion of~\eqref{eq:fluxJ:r}--\eqref{eq:u:r} and~\eqref{eq:p:N_2}, defines a system of $|R||K| + 2|R| + 1$
equations for the $|R||K| + 2|R| + 1$ unknown scalar fields. Appropriate
initial conditions, boundary conditions, and importantly membrane mechanisms
close the system.

\subsection{Membrane mechanisms}
\label{sec:sub:memmech}
The transmembrane water flux $w_{rR}$ is driven by a combination of mechanical
and osmotic pressure, and can be expressed as:
\begin{equation}
    \label{eq:w:r}
    w_{rR} = \eta_r\left(p_r - p_R + RT \left( \frac{a_R}{\alpha_R} + \sum_k [k]_R
    - \frac{a_r}{\alpha_r} - \sum_k [k]_r \right) \right),
\end{equation}
where $\eta_r$ (m$^4$/(mol s)) denotes the \emph{water permeability}.
The compartmental and extracellular potentials across the membrane are coupled
by defining $\phi_{rR}$ as the jump in the electrical potential across the membrane
between compartment $r$ and the ECS (c.f.~\eqref{eq_phi_rR}),
\begin{equation}
    \label{eq:phi:M}
    \phi_{rR} = \phi_r - \phi_R.
\end{equation}
The transmembrane ion flux $J_{rR}^k$ between compartment $r$ and the ECS of
ion species $k$ is subject to modelling, and will typically take the form:
\begin{equation}
    \label{eq:passive:active}
    J_{rR}^k = a_{rR}^k(\phi_{rR},[k]_r) + p_{rR}^k(\phi_{rR}, [k]_r, s_1, \ldots, s_M).
\end{equation}
Here, $a_{rR}^k$ represents active membrane mechanisms (e.g.~ionic pumps)
and $p_{rR}^k$ denotes passive membrane mechanisms (e.g.~leak or voltage
gated ion channels, cotransporters). The passive membrane mechanisms typically
depend on (unitless) gating variables $s_m = s_m(\phi_{rR})$ for $m \in {1,
\dots, M}$ governed by an ODE system of the form:
\begin{equation} \label{eq:ODEs}
    \frac{\partial s_m}{\partial t} = \alpha_m(1 - s_m) - \beta_m s_m,
\end{equation}
where $\alpha_m = \alpha_m(\phi_{rR})$ (1/s) and
$\beta_m = \beta_m(\phi_{rR})$ (1/s) are rate coefficients\cite{sterratt2011principles}.

\subsection{Boundary conditions}
The boundary conditions will strongly depend on the problem of interest. If not otherwise stated, we assume
that no ion flux or fluid leaves on the boundary $\partial \Omega$, that
is,
\begin{equation} \label{bc_ionflux_u}
J^{k}_r(x,t) = 0 \quad \text{and} \quad u_r(x,t) = 0, \quad \text{on } \partial \Omega,
\end{equation}
for $r=1, \dots, R$.

\subsection{Effective diffusion coefficients}
\label{sec:sub:diffusion}
We model the effective diffusion coefficients $D_r^k$ for each ion species
$k \in K$ by:
\begin{equation} \label{eq:diffusion:r}
D_r^k = \alpha_r \chi_r D^k \quad \text{and} \quad D_R^k = \alpha_R D^k,
\end{equation}
for $r=1, \dots, R-1$, where $D^k$ (m$^2$/s) denotes the diffusion coefficient in water 
and $\chi_r$ (dimensionless) reflects the cellular \emph{gap junction connectivity}.

\subsection{Zero flow limit}
\label{sec:zeroflowlim}
We first consider a simplified version of the mathematical model presented in
Section~\ref{sec:mathmodel}, where the compartmental fluid velocity
$u_r$ for $r = 1, \dots, R$ is assumed to be zero. Thus, the
advective terms in~\eqref{eq:alpha:r} and~\eqref{eq:fluxJ:r} vanish
and~\eqref{eq:p:N_2} is decoupled from the rest of the system. The
remaining equations, \eqref{eq:alpha:r}
and~\eqref{eq:alpha:R:sum}--\eqref{eq_phi_rR} with the insertion of 
of~\eqref{eq:fluxJ:r}, describe the dynamics of the volume fractions $\alpha_r$,
the ion concentrations $[k]_r$ and the electrical potentials $\phi_r$, for $r
\in \{1, \dots, R\}$ and for each $k \in K$.

\section{Numerical schemes for the zero flow limit}
\label{sec:scheme}
Below, we present two numerical schemes for the zero flow limit model based on
the finite element method. 

\subsection{Spatial discretization} \label{sec_spatial_discretization}
To obtain a variational formulation, we
multiply~\eqref{eq:alpha:r}
and~\eqref{eq_conservation_rR}--\eqref{eq_phi_rR} with suitable test
functions, integrate over the domain $\Omega$, integrate terms
with higher order derivatives by parts, and insert the boundary
condition~\eqref{bc_ionflux_u}. Below, we let $\inner{u}{v} =
\int_{\Omega} u v \dx$. Let $S_r \subset L^2(\Omega)$, $V_r^k \subset
H^1(\Omega)$, $V_R^k \subset H^1(\Omega)$, $T_r \subset H^1(\Omega)$,
and $T_R \subset H^1(\Omega)$ for $r = 1, \dots, R-1$ be spaces of
functions for $k = 1, \dots, |K|$.  The resulting system reads: find
$\alpha_r \in S_r$, $[k]_r \in V_r^k$, $\phi_r \in T_r$ ($r = 1,
\dots, R-1$), $[k]_R \in V_R^k$, and $\phi_R \in T_R$ such that:
\begin{subequations} \label{eq:varform:summary}
\begin{align}
    \label{eq:varform:summary:alpha:r}
    \inner{\frac{\partial \alpha_r}{\partial t}}{s_r}
    +  \gamma_{rR} \inner{w_{rR}}{s_r}
    &= 0 , \\
    \label{eq:varform:summary:k:r}
    \inner{ \frac{\partial \alpha_r [k]_r}{\partial t}}{v_r^k}
    - \inner{J_r^k}{\nabla v_r^k}
    +\gamma_{rR}\inner{{J}_{rR}^k}{v_r^k}
    &= 0 ,\\
    \label{eq:varform:summary:k:R}
    \inner{ \frac{ \partial \alpha_R [k]_R}{\partial t}}{v_R^k}
    - \inner{J_R^k}{\nabla v_R^k}
    - \sum_r \gamma_{rR}\inner{{J}_{rR}^k}{v_R^k}
    &= 0, \\
    \label{eq:varform:summary:phi:r}
    \inner{\gamma_{rR} C_{rR} \phi_{rR}}{t_r}
    - \inner{z^0 F a_r}{t_r}
    - \inner{F\alpha_r \sum_k z^k [k]_r}{t_r} &= 0,  \\
    \label{eq:varform:summary:phi:R}
    - \sum_r \inner{\gamma_{rR} C_{rR}\phi_{rR}}{t_R}
    - \inner{z^0 F a_R}{t_R} - \inner{ F\alpha_R \sum_k z^k [k]_R}{
        t_R} &= 0,
\end{align}
\end{subequations}
for all $s_r \in S_r$, $v_r^k \in V_r^k$, $v_R^k \in V_R^k$, $t_r \in T_r$,
$t_R \in T_R$. The compartmental ion flux $J_{r}^{k}$ is given
by~\eqref{eq:fluxJ:r}, the transmembrane water flux $w_{rR}$ is given
by~\eqref{eq:w:r}, while the extracellular volume fraction $\alpha_R$ is given
by~\eqref{eq:alpha:R:sum}. The transmembrane ion fluxes $J^k_{rR}$ will depend
on the membrane mechanisms of interest. Note that in the above formulation, the
potentials $\phi_r$ for $r = 1, \dots, R$ are only determined up to a constant.

Next, we discretize the domain $\Omega$ by a simplicial mesh
$\mathcal{T}_h$, where $h$ denotes the mesh size. We introduce
separate finite element spaces for approximating the unknown fields in
the weak formulation~\eqref{eq:varform:summary}. Specifically, we
approximate the volume fractions $\alpha_r$ using piecewise constants ($P_0$),
and the ion concentrations $[k]_r$, $[k]_R$ and potentials $\phi_r$,
$\phi_R$ using continuous piecewise linear polynomials $P_1^c$ defined
relative to the mesh $\mathcal{T}_h$.

\subsection{Temporal PDE discretization: BDF2 and CN}
\label{sec:sub:scheme:BDF2:CN}

We apply a second order Backward Differentiation Formula (BDF2) in time,
resulting in the following system: given
$\alpha_r^n$, $[k]_r^n$, and $[k]_R^n$ at time level $n$ and $\alpha_r^{n-1}$,
$[k]_r^{n-1}$, and $[k]_R^{n-1}$ at time level $n-1$, find at time level $n+1$
the volume fractions $\alpha_r
\in S_r$, the ion concentrations $[k]_r \in V_r^k$, $[k]_R \in V_R^k$ and the
potentials $\phi_r \in T_r$, $\phi_R \in T_R$ such that:
\begin{subequations} \label{eq:BDF2:summary}
\begin{align}
    \label{eq:BDF2:summary:alpha:r}
    \frac{1}{\Delta t}\inner{\alpha_r - \frac{4}{3}\alpha_r^n + \frac{1}{3}\alpha_r^{n-1}}{s_r}
    +  \frac{2}{3}\gamma_{rR} \inner{w_{rR}}{s_r}
    &= 0 , \\
    \label{eq:BDF2:summary:k:r}
    \begin{split}
     \frac{1}{\Delta t} \inner{\alpha_r [k]_r - \frac{4}{3}\alpha_r^n [k]_r^n 
     + \frac{1}{3}\alpha_r^{n-1} [k]_r^{n-1}}{v_r^k}
    - \frac{2}{3}\inner{J_r^{k,n}}{\nabla v_r^k}&\\
    + \frac{2}{3}\gamma_{rR}\inner{{J}_{rR}^{k,n}}{v_r^k}
    &= 0 ,
    \end{split}\\
    \label{eq:BDF2:summary:k:R}
    \begin{split}
    \frac{1}{\Delta t} \inner{\alpha_R [k]_R - \frac{4}{3}\alpha_R^n [k]_R^n
    + \frac{1}{3}\alpha_R^{n-1} [k]_R^{n-1} }{v_R^k}
    - \frac{2}{3}\inner{J_R^{k,n}}{\nabla v_R^k}&\\
    - \frac{2}{3}\sum_r \gamma_{rR}\inner{{J}_{rR}^{k,n}}{v_R^k}
    &= 0, 
    \end{split}\\
    \label{eq:BDF2:summary:phi:r}
    \inner{\gamma_{rR} C_{rR} \phi_{rR}}{t_r}
    - \inner{z^0 F a_r}{t_r}
    - \inner{F\alpha_r \sum_k z^k [k]_r}{ t_r} &= 0,  \\
    \label{eq:BDF2:summary:phi:R}
    - \sum_r \inner{\gamma_{rR} C_{rR}\phi_{rR}}{ t_r}
    - \inner{z^0 F a_R}{t_R} - \inner{ F\alpha_R \sum_k z^k [k]_R}{
        t_R} &= 0,
\end{align}
\end{subequations}
for all $s_r \in S_r$, $v_r^k \in V_r^k$, $v_R^k \in V_R^k$, $t_r \in T_r$,
$t_R \in T_R$. The solutions at time level
$1$ and $0$ are given by a backward Euler step and the initial conditions, respectively. 
Note that the passive part of the membrane flux ${J}_{rR}^{k,n}$
is treated implicitly, whereas the active part is treated explicitly, both in
the BDF2 and the BE time stepping schemes (see~\eqref{eq:passive:active}).

We will also compare with the following Crank-Nicholson (CN) scheme in time: given
$\alpha_r^n$, $[k]_r^n$, and $[k]_R^n$ at time level $n$, find at time
level $n+1$ the volume fractions $\alpha_r \in S_r$, the ion
concentrations $[k]_r \in V_r^k$, $[k]_R \in V_R^k$, and the
potentials $\phi_r \in T_r$, $\phi_R \in T_R$, such that
\begin{subequations}
\begin{align}
    \label{eq:CN:summary:alpha:r}
    \frac{1}{\Delta t}\inner{  \alpha_r - \alpha_r^n}{s_r}
    +  \gamma_{rR} \inner{w_{rR}^{n+1/2}}{s_r}
    &= 0 , \\
    \label{eq:CN:summary:k:r}
     \frac{1}{\Delta t} \inner{ \alpha_r [k]_r - \alpha_r^n [k]_r^n}{v_r^k}
    - \inner{J_r^{k,n+1/2}}{\nabla v_r^k}
    +\gamma_{rR}\inner{{J}_{rR}^{k, n+1/2}}{v_r^k}
    &= 0 ,\\
    \label{eq:CN:summary:k:R}
    \frac{1}{\Delta t} \inner{\alpha_R [k]_R - \alpha_R^n [k]_R^n}{v_R^k}
    - \inner{J_R^{k, n+1/2}}{\nabla v_R^k}
    - \sum_r \gamma_{rR}\inner{{J}_{rR}^{k, n+1/2}}{v_R^k}
    &= 0, 
\end{align}
\end{subequations}
together with~\eqref{eq:BDF2:summary:phi:r} and~\eqref{eq:BDF2:summary:phi:R},
for all $s_r \in S_r$, $v_r^k \in V_r^k$, $v_R^k \in V_R^k$, $t_r \in T_r$,
$t_R \in T_R$ Here, $f^{n+1/2}$ denotes the solution at time level $n+1/2$ and
is approximated by $(f^n + f^{n+1})/2$, for $f \in \{w_{rR}, J^k_{rR},
J^k_r\}$.

\subsection{Strang splitting scheme}
\label{sec:sub:scheme:strang}

We use a second order Strang splitting scheme where we solve the coupled
systems of ODEs and PDEs step-wise in the following manner:
\begin{enumerate}
    \item Insert the previous solution of $\phi_{rR}$, $[k]_r$, $[k]_R$
        into the system of ODEs~\eqref{eq:ODEs}, and solve ODEs for a half timestep
        $\Delta t/2$;
    \item Solve the system of PDEs
        (e.g.~\eqref{eq:BDF2:summary})
        with values for the gating variables $s_1$, $s_2$, \dots,
        $s_M$  from step (1) for one time step $\Delta t$;
    \item Insert solution for $\phi_{rR}$, $[k]_r$, $[k]_R$ from step (2)
        into the system of ODEs~\eqref{eq:ODEs}, and solve the ODEs for a half timestep
        $\Delta t/2$;
\end{enumerate}
and continue steps (1)--(3) until the global end time is reached. We
compare the Strang splitting scheme with a first order Godunov scheme;
see e.g.~Sundnes et al.\cite{sundnes2006solving} for a detailed
description the Godunov method..

\subsection{ODE solvers}
\label{sec:sub:scheme:ODE}
We consider three different schemes for the ODE time-stepping: a fourth order
Runge-Kutta (RK4) method, a fourth order explicit singly diagonal implicit
Runge-Kutta (ESDIRK4) method or a second order backward Euler (BE) method.  The
ODEs are solved with the RK4 method unless otherwise stated. See
e.g.~Langtangen and Linge\cite{petter2017finite} for details.

\section{Numerical convergence study: smooth analytical solution}
\label{sec:MMS}
To evaluate the numerical accuracy of the various schemes presented
above, we begin by constructing a smooth analytical solution using the method of
manufactured solutions\cite{roache1998verification}.

\subsection{Problem description}
\label{sec:sub:MMS:model}
We consider a two-compartment version of the model in the zero flow limit and a
neuronal ($n$, $r=1$) and an extracellular ($e$, $r=2$) compartment. We use
$1,2$ and $n,e$ interchangeably for subscripts of our variables and model
parameters.  In each compartment, we model the movement of potassium (K$^+$),
sodium (Na$^+$) and chloride (Cl$^-$).  In the numerical experiments of this
test, we consider a one dimensional domain $\Omega = [0,1]$ uniformly meshed
with $N \in \{8,16,32,64,128\}$ elements. We initially set $\Delta t = 10^{-3}$
s, and then halve the timestep with each spatial refinement. The errors are
evaluated at $t = 2 \times 10^{-3}$ s. Further, we assume that the
transmembrane ion flux $J_{ne}^k$ depends on the gating variables $m$, $h$ and
$g$ governed by the following ODEs:
\begin{equation}
    \label{eq:MMS:gating}
    \frac{\partial m}{\partial t} = \phi_{ne}, \quad
    \frac{\partial h}{\partial t} = \phi_{ne}, \quad
    \frac{\partial g}{\partial t} = \phi_{ne},
\end{equation}
where $\phi_{ne} = \phi_n - \phi_e$ is the membrane potential. The analytical
solution to the PDE system is given by:
\begin{equation}
    \label{eq:MMS:compartmental}
\begin{aligned}
    \alpha_n &= 0.3 - 0.1\sin(2\pi x)\exp(-t), && \\
    [\Na^+]_n &= 0.7 + 0.3\sin(\pi x)\exp(-t), & [\Na^+]_e &= 1.0 + 0.6\sin(\pi x)\exp(-t), \\
    [\K^+]_n &= 0.3 + 0.3\sin(\pi x)\exp(-t), & [\K^+]_e &= 1.0 + 0.2\sin(\pi x)\exp(-t), \\
    [\Cl^-]_n &= 1.0 + 0.6\sin(\pi x)\exp(-t), & [\Cl^-]_e &= 2.0 + 0.8\sin(\pi x)\exp(-t), \\
    \phi_n &= \sin(2\pi x)\exp(-t), & \phi_e &= \sin(2\pi x)(1 + \exp(-t)),
\end{aligned}
\end{equation}
and the solution to the system of ODEs~\eqref{eq:MMS:gating} is
\begin{equation} \label{eq:MMS:compartmental_2}
m = \cos(t) \cos(\pi x), \quad h = \cos(t) \cos(\pi x), \quad g = \cos(t) \cos(\pi x).
\end{equation}
Parameter values are given in Table~\ref{tab:params:n}. Initial and boundary
conditions are governed by the exact solutions~\eqref{eq:MMS:compartmental}.

\begin{table}[ht]
\begin{center}
    \begin{tabular}{  lcccc }
        \toprule
        Parameter & Symbol &  Value & Unit & Ref. \\
        \midrule
        Temperature &$T$ & $310$ & K & --\\
        Faraday's constant &$F$ & $96485 $ & C/mol & --  \\
        Gas constant &$R$ & $8.3144598$ & J/(mol K) & -- \\
        Membr. area-to-volume neuron &$\gamma_{ne}$ &
        $6.3849\times 10^{5}$ & 1/m &  \cite{kager2000simulated} \\
        Membr. capacitance neuron &$C_{ne}$ & $7.5 \times 10^{-3}$ & F/m$^2$  &
        \cite{kager2000simulated} \\
        Membr. water permeability neuron &$\eta_{ne}$ &  $5.4\times 10^{-10}$
        & m$^4$/(mol s) & \cite{o2016effects} \\
        Gap junction connectivity neuron & $\chi_n$ & $0$ &  & \cite{mori2015multidomain} \\
        Diffusion coefficient $\Na^+$ &$D^{\Na}$ & $1.33\times 10^{-9}$ & m$^2$/s &  \cite{hille2001} \\
        Diffusion coefficient $\K^+$ &$D^{\K}$ & $1.96\times 10^{-9}$ & m$^2$/s &  \cite{hille2001} \\
        Diffusion coefficient $\Cl^+$ &$D^{\Cl}$ & $2.03\times 10^{-9}$ & m$^2$/s &  \cite{hille2001} \\
        Valence $\Na^+$ &$z^{\Na}$ & $1$ &  & -- \\
        Valence $\K^+$ &$z^{\K}$ & $1$ &  &   -- \\
        Valence $\Cl^-$ &$z^{\Cl}$ & $-1$ &  & -- \\
        \bottomrule
\end{tabular}
\end{center}
    \caption{Physical model parameters. We use SI base units, that is, Kelvin
    (K), Coulomb (C), mole (mol), meter (m), second (s), and Joule (J).  The
    values are collected from Hille et al.\protect\cite{hille2001}, Kager et
    al.\protect\cite{kager2000simulated},
    Mori\protect\cite{mori2015multidomain}, and O'Connell et
    al.\protect\cite{o2016effects}. -- indicates that a standard value is
    used.}
    \label{tab:params:n}
\end{table}

\subsection{Convergence and convergence rates under refinement}
Based on the approximation spaces and the time discretization, we expect 
the optimal theoretical rate of convergence to be $1$ in the $H^1$-norm and $2$ in
the $L^2$-norm for the concentrations $[k]_n$, $[k]_e$, the potentials
$\phi_n$, $\phi_e$ and the gating variables $m,h,g$, and $1$ in the $L^2$-norm
for the volume fraction $\alpha_n$. Our numerical observations are in agreement
with these optimal rates, both for the BDF2 scheme
(Table~\ref{tab:MMS:BDF2}A) and for the CN scheme (Table~\ref{tab:MMS:BDF2}B). We
observe second order convergence in the $L^2$-norm for the approximation of the
extracellular and intracellular concentrations and potentials, and first order
convergence in the $H^1$-norm. The neuron potential approximated by the CN
scheme displays superconvergence between $N=64$ and $N=128$
(Table~\ref{tab:MMS:BDF2}B). For the volume fractions, we observe a convergence
rate of $1$ in the $L^2$-norm.

\begin{table}[ht]
    \textbf{A} 
    \begin{center}
    \begin{tabular}{ c  c  c  c  c}
        \toprule
        $N$
        & $\norm{[{\rm Na}]_e - {[{\rm Na}]_e}_h}_{L2}$
        & $\norm{\phi_n - {\phi_n}_h}_{L2}$
        & $\norm{\alpha_n - {\alpha_n}_h}_{L^2}$
        & $\norm{m - {m}_h}_{L^2}$ \\
        \midrule
        8 & 5.74E-03(-----) & 7.66E-02(-----) & 1.58E-02(-----) & 9.99E-03(-----)\\
        16 & 1.45E-03(1.99) & 1.90E-02(2.01) & 7.97E-03(0.98) & 2.50E-03(2.00)\\
        32 & 3.63E-04(1.99) & 4.72E-03(2.01) & 4.00E-03(1.00) & 6.26E-04(2.00)\\
        64 & 9.11E-05(2.00) & 1.17E-03(2.01) & 2.00E-03(1.00) & 1.56E-04(2.00)\\
        128 & 2.28E-05(2.00) & 2.93E-04(2.00) & 1.00E-03(1.00) & 3.91E-05(2.00)\\
        \midrule

        $N$
        & $\norm{[{\rm Na}]_e - {[{\rm Na}]_e}_h}_{H^1}$
        & $\norm{\phi_n - {\phi_n}_h}_{H^1}$ \\
        \cmidrule{1-3}
        8 & 1.50E-01(-----) & 1.02E+00(-----)\\
        16 & 7.53E-02(1.00) & 5.05E-01(1.02)\\
        32 & 3.77E-02(1.00) & 2.52E-01(1.00)\\
        64 & 1.88E-02(1.00) & 1.26E-01(1.00)\\
        128 & 9.43E-03(1.00) & 6.28E-02(1.00)\\
        \cmidrule{1-3}
    \end{tabular}
    \end{center}
    \textbf{B}
    \begin{center}
    \begin{tabular}{ c  c  c  c  c}
        \toprule
        $N$
        & $\norm{[{\rm Na}]_e - {[{\rm Na}]_e}_h}_{L^2}$
        & $\norm{\phi_n - {\phi_n}_h}_{L^2}$
        & $\norm{\alpha_n - {\alpha_n}_h}_{L^2}$
        & $\norm{m - {m}_h}_{L^2}$ \\
        \midrule
        8 & 2.44E-03(-----) & 1.32E-01(-----) & 1.52E-02(-----) & 4.27E-03(-----)\\
        16 & 6.05E-04(2.01) & 3.43E-02(1.94) & 7.78E-03(0.97) & 1.04E-03(2.03)\\
        32 & 1.52E-04(2.00) & 8.68E-03(1.98) & 3.91E-03(0.99) & 2.59E-04(2.01)\\
        64 & 3.83E-05(1.99) & 2.17E-03(2.00) & 1.96E-03(1.00) & 6.46E-05(2.00)\\
        128 & 9.83E-06(1.96) & 6.54E-05(5.05) & 9.82E-04(1.00) & 1.61E-05(2.00)\\
        \midrule
        $N$
        & $\norm{[{\rm Na}]_e - {[{\rm Na}]_e}_h}_{H^1}$
        & $\norm{\phi_n - {\phi_n}_h}_{H^1}$ \\
        \cmidrule{1-3}
        8 & 1.47E-01(-----) & 1.17E+00(-----)\\
        16 & 7.38E-02(1.00) & 5.21E-01(1.17)\\
        32 & 3.70E-02(1.00) & 2.50E-01(1.06)\\
        64 & 1.85E-02(1.00) & 1.24E-01(1.02)\\
        128 & 9.30E-03(0.99) & 6.17E-02(1.00)\\
        \cmidrule{1-3}
    \end{tabular}
    \end{center}
  \caption{Selected $L^2$ (upper panel) and $H^1$-errors (lower panel) and
    convergence rates (in parenthesis) for the BDF2 (\textbf{A}) and CN (\textbf{B}) schemes at time $t=0.002$ s.
    The test was run on the unit interval, and we initially let $\Delta t = 1 \times 10^{-3}$ s,
    and then halve the timestep with each mesh refinement. The spatial
    discretization consists of $N$ intervals.} 
    \label{tab:MMS:BDF2}
\end{table}

\section{Numerical convergence study: physiological CSD wave}
\label{sec:convergence}

Next, we consider the simulation of cortical spreading depression with
a sharp wave front. This is a more challenging problem, with
characteristics quite different from the previous smooth MMS case. We
numerically study the effect of splitting scheme, time-discretization
of the PDEs, and discretization of the ODEs.

\subsection{Problem description}
\label{sec:sub:convergence:model}

\begin{table}[ht]
\begin{center}
    \begin{tabular}{  l  c  c  c  c }
        \toprule
        Parameter & Symbol & Value & Unit & Ref. \\
        \midrule
        $\Na^+$ leak conductance neuron&$g_{n,\T{leak}}^{\Na}$ & $0.2$ &
        S/m$^{2}$ &\cite{kager2000simulated}\\
        $\K^+$ leak conductance neuron&$g_{n,\T{leak}}^{\K}$ & $0.7$ &
        S/m$^{2}$ &\cite{kager2000simulated}\\
        $\Cl^-$ leak conductance neuron &$g_{n,\T{leak}}^{\Cl}$ & $2.0$ & S/m$^{2}$ &\cite{o2016effects}\\
        Maximum pump rate neuron &$\hat{I}_n$ &$0.1372$ & A/m$^{2}$ &\cite{yao2011}\\
        Threshold for pump $[\Na^+]_r$ &$m_{\Na}$ & $7.7$ & mol/m$^{3}$ &\cite{yao2011}\\
        Threshold for pump $[\K^+]_R$ &$m_{\K}$ & $2.0$ & mol/m$^{3}$ &\cite{yao2011}\\
        \bottomrule
\end{tabular}
 \caption{Physical parameters for the neuron membrane mechanisms. We use SI
    base units, i.e.~meter (m), mole (mol), Siemens (S) and ampere (A). The
    values are collected from Kager et al.\protect\cite{kager2000simulated}, O'Connell
    et al.\protect\cite{o2016effects}, and Yao et al.\protect\cite{yao2011}.}
    \label{tab:params:mem:n}
\end{center}
\end{table}

We define a more physiological version of the mathematical
model considered in Section~\ref{sec:sub:MMS:model} (two compartments, zero
flow limit, and a neuronal ($n$, $r=1$) and an extracellular compartment ($e$,
$r=2$)). In each compartment, we again model the movement of potassium (K$^+$),
sodium (Na$^+$) and chloride (Cl$^-$). We consider a 1D domain of length $0.01$
m ($10$ mm), and now apply physiologically relevant neuronal membrane
mechanisms as described below, notably including a system of ODEs describing
the gating variables of voltage-gated sodium and potassium channels. The domain
and the transmembrane ion flux densities are taken from the original Mori study\cite{mori2015multidomain}.

Concretely, the transmembrane ion flux
densities $J_{ne}^k$ for $k = \{\Na^+, \K^+, \Cl^- \}$ in~\eqref{eq_conservation_rR} are modelled as:
\begin{subequations}
    \label{eq:totmem:fluxes}
\begin{align}
    J^{\Na}_{ne} &= \frac{1}{F z^{\Na}}\left( I^{\Na}_{n,\T{leak}} + I_{\NaP} +
    3I_{n,\ATP} + I^{\Na}_{\rm ex} \right), \\
    J^{\K}_{ne} &=  \frac{1}{F z^{\K}}\left(I^{\K}_{n,\T{leak}} + I_{\KDR} + I_{\KA}  -
    2I_{n,\ATP} + I^{\K}_{\rm ex} \right), \\
    J^{\Cl}_{ne} &= \frac{1}{F z^{\Cl}}\left(I^{\Cl}_{n,\T{leak}} + I^{\Cl}_{\rm ex}\right),
\end{align}
\end{subequations}
where $I_{n,\T{leak}}^{\Na}$, $I_{n,\T{leak}}^{\K}$, $I_{n,\T{leak}}^{\Cl}$,
$I_{\NaP}$, $I_{\KDR}$, $I_{\KA}$, and $I_{n,\ATP}$ denote the sodium leak current, the
potassium leak current, the chloride leak current, the persistent sodium
current, the potassium delayed rectifier current,
the transient potassium current,
and the Na/K/ATPase current, respectively. Further, $I^{\Na}_{\rm ex}$,
$I^{\K}_{\rm ex}$, and $I^{\Cl}_{\rm ex}$ are excitatory fluxes used to trigger
a cortical spreading depression wave and whose expressions are given by~\eqref{eq:trigger} in Section~\ref{sec_init_CSD_wave} below.  Note
that the currents (A/m$^{2}$) are converted to ion fluxes (mol/(m$^{2}$s)) by
dividing by Faraday's constant $F$ times the valence $z^k$.

The leak currents (A/m$^2$) of ion species $k$
over the membrane between compartment $r$ (here with $r = n$) and $R$
are modelled as:
\begin{equation}
    \label{eq:mem:leak}
    I_{r,\T{leak}}^k = g_{r,\T{leak}}^k (\phi_{rR} - E^k_r), \qquad
    E^{k}_r = \frac{RT}{F z^k} \ln\left(\frac{[k]_R}{[k]_r}\right),
\end{equation}
where $E_k^r$ (V) denotes the Nernst potential. The values for the neuronal
leak conductances $g_{n,\T{leak}}^k$ are listed in Table~\ref{tab:params:mem:n}.

The current-voltage relation for the voltage-gated currents $I_{\NaP}$, $I_{\KDR}$ and $I_{\KDR}$
(A/m$^{2}$) are described by the Goldman-Hodgkin-Katz (GHK) current equation:
\[I_{\GHK}^k = g^k  m^p h^q
    \frac{F \mu \left([k]_e - [k]_n e^{-\mu}\right)}
    {1 - e^{-\mu}}.\]
Here, $\mu = F \phi_{ne}/(RT)$ is dimensionless, while $g^k$ (m/s)
denotes the product of membrane permeability and conductance.
The gating variables $m$ and $h$ describe the proportion of open
voltage--gated ion channels, and are governed by the following ODEs:
\begin{subequations}
\begin{align}
    \label{ODE:gating}
    \frac{\partial m}{\partial t} &= \alpha_m(\phi_{ne})(1 - m) - \beta_m(\phi_{ne})m , \\
    \frac{\partial h}{\partial t} &= \alpha_h(\phi_{ne})(1 - h) - \beta_h(\phi_{ne})h ,
\end{align}
\end{subequations}
where the activation rate functions $\alpha_m : \R \rightarrow \R$ and
$\beta_m: \R \rightarrow \R$, and the inactivation rate functions $\alpha_h :
\R \rightarrow \R$ and $\beta_h: \R \rightarrow \R$, are specified for each
current in Table~\ref{tab:params:currents:mem}. The initial conditions are
given in Table~\ref{tab:init}A (Supplementary Tables).

\begin{table}[ht]
\begin{center}
    \begin{tabular}{cccl}
        \toprule
        Current & Permeability & Gates & Voltage dependent rate constants \\
        (A/m$^{2}$) & (m/s) &  &  \\
        \midrule
        $I_{\NaP}$ & $2.0\times 10^{-7}$ & m$^2$ h &
        $\alpha_m = 1/(6 + 6\exp(-0.143 \phi_{ne} - 5.67))$ \\
        & & & $\beta_m = 1/6 - \alpha_m$ \\
        & & & $\alpha_h = 5.12\times 10^{-8} \exp(-0.056\phi_{ne} - 2.94)$ \\
        & & & $\beta_h = 1.6 \times 10^{-6}/(1 + \exp(-0.2 \phi_{ne} - 8))$ \\
        \midrule
        $I_{\KDR}$ & $1.0\times 10^{-5}$  & m & $\alpha_m = 0.016\frac{- \phi_{ne}
        - 34.9} {\exp(-0.2 \phi_{ne} - 6.98)- 1}$ \\
        & & & $\beta_m = 0.25 \exp(-0.025\phi_{ne} - 1.25)$ \\
        \midrule
        $I_{\KA}$  & $2.0\times 10^{-6}$ & m$^2$ h & $\alpha_m = 0.02\frac{-
        \phi_{ne} - 56.9} {\exp(-0.1 \phi_{ne} - 56.9) - 1}$ \\
        & & & $\beta_m = 0.0175 \frac{\phi_{ne} + 29.9} {\exp(0.1\phi_{ne} + 29.9) - 1}$ \\
        & & & $\alpha_h = 0.016\exp(-0.05\phi_{ne} - 4.61)$ \\
        & & & $\beta_h = 0.5/(\exp(-0.2\phi_{ne} - 11.98) + 1)$ \\
        \bottomrule
\end{tabular}
\end{center}
 \caption{Permeability, gates and voltage dependent expressions for the
    activation rates ($\alpha$) and the inactivation rates ($\beta$) for the
    persistent sodium current ($I_\NaP$), the potassium delayed rectifier
    current ($I_\KDR$) and the transient potassium current ($I_\KA$) in the
    neuron membrane. The values are collected from Kager et
    al.\protect\cite{kager2000simulated} and Yao et al.\protect\cite{yao2011}.}
    \label{tab:params:currents:mem}
\end{table}

The Na/K/ATPase pump exchanges 2 potassium ions for 3 sodium ions, and
the pump current $I_{r,\ATP}$ (A/m$^{2}$) over the membrane between compartment
$r$ (here with $r=n$) and $R$ is modelled as:
\begin{equation}
    \label{eq:mem:ATP}
    I_{r,\ATP} = \frac{\hat{I}_r}{(1 + \frac{m_{\K}}{[\K^+]_R})^2
    (1 + \frac{m_{\Na}}{[\Na^+]_r})^3},
\end{equation}
where $\hat{I}_r$ is the maximum pump rate, and $m_{\Na}$ and $m_{\K}$ denote
the sodium and potassium pump threshold, respectively. The values for the
neuron pump parameters are listed in Table~\ref{tab:params:mem:n}.

\subsubsection{Initiation of the CSD wave} \label{sec_init_CSD_wave}

Following the original Mori study\cite{mori2015multidomain}, we initiate a CSD wave by adding
excitatory fluxes to the transmembrane fluxes defined by~\eqref{eq:totmem:fluxes} in the following manner:
\begin{equation} 
    \label{eq:trigger}
    \begin{aligned}
    I_{\rm ex}^k &= G_{\rm ex}(E^k_r - \phi_{nR}), \\
    G_{\rm ex}(x, t) &= 
    \begin{cases*}
    G_{\max}\cos^2(\pi x / 2L_{\rm ex}) \sin(\pi t/T_{\rm ex}) & if $ x \leq L_{\rm ex}$ and $t \leq T_{\rm ex}$,\\
    0 & otherwise,
    \end{cases*}
    \end{aligned}
\end{equation}
for $k = \{\Na^+, \K^+, \Cl^-\}$. We set $L_{\rm ex} = 2.0 \times
10^{-5}$ m, $T_{\rm ex} = 2$ s, and $G_{\max} = 5.0$ S/m$^2$. 

\subsection{CSD wave characteristics} \label{sec_CSD_wave_charac}

\begin{figure}
\centering
  \includegraphics[width=\linewidth]{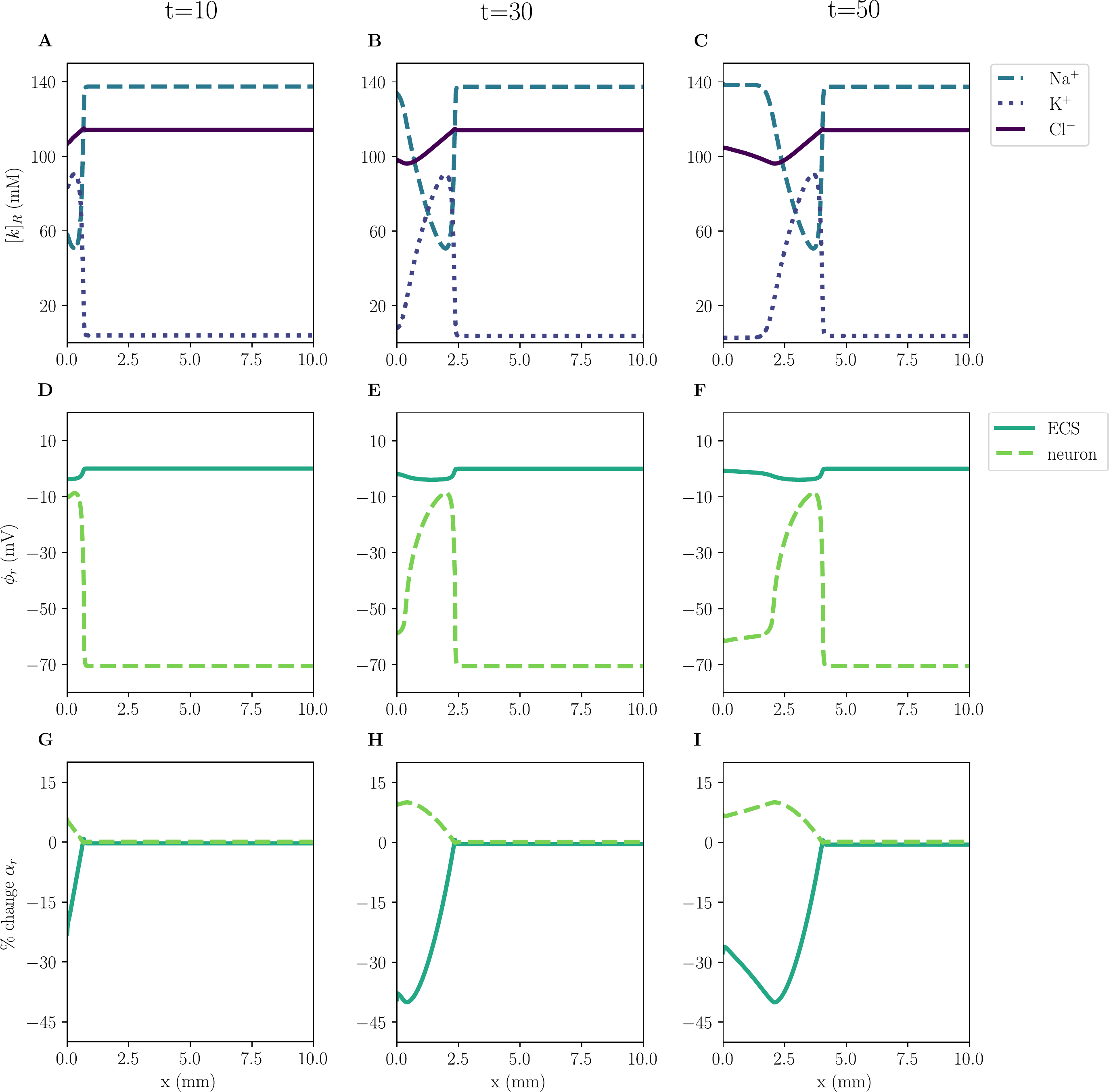}
  \caption{Snapshot of a CSD wave spreading through the tissue domain --
    extracellular sodium (Na$^+$), potassium (K$^+$) and chloride (Cl$^-$)
    concentrations (\textbf{A}, \textbf{B}, \textbf{C}), neuron and extracellular
    potentials (\textbf{D}, \textbf{E}, \textbf{F}), and neuron and extracellular
    changes in volume fractions (\textbf{G}, \textbf{H}, \textbf{I}) -- at, from
    left to right, $t = 10, 30, 50$ s. Numerical scheme and resolutions: BDF2,
    Strang, RK4 with $N=32000$ and $\Delta t = 0.195$ s.}
  \label{fig:zerolimit:wave}
\end{figure}

The excitatory flux stimulation leads to a wave of neuronal depolarization,
ionic concentration changes and neuronal swelling spreading through the tissue
domain (Figure~\ref{fig:zerolimit:wave}). We observe a dramatic depolarization
of the neuron potential: from $-71$ mV to $-8.7$ mV, accompanied by a small drop
in the extracellular potential: from $0$ to $-3.9$ mV. This latter drop is known
as the DC shift, a key characteristic associated with cortical spreading
depression\cite{pietrobon2014chaos}. The neuronal depolarization wave is
followed by a complete break-down of the ionic homeostatis: a substantial
increase in the concentrations of extracellular potassium, and decreases in
extracellular sodium and chloride. In response to the ionic shifts, the neurons
swell with an increase in volume fraction of up to 10\%, while the
extracellular space shrinks correspondingly. We observe that although we have not
enforced positivity of $\alpha_r$ explicitly, $\alpha_r \geq 0$ holds throughout our numerical experiments.

Experimental studies typically focus on the speed and duration of the
CSD wave. We thus define the following quantities of interest to be
studied further quantitatively in terms of numerical convergence.
\begin{itemize}[leftmargin=0cm,itemindent=.5cm,noitemsep]
\item
    The point of $x_{\T{peak}, 50}$ of peak neuron potential $\phi_n$ at
    $t = 50$ s.
\item
    The mean CSD wave speed $\bar{v}_{\rm CSD}$ is computed by the
    distance between the peak neuronal potential at two points (as
    long as $\phi_n > -20$ mV) divided by the time elapsed to cover
    that distance.
\item
    The (temporal) duration $d_{\rm CSD}$ of the CSD wave in terms of elevated extracellular
    potassium levels at $x = 1$ mm, where
    \begin{align*}
        \mathbb{T} & =\{t \ | \ [\K^+]_R(1, t) > k_{\rm thres}\} \\
        d_{\rm CSD} &= \max \mathbb{T} - \min\mathbb{T}
    \end{align*}
    with $k_{\rm thres} = 10$ mM for $t \in (0, T]$.
\item
    The (spatial) width $w_{\rm CSD}$ of the extracellular potassium wave at
    $t= 50$ s, where
    \begin{align*}
        \mathbb{X} &= \{x \ | \ [\K^+]_R(x, 50) > k_{\rm thres}\} \\
        w_{\rm CSD} &= \max\mathbb{X} - \min\mathbb{X}
    \end{align*}
    with $k_{\rm thres} = 10$ mM for $x \in \Omega$.
\end{itemize}

\subsection{Convergence of the CSD wave characteristics}
\label{sec:sub:convergence:res}
\begin{figure}[pb]
    \textbf{A} 
    \vspace{0.25cm}
    \begin{center}
    \begin{overpic}[width=0.12\linewidth]{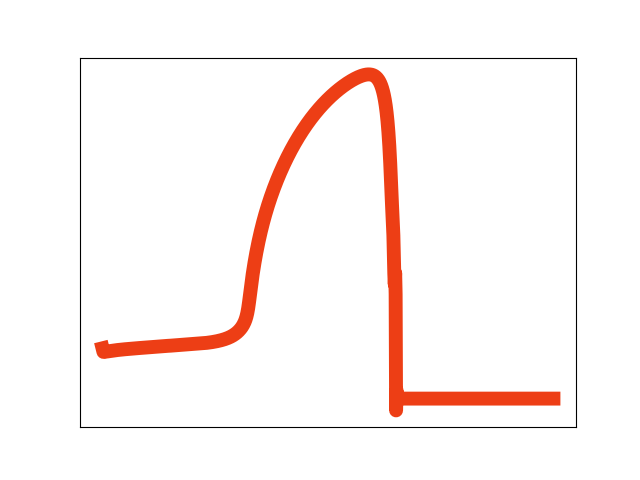}
        \put(-3,82){\color{black}\vector(1,0){50}}
        \put(55,82){Temporal refinement}
        \put(-21,-155){\rotatebox{90}{Spatial refinement}}
        \put(-3,82){\color{black}\vector(0,-1){50}}
    \end{overpic}
    \includegraphics[width=0.12\linewidth]{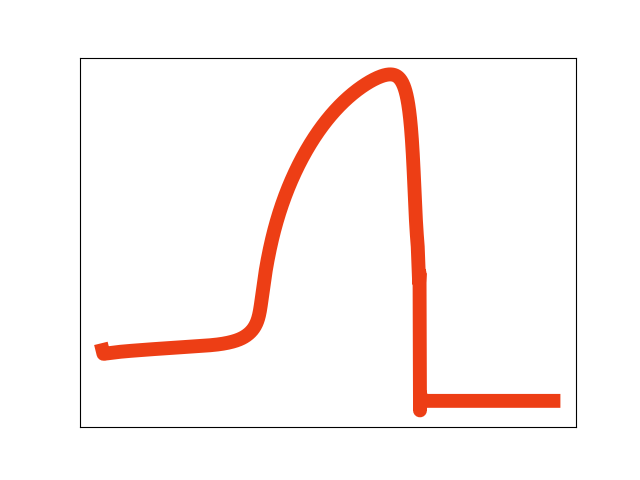}
    \includegraphics[width=0.12\linewidth]{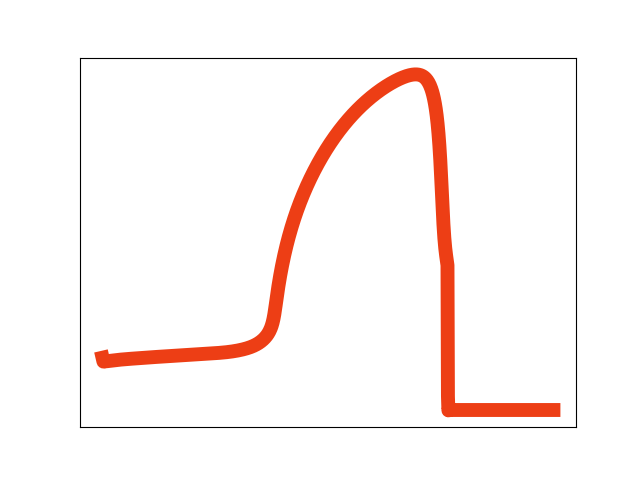}
    \includegraphics[width=0.12\linewidth]{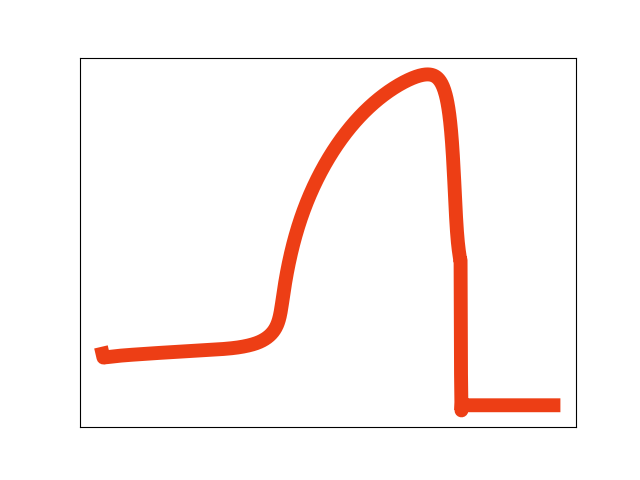}
    \includegraphics[width=0.12\linewidth]{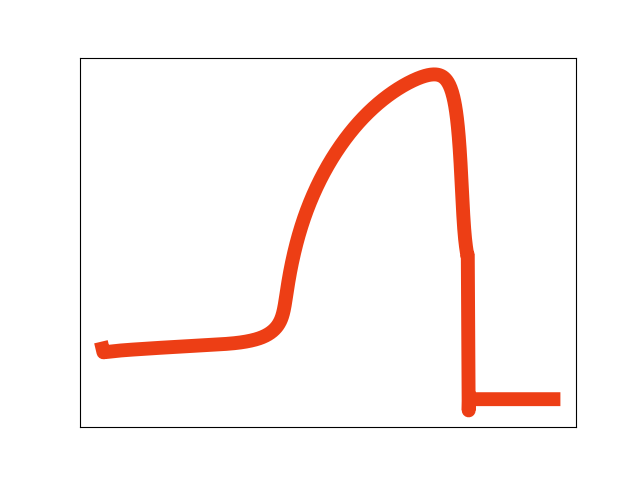}
    \includegraphics[width=0.12\linewidth]{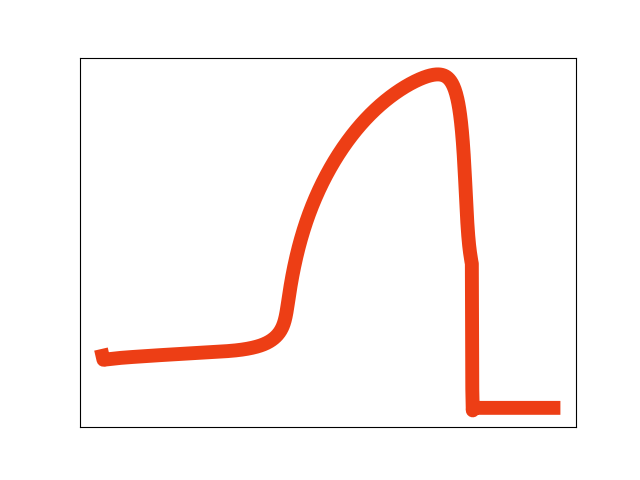}
    \includegraphics[width=0.12\linewidth]{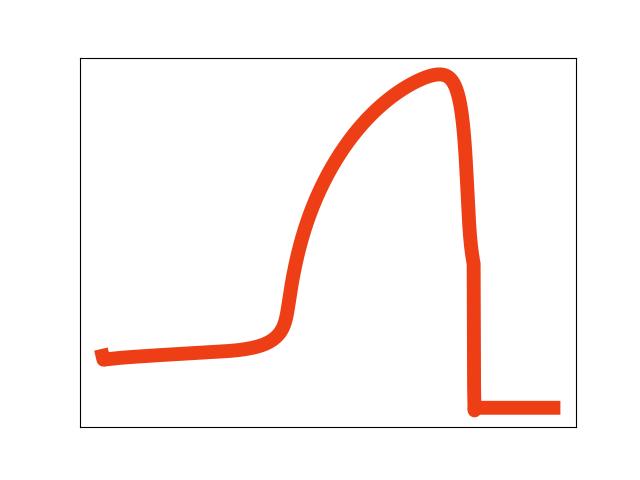}
    \\
    \includegraphics[width=0.12\linewidth]{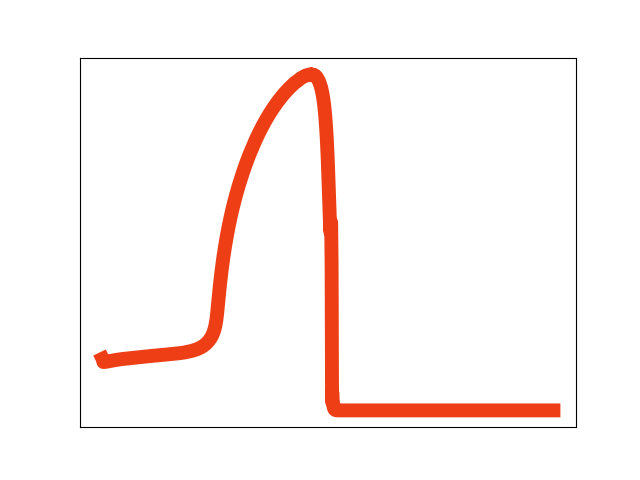}
    \includegraphics[width=0.12\linewidth]{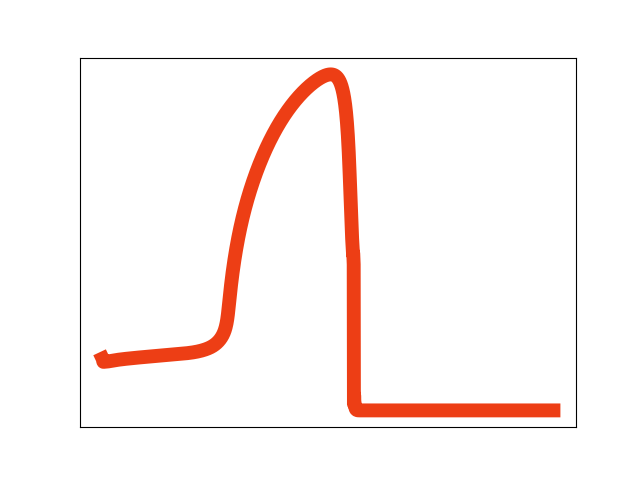}
    \includegraphics[width=0.12\linewidth]{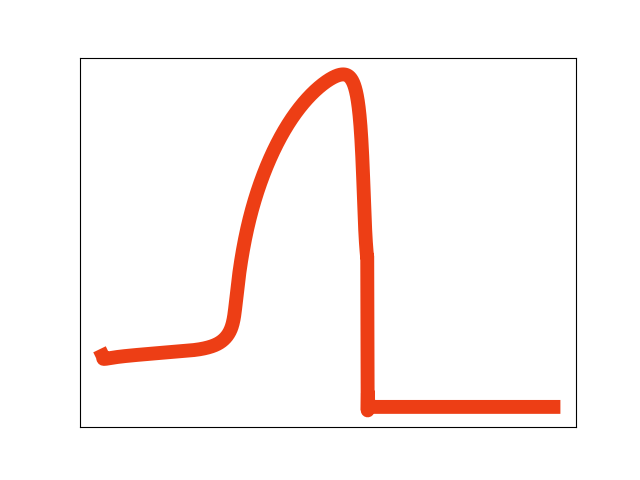}
    \includegraphics[width=0.12\linewidth]{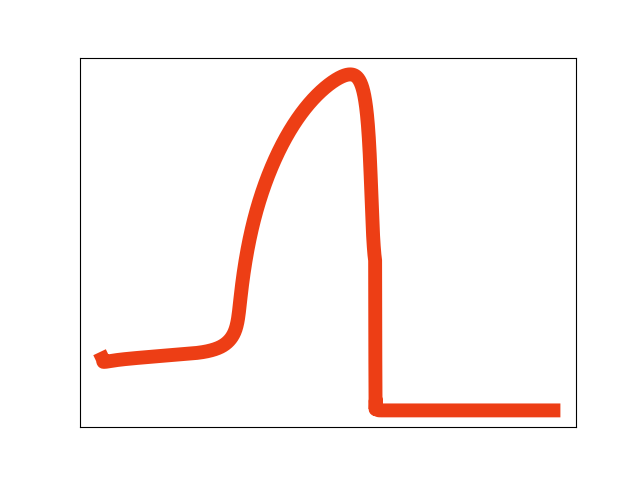}
    \includegraphics[width=0.12\linewidth]{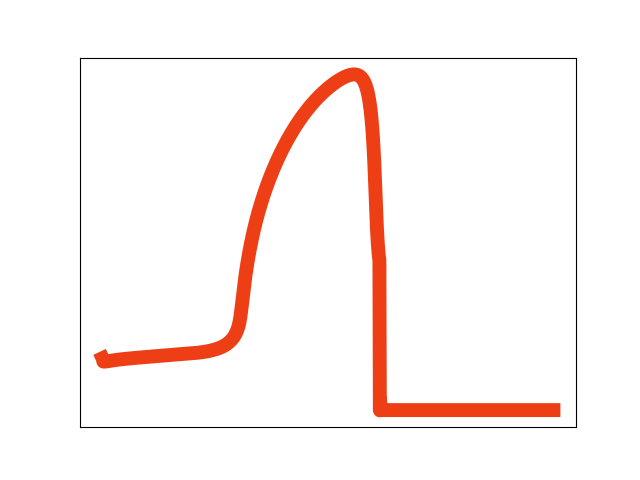}
    \includegraphics[width=0.12\linewidth]{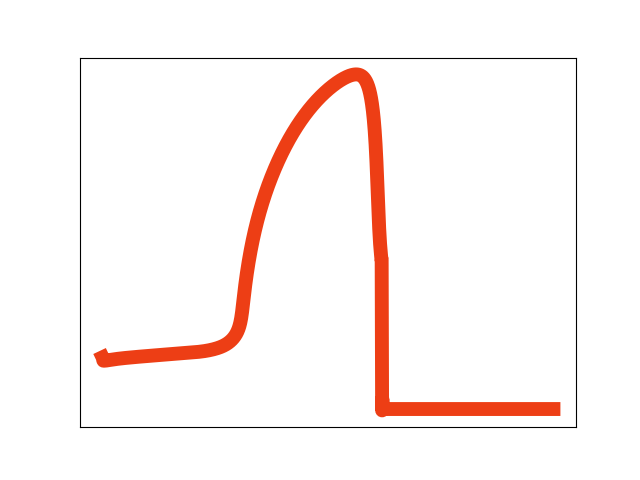}
    \includegraphics[width=0.12\linewidth]{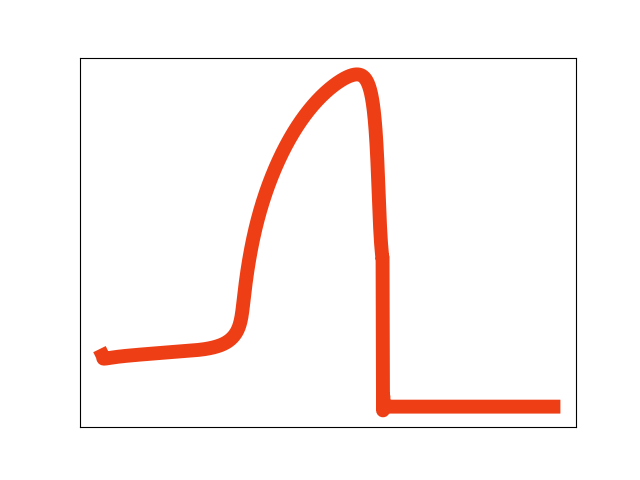} \\

    \includegraphics[width=0.12\linewidth]{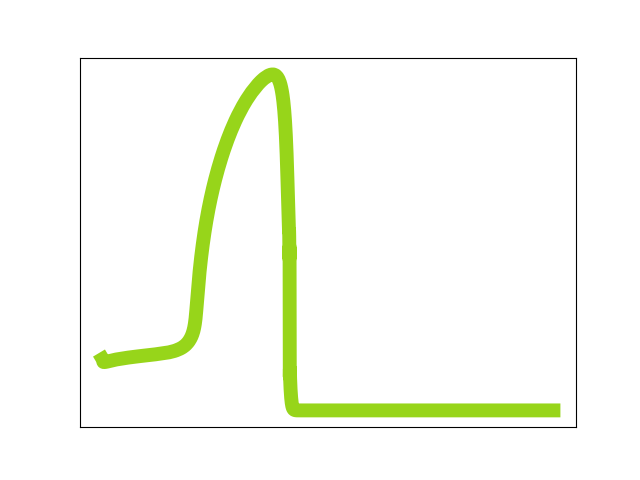}
    \includegraphics[width=0.12\linewidth]{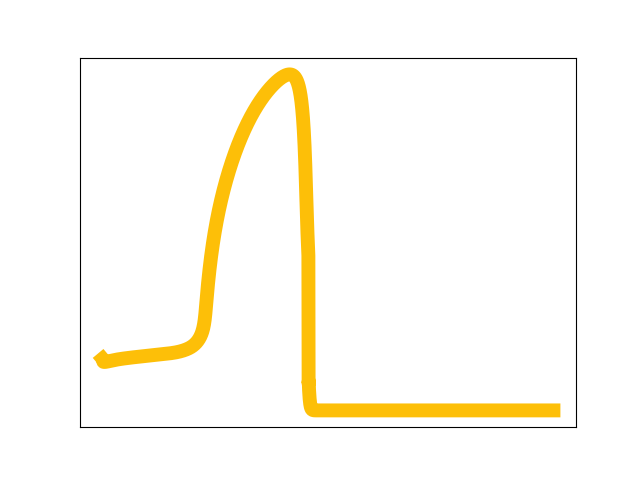}
    \includegraphics[width=0.12\linewidth]{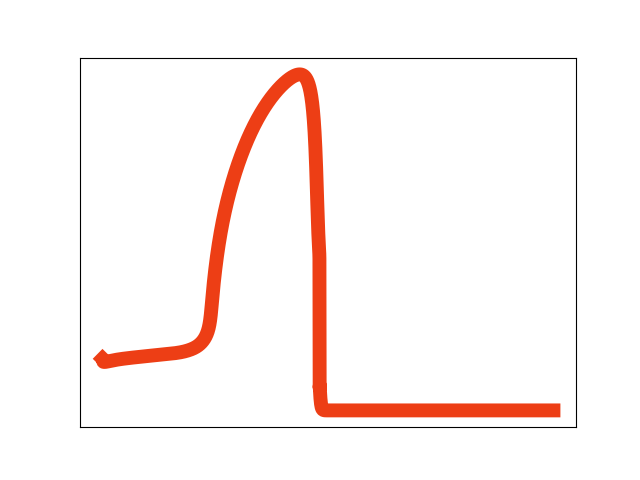}
    \includegraphics[width=0.12\linewidth]{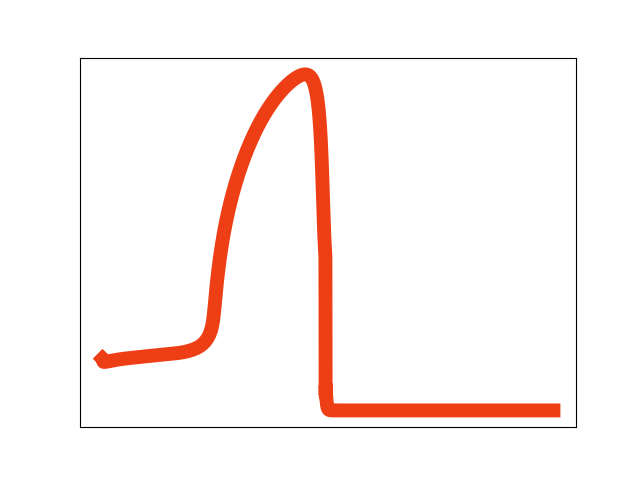}
    \includegraphics[width=0.12\linewidth]{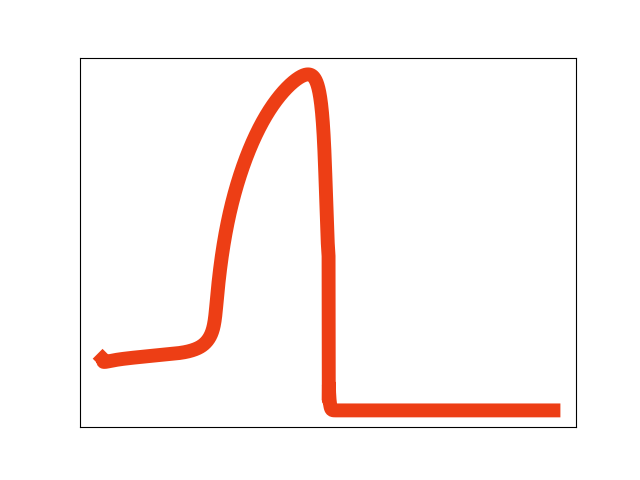}
    \includegraphics[width=0.12\linewidth]{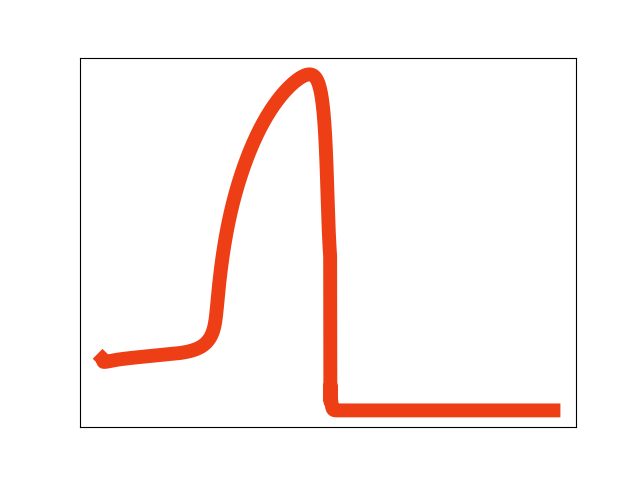}
    \includegraphics[width=0.12\linewidth]{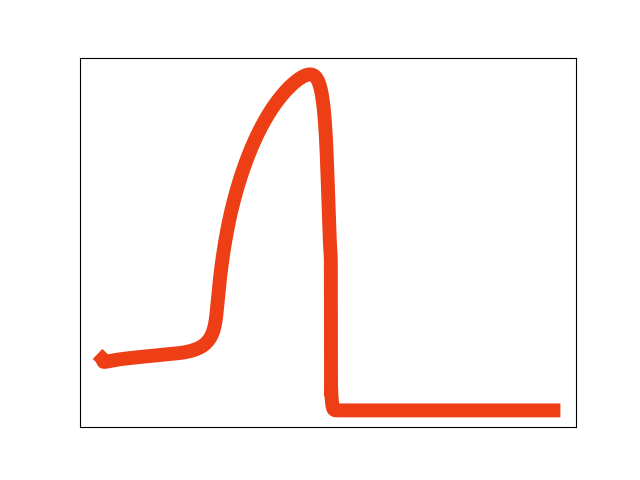} \\

    \includegraphics[width=0.12\linewidth]{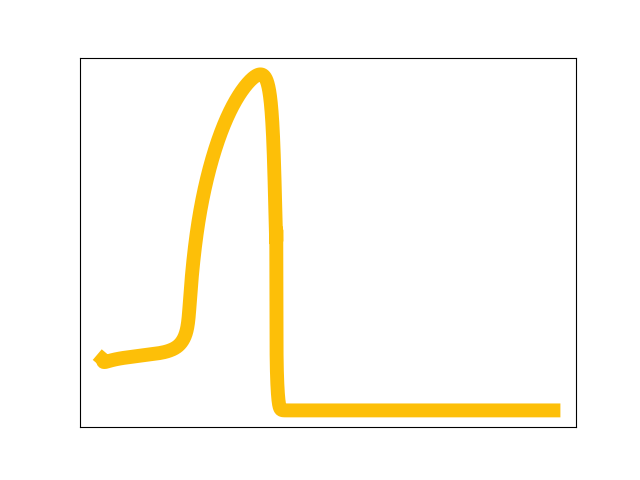}
    \includegraphics[width=0.12\linewidth]{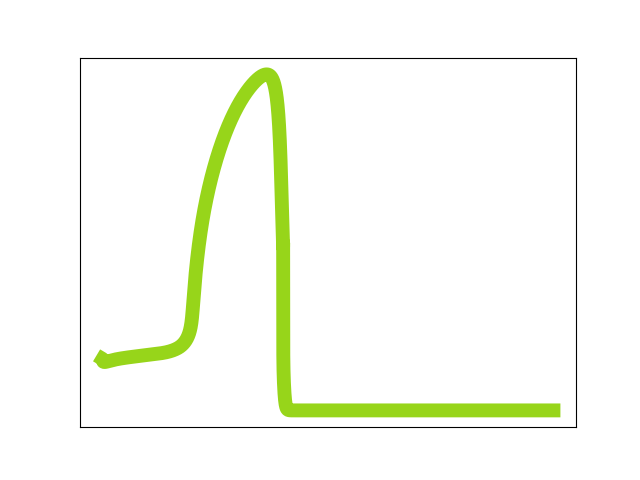}
    \includegraphics[width=0.12\linewidth]{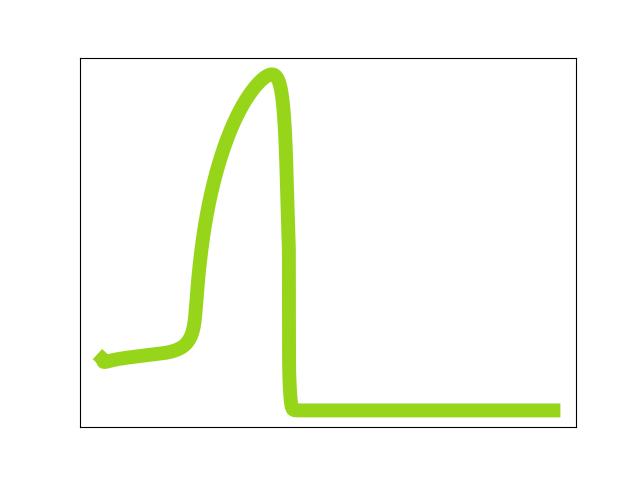}
    \includegraphics[width=0.12\linewidth]{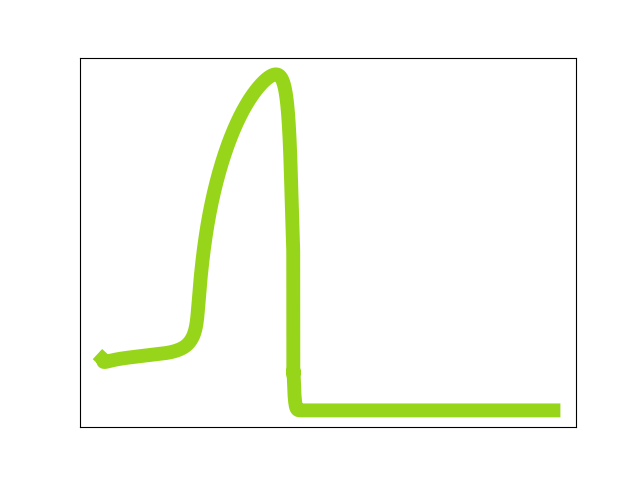}
    \includegraphics[width=0.12\linewidth]{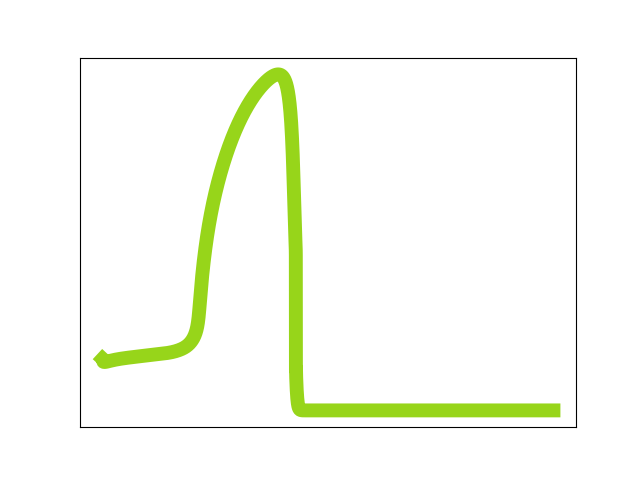}
    \includegraphics[width=0.12\linewidth]{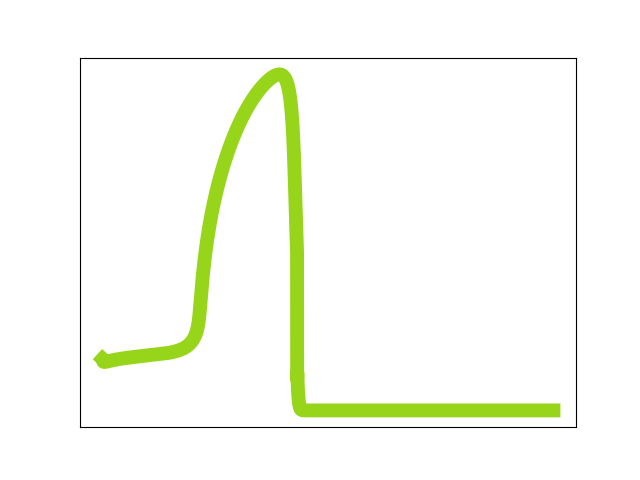}
    \includegraphics[width=0.12\linewidth]{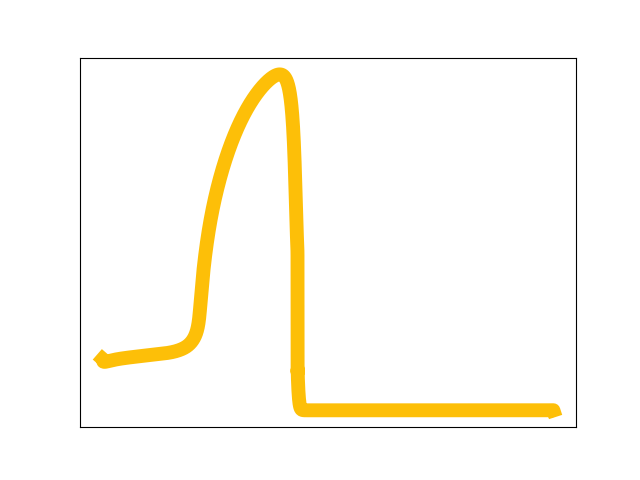} \\

    \includegraphics[width=0.12\linewidth]{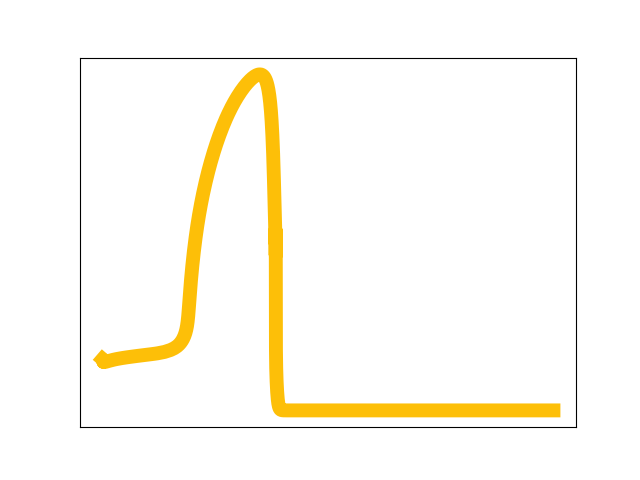}
    \includegraphics[width=0.12\linewidth]{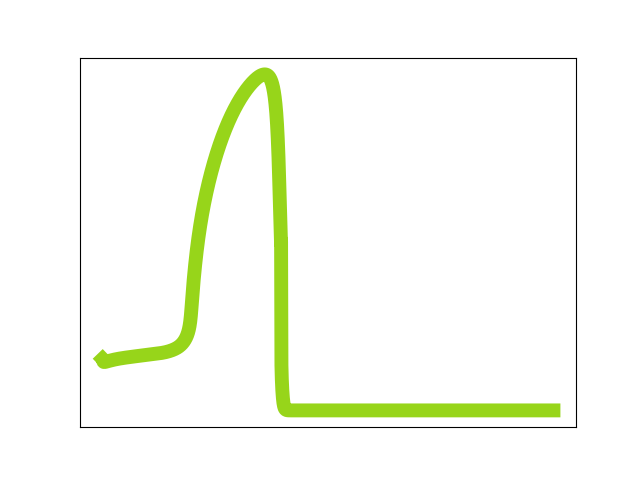}
    \includegraphics[width=0.12\linewidth]{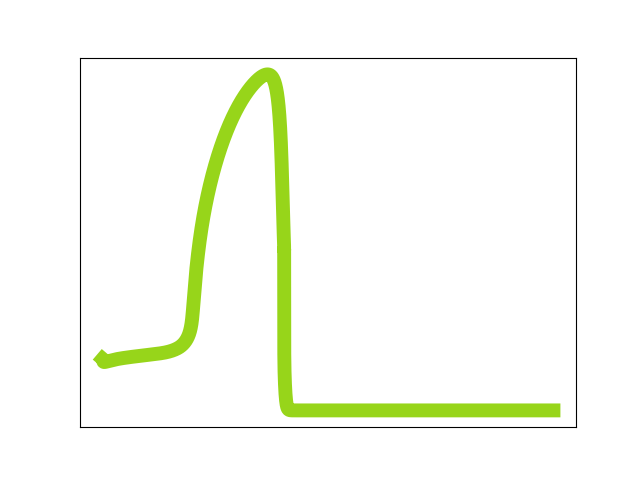}
    \includegraphics[width=0.12\linewidth]{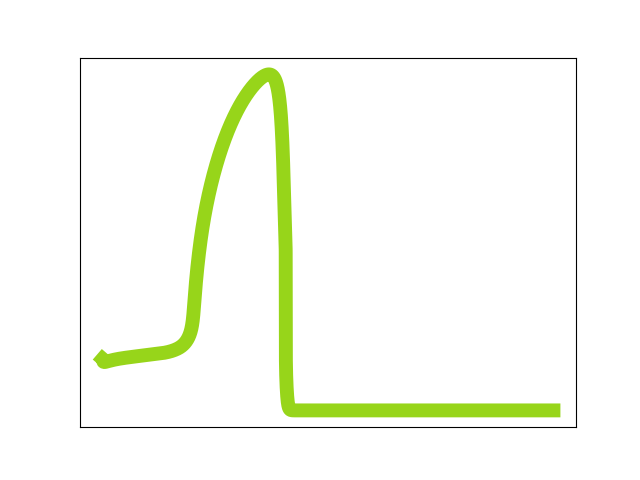}
    \includegraphics[width=0.12\linewidth]{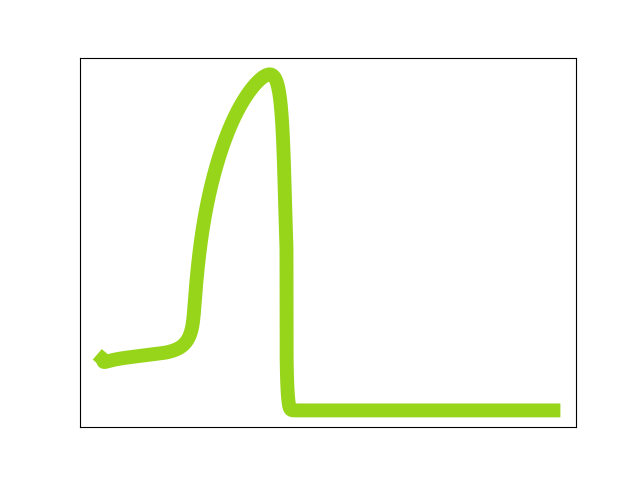}
    \includegraphics[width=0.12\linewidth]{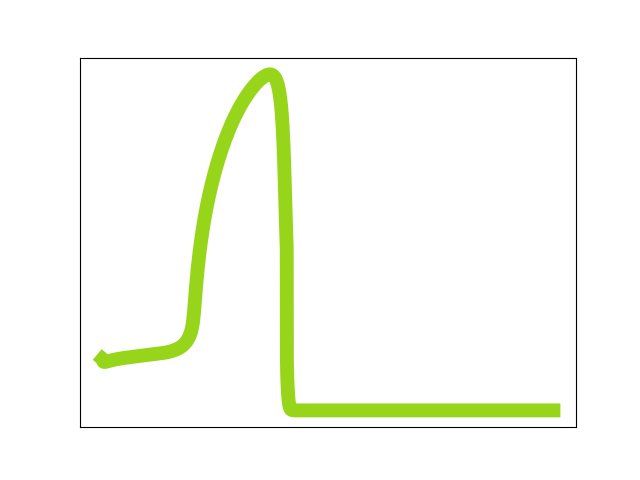}
    \includegraphics[width=0.12\linewidth]{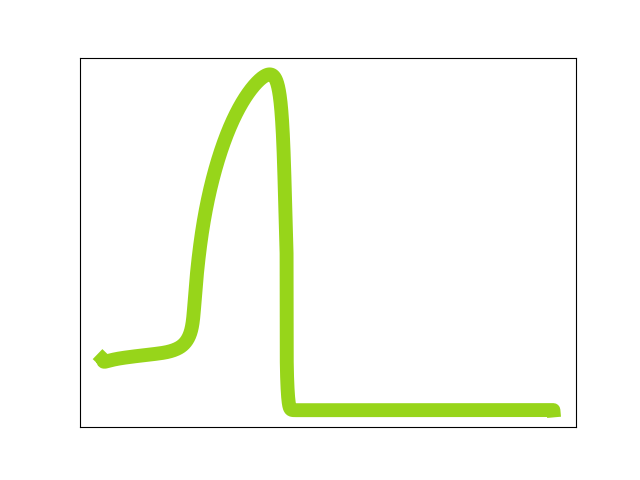} \\

    \includegraphics[width=0.12\linewidth]{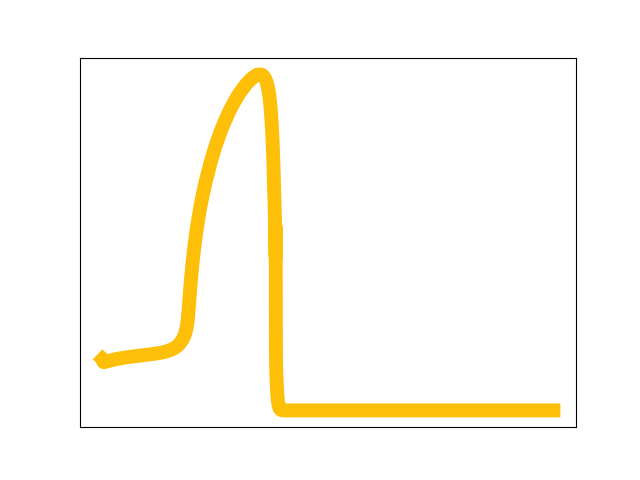}
    \includegraphics[width=0.12\linewidth]{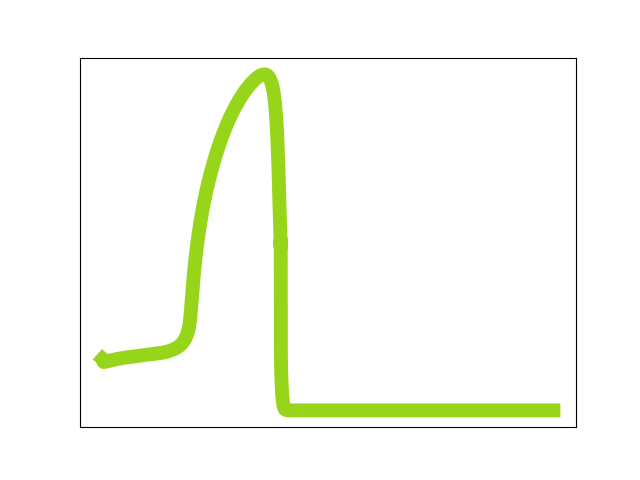}
    \includegraphics[width=0.12\linewidth]{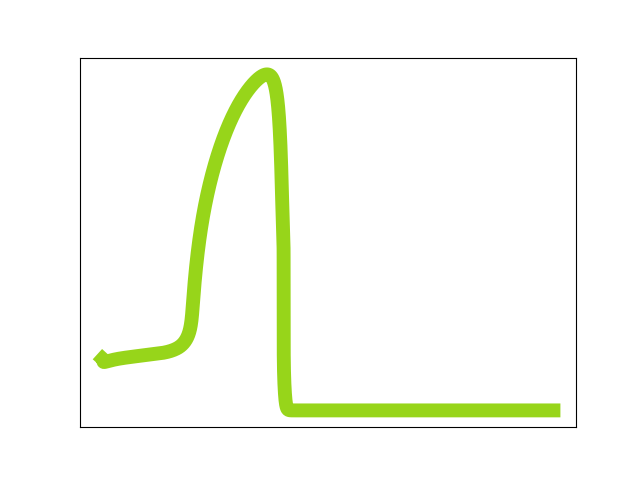}
    \includegraphics[width=0.12\linewidth]{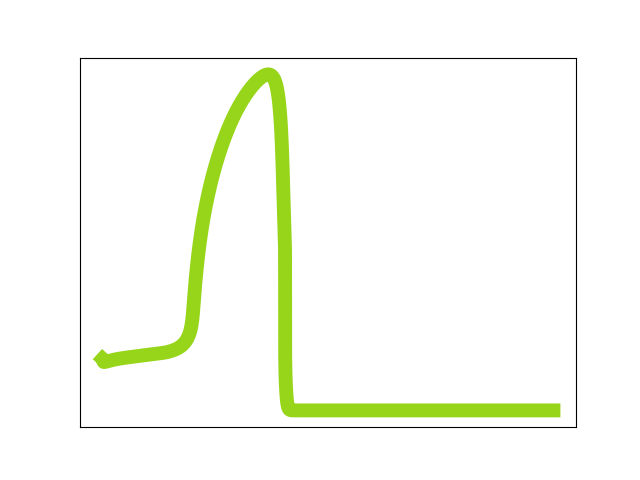}
    \includegraphics[width=0.12\linewidth]{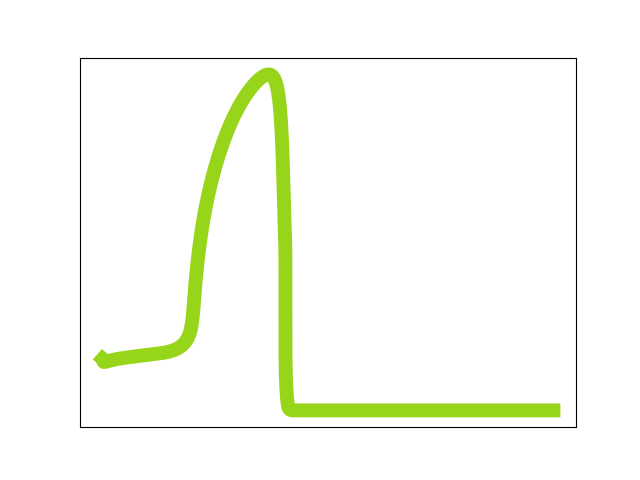}
    \includegraphics[width=0.12\linewidth]{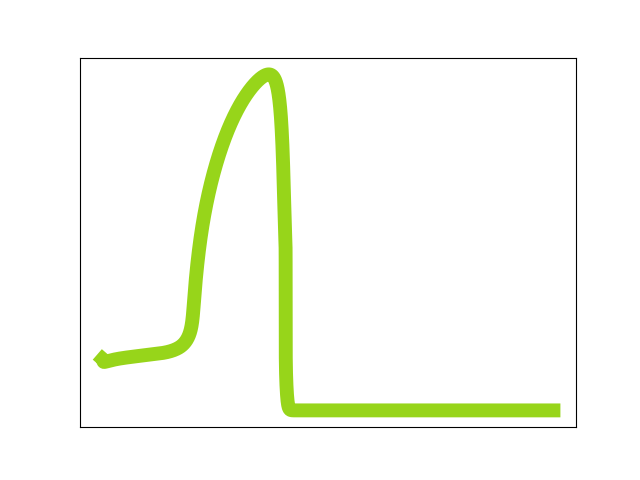}
    \includegraphics[width=0.12\linewidth]{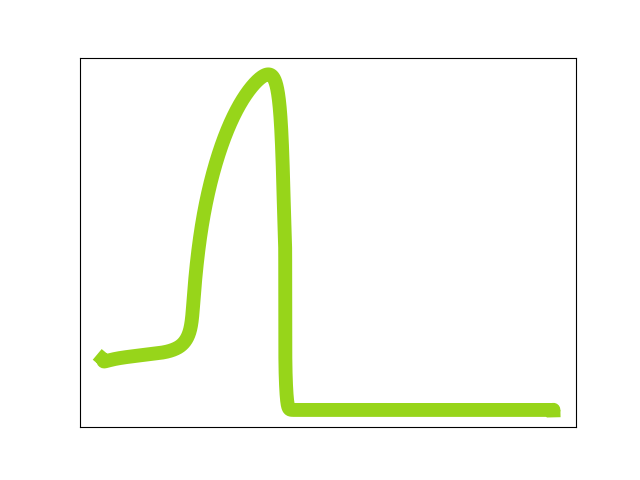}
    \end{center}
    \textbf{B} 
    \begin{center}
    \begin{tabular}{ c | c c c c c c c | c}
        \toprule
        \diagbox{$N$}{$\Delta t$} & 12.5 & 6.25 & 3.125 & 1.563 & 0.781 & 0.391 & 0.195 & $\Delta \bar{v}_{\rm CSD}$ \\
        \midrule
        1000  &7.985 &8.615 &9.385 &9.738 &9.938 &10.046 &10.092 & --\\
        2000  &6.285 &6.854 &7.223 &7.431 &7.538 &7.600 &7.638 &2.454\\
        4000  &5.135 &5.636 &5.938 &6.096 &6.181 &6.227 &6.246 &1.392\\
        8000  &4.796 &4.980 &5.128 &5.242 &5.311 & 5.349 & 5.366 & 0.880 \\
        16000 &4.793 &4.935 &5.020 &5.060 & 5.080 & 5.091 & 5.096 & 0.270 \\
        32000 &4.798 &4.942 &5.015 &5.053 & 5.063 & 5.076 & 5.084 & 0.012 \\ 
        \midrule
        $\Delta \bar{v}_{\rm CSD}$ & -- &0.144 & 0.073 &0.038 &0.10 & 0.013& 0.008& \\
        \bottomrule
    \end{tabular}
    \end{center}
    \textbf{C} 
    \begin{center}
    \begin{tabular}{ c | c c c c c c c | c}
        \toprule
        \diagbox{$N$}{$\Delta t$} & 12.5 & 6.25 & 3.125 & 1.563 & 0.781 & 0.391 & 0.195 & $\Delta w_{\rm CSD}$ \\
        \midrule
        1000  &3.500 &3.780 &4.120 &4.270 &4.360 &4.400 &4.420 & -- \\
        2000  &2.745 &3.005 &3.170 &3.260 &3.310 &3.335 &3.350 & 1.070 \\
        4000  &2.260 &2.480 &2.607 &2.680 &2.715 &2.735 &2.745 & 0.605 \\
        8000  &2.121 &2.199 &2.261 &2.310 &2.340 &2.355 &2.362 & 0.383 \\
        16000 &2.119 &2.183 &2.217 &2.234 &2.242 &2.247 &2.249 & 0.113 \\
        32000 &2.121 &2.183 &2.215 &2.231 &2.236 &2.242 &2.244 & 0.005 \\ 
        \midrule
        $\Delta w_{\rm CSD}$ & -- &0.062 &0.032 &0.016 &0.005 &0.006 &0.002 & \\
        \bottomrule
    \end{tabular}
    \end{center}
    \caption{CSD wave characteristics and quantities of interest under refinement
        in space (rows) and time (columns). \textbf{A}: Neuron potential $\phi_n(x,
        50)$ (mV) versus $x \in \Omega$ (mm), where green, yellow, and red represent
        wave speeds differing by respectively $\pm 5\%$, $\pm 15\%$, and more than $\pm
        15\%$ from our estimated wave speed of $5.1$ mm/min.  \textbf{B}: CSD
        mean wave speed $\bar{v}_{\rm CSD}$ (mm/min) and difference $\Delta
        \bar{v}_{\rm CSD}$ between consecutive refinements. \textbf{C}: Wave width
        $w_{\rm CSD}$ (mm) at $t=50$ s and difference $\Delta w_{\rm CSD}$ between
        consecutive refinements. Numerical scheme: Strang splitting, BDF2, ESDIRK4.}
        \label{fig:PDE_BDF2:ODE_ESDIRK4}
\end{figure}

We begin by considering a reference scheme -- based on Strang splitting, a BDF2
method for the PDE time-discretization, and ESDIRK4 for the ODE time-stepping
-- to compute the mean speed, the (spatial) width, and the (temporal) duration
of the CSD wave (cf.~Section~\ref{sec_CSD_wave_charac}) during temporal and
spatial refinement, before discussing how variations in terms of splitting
scheme, time-stepping and higher order spatial discretization affect accuracy
and convergence. Specifically, we apply the reference scheme and calculate the
quantities of interest for different mesh resolutions and time steps: $\Delta
x_N = 10/N$ mm for $N = 1000, 2000, 4000, 8000, 16000, 32000$ and $\Delta t_i =
12.5/i$ s for $i = 1, 2, 4, 8, 16, 32, 64$. The results are presented in
Figure~\ref{fig:PDE_BDF2:ODE_ESDIRK4}.  Wave speeds are converted from the
native m/s to mm/min for interpretability.

In general, the computed mean wave speed and width increases with decreasing
$\Delta t$, and decreases with decreasing $\Delta x$: the smaller the time
step, the faster and wider the wave, while the smaller the mesh size, the
slower and narrower the wave (Figure~\ref{fig:PDE_BDF2:ODE_ESDIRK4}A).
Regarding the mean wave speed (Figure~\ref{fig:PDE_BDF2:ODE_ESDIRK4}B), we
observe that the computed values vary substantially, ranging from $5.063$ to
$10.090$ mm/min. We observe that the spatial errors, estimated by proxy by the
difference between consecutive spatial refinements, dominate the temporal
errors/differences. In particular, the differences $\Delta \bar{v}_{\rm CSD}$
in wave speed between the coarsest mesh sizes $N = 1000$ and $2000$ are in the
range $1.70-2.45$ mm/min.  Conversely, the differences in wave speed between
the coarsest time steps $\Delta t = 12.5$ and $6.25$ are in the range
$0.51-0.64$ mm/min. For the ultimate spatial and temporal refinement level, we
observe a difference of $0.012$ and $0.008$ mm/min, respectively. Finally, we
observe that the wave speed seems to converge in space and time with an
estimated wave speed of $5.1$ mm/min: the differences $\Delta \bar{v}_{\rm
CSD}$ reduce as $\Delta x$ and $\Delta t$ is reduced. There is however no clear
rate of convergence.

Similarly, we observe large variations in the computed (spatial) CSD wave
width, ranging from $1.121$ to $4.420$ mm
(Figure~\ref{fig:PDE_BDF2:ODE_ESDIRK4}C). The difference ($\Delta w_{\rm CSD}$)
varies in the range $0.755-1.07$ mm and $0.062-0.280$ mm between respectively
the coarsest mesh sizes $N=1000$ and $2000$ and the coarsest time steps $\Delta
t = 12.5$ and $6.25$ s. For the finest time and mesh discretizations, we
observe differences of $0.005$ mm and $0.002$ mm, respectively. Moreover, the
differences $\Delta w_{\rm CSD}$ reduce as the mesh size and time step are reduced.
As with the mean wave speed, there is no clear rate of convergence. In contrast, the (temporal) duration $d_{\rm CSD}$ of the CSD wave
does not change substantially during refinement in space and time: we
observe a duration of elevated extracellular potassium levels of
$17-18$ s for all the spatial and temporal resolutions considered
(results not shown). Finally, we observe that this implicit
higher-order (reference) scheme behaves qualitatively similar to the
implicit lower-order scheme -- based on Godunov splitting, a BE method
for the PDE time-discretization, and BE for the ODE time-stepping
as shown in Figure~\ref{fig:intro}.

\subsection{Choice of discretizations}
\label{sec_choice_disc}

We turn to evaluate the effect of discretization choices in terms of splitting
scheme, time-stepping, and higher order spatial discretization.

\subsubsection*{Splitting scheme}
To evaluate the second order Strang splitting scheme, we compared the
computed CSD wave characteristics, peak neuron potentials, and wave
speeds with those computed using a first order Godunov splitting
scheme (cf.~Section~\ref{sec:sub:scheme:strang}). The computed wave
speeds are comparable at given resolutions, and we observe a similar
difference decay ($\Delta \bar{v}$) for the Strang and Godunov
splitting schemes (see Table~\ref{tab:godunov} in Supplementary Tables).

\subsubsection*{Time stepping}
To assess how the choice of PDE time stepping affects the convergence of the
numerical schemes, we repeat the convergence study presented above in
Section~\ref{sec:sub:convergence:res} replacing the BDF2 scheme by a
Crank-Nicolson (CN) scheme and ESDIRK4 by RK4 for the ODE time stepping
(Table~\ref{tab:cn}). We note that choosing an explicit ODE time stepping
scheme (RK4) here is based on the observation that CN for the PDE time stepping in
combination with implicit ODE time stepping schemes results in a diverging
(non-linear) ODE solver.

Importantly, we note that the PDE Newton solver fails to
converge for $\Delta t \geq 1.563$ ms for this scheme. Again, we observe that
the computed speed and width of the wave decreases with decreasing $\Delta x$,
and is essentially constant with decreasing $\Delta t$ (likely due to the
already fine timestep required for convergence of the PDE Newton solver).  The
spatial errors, estimated by proxy by the difference between consecutive
spatial refinements, are comparable to those reported for BDF2, and again we
obtain an estimated wave speed of $5.1$ mm/min for this model scenario.
\begin{table}[ht]
    \textbf{A} 
    \begin{center}
    \begin{tabular}{ c | c c c c c c c | c}
        \toprule
        \diagbox{$N$}{$\Delta t$} & 12.5 & 6.25 & 3.125 & 1.563 & 0.781 & 0.391 & 0.195 & $\Delta \bar{v}_{\rm CSD}$ \\
        \midrule
        1000  & $*$ & $*$ & $*$ & $*$ &10.154 &10.154 &10.154 & -- \\
        2000  & $*$ & $*$ & $*$ & $*$ &7.669 &7.669 &7.669 &2.485 \\
        4000  & $*$ & $*$ & $*$ & $*$ &6.269 &6.269 &6.269 &1.400 \\
        8000  & $*$ & $*$ & $*$ & $*$ &5.384 &5.386 &5.386 &0.883 \\
        16000 & $*$ & $*$ & $*$ & $*$ &5.102 &5.103 &5.103 &0.283 \\
        32000 & $*$ & $*$ & $*$ & $*$ &5.090 &5.090 &5.091 &0.012 \\
        \midrule
        $\Delta \bar{v}_{\rm CSD}$ & -- & -- & -- & -- & -- &0.000 &0.001 & \\
        \bottomrule
    \end{tabular}
    \end{center}
    \textbf{B} 
    \begin{center}
    \begin{tabular}{ c | c c c c c c c | c}
        \toprule
        \diagbox{$N$}{$\Delta t$} & 12.5 & 6.25 & 3.125 & 1.563 & 0.781 & 0.391 & 0.195 & $\Delta w_{\rm CSD}$ \\
        \midrule
        1000  & $*$ & $*$ & $*$ & $*$ &4.450 &4.450 &4.450 & -- \\
        2000  & $*$ & $*$ & $*$ & $*$ &3.365 &3.365 &3.365 &1.085 \\
        4000  & $*$ & $*$ & $*$ & $*$ &2.752 &2.755 &2.755 &0.610 \\
        8000  & $*$ & $*$ & $*$ & $*$ &2.370 &2.371 &2.371 &0.384 \\
        16000 & $*$ & $*$ & $*$ & $*$ &2.252 &2.252 &2.252 &0.119 \\
        32000 & $*$ & $*$ & $*$ & $*$ &2.247 &2.247 &2.247 &0.005 \\
        \midrule
        $\Delta w_{\rm CSD}$ & -- & -- & -- & -- & -- &0.000 &0.000 & \\
        \bottomrule
    \end{tabular}
    \end{center}
    \caption{CSD wave characteristics and quantities of interest under refinement
        in space (rows) and time (columns). \textbf{A}: CSD
        mean wave speed $\bar{v}_{\rm CSD}$ (mm/min) and difference $\Delta
        \bar{v}_{\rm CSD}$ between consecutive refinements. \textbf{B}: Wave width
        $w_{\rm CSD}$ (mm) at $t=50$ s and difference $\Delta w_{\rm CSD}$ between
        consecutive refinements. Numerical scheme: Strang splitting, CN, RK4.
        $*$ indicates
        that the solver failed to converge.}
    \label{tab:cn}
\end{table}

To assess how the choice of ODE time stepping affects the convergence
of the numerical schemes, we repeat the convergence study presented
above in Section~\ref{sec:sub:convergence:res} replacing ESDIRK4 by a
first-order Backward Euler (BE) scheme or an explicit 4th order
Runge-Kutta (RK4) method.  For the RK4 method, we observe nearly
indistinguishable results as for ESDIRK4 with the important
distinction that the RK4 solvers fail to converge for $\Delta t \geq
3.125$ ms (see Supplementary Tables, Table~\ref{tab:ode}A). Also for
the BE scheme, the computed wave speed values are similar, but with
slightly bigger differences in terms of differences between
resolutions (see Supplementary Tables, Table~\ref{tab:ode}B).

\subsubsection*{Higher order spatial discretization}
To assess how the polynomial degree of the finite elements affects the
convergence of the numerical schemes, we repeat the convergence study
presented in Section~\ref{sec:sub:convergence:res} replacing the
lowest order finite element pairings by higher order pairings. We
consider the model described in Section~\ref{sec:mathmodel}
discretized with discontinuous piecewise linear polynomials ($P_1$) for
the volume fractions $\alpha_r$ and continuous piecewise linear
polynomials of degree 2 ($P_2^c$) for the ion concentrations and
potentials. The finite element software used within this work
(FEniCS\cite{logg2012automated}) has automated functionality for
solving coupled PDE-ODE systems via the
PointIntegralSolver\cite{farrell2019automated} class. The
PointIntegralSolver solves the ODEs at the vertices of the elements
and only supports the use of elements with degrees of freedom located
at the vertices. Higher order elements are thus not supported. To assess the
accuracy and convergence of the higher order scheme, we here consider an
alternative implementation where the ODEs are solved at each degree of
freedom of the spatial discretization using a backward Euler scheme
and a Newton solver.

\begin{table}[ht]
    \textbf{A} 
    \begin{center}
    \begin{tabular}{ c | c c | c c}
         \toprule
         &  \multicolumn{2}{c}{$P_0$-$P_1^c$} & \multicolumn{2}{c}{$P_1$-$P_2^c$}\\
         \midrule
         $N$ & $\bar{v}_{\rm CSD}$ & $\Delta \bar{v}_{\rm CSD}$ & $\bar{v}_{\rm CSD}$ & $\Delta \bar{v}_{\rm CSD}$\\
         \midrule
         1000 &7.585 & -- & 6.415 & --\\
         2000 &6.346 &1.239 & 5.416 &0.999\\
         4000 &5.412 &0.934&5.005 &0.411\\
         8000 &5.005 &0.407&4.987 &0.018\\
         16000 &4.988 &0.017 &4.989 &-0.002\\
         32000 &4.989 &-0.001 & --& --\\
        \bottomrule
    \end{tabular}
    \end{center}
    \textbf{B} 
    \begin{center}
    \begin{tabular}{ c | c c | c c}
        \toprule
        &  \multicolumn{2}{c}{$P_0$-$P_1^c$} & \multicolumn{2}{c}{$P_1$-$P_2^c$}\\
        \midrule
        $N$ & $w_{\rm CSD}$ & $\Delta w_{\rm CSD}$ & $w_{\rm CSD}$ & $\Delta w_{\rm CSD}$\\
        \midrule
        1000 &3.330 & -- & 2.800& --\\
        2000 &2.785 &0.545 &2.380 &0.420\\
        4000 &2.380 &0.405  &2.210 &0.170\\
        8000 &2.210 &0.170 &2.202 &0.008\\
        16000 &2.203 &0.007 &2.203 &-0.001\\
         32000 &2.203 &0.000 & --& --\\
        \bottomrule
    \end{tabular}
    \end{center}
    \caption{Comparison of CSD wave characteristics and quantities of interest
    under refinement in space between the $P_0$-$P_1^c$ spatial discretization
    described in Section~\ref{sec_spatial_discretization} and a $P_1$-$P_2^c$
    discretization (both solved using the implementation allowing for higher
    order elements). \textbf{A}: CSD mean wave speed $\bar{v}_{\rm CSD}$
    (mm/min) and difference $\Delta \bar{v}_{\rm CSD}$ between consecutive
    refinements. \textbf{B}: Wave width $w_{\rm CSD}$ (mm) at $t=50$ s and
    difference $\Delta w_{\rm CSD}$ between consecutive refinements. Numerical
    scheme: Strang splitting, BDF2, BE using $\Delta t =3.125$ ms.}
    \label{tab:CG2}
\end{table}

In Table~\ref{tab:CG2}, we compare the CSD wave characteristics (wave
speed $\bar{v}_{\rm CSD}$ and wave width $w_{\rm CSD}$) computed after solving
the system with the two different spatial discretizations: $P_0$-$P_1^c$ and $P_1$-$P_2^c$.
The numerical scheme consists of a second order Strang splitting, together with
a BDF2 time-stepping scheme for the PDE and a backward Euler (BE) scheme for
the system of ODEs, with a timestep of $\Delta t = 3.125$ ms. We observe that
the wave speed and wave width computed using a second order $P_1$-$P_2^c$ spatial
discretization with a mesh of $N$ elements (where
$N\in\{1000,2000,4000,8000,16000\}$) are similar to the quantities obtained by
a $P_0$-$P_1^c$ discretization with $2N$ elements (i.e.~the same number of degrees of
freedom).

We also study the approximation of the neuron potential by arbitrary
higher-order polynomials to compare low-degree finite element spatial
discretization against higher order discretizations or global spectral
methods\cite{gottlieb1977numerical}. First, the system of equations is solved
using a $P_0$-$P_1^c$ finite element scheme with a large spatial discretization of
$N=32000$. We then use the Chebfun software to approximate the neuron potential
with Chebyshev polynomials\cite{driscoll2014chebfun,floater2007barycentric}.
\begin{figure}
\centering
        \begin{overpic}[width=0.45\linewidth]{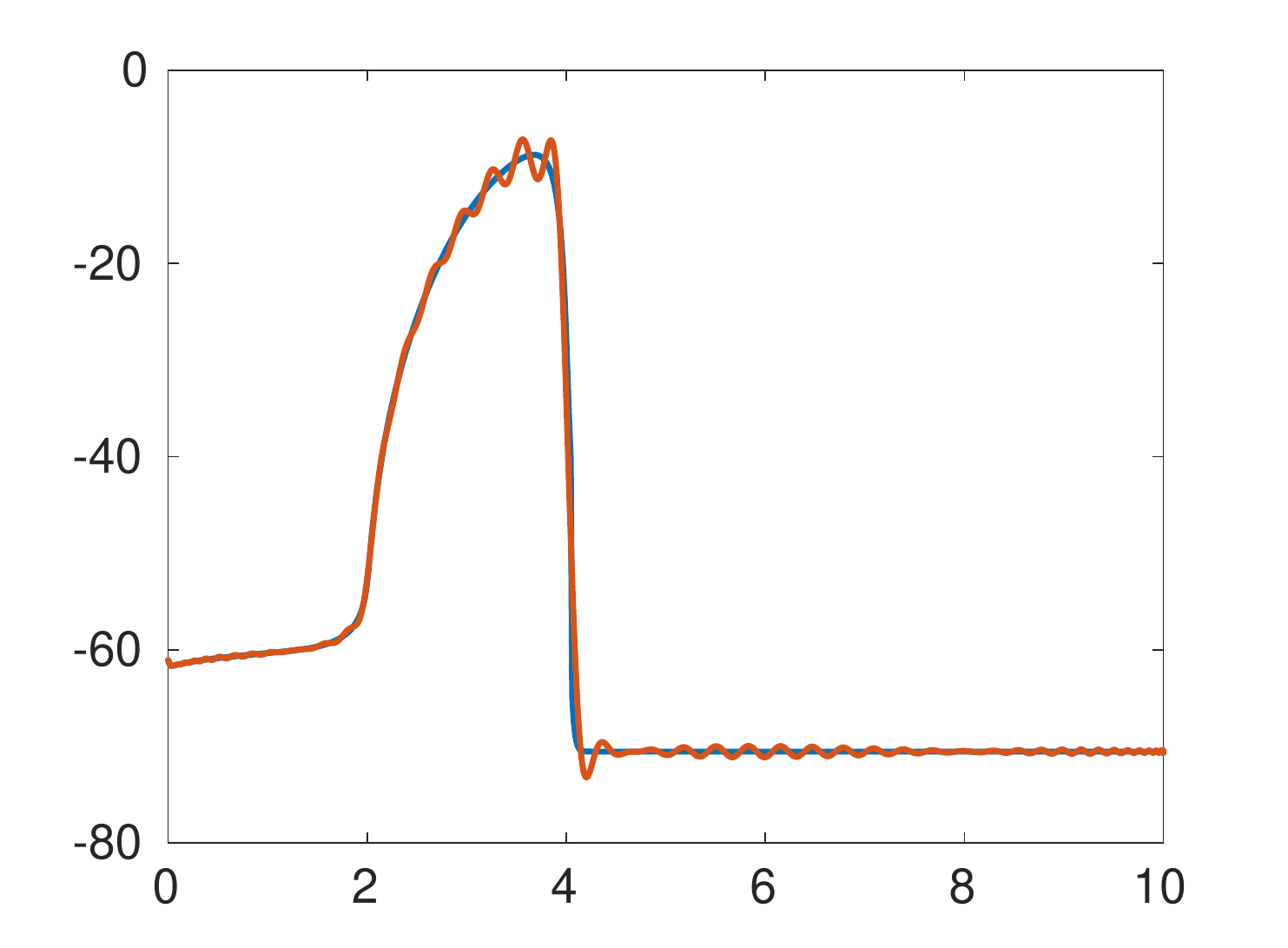}
        \put(-3,25){\rotatebox{90}{$\phi_N$ (mV)}}
        \put(42,-3){$x$ (mm)}
        \end{overpic}
        \hspace{0.4cm}
        \begin{overpic}[width=0.45\linewidth]{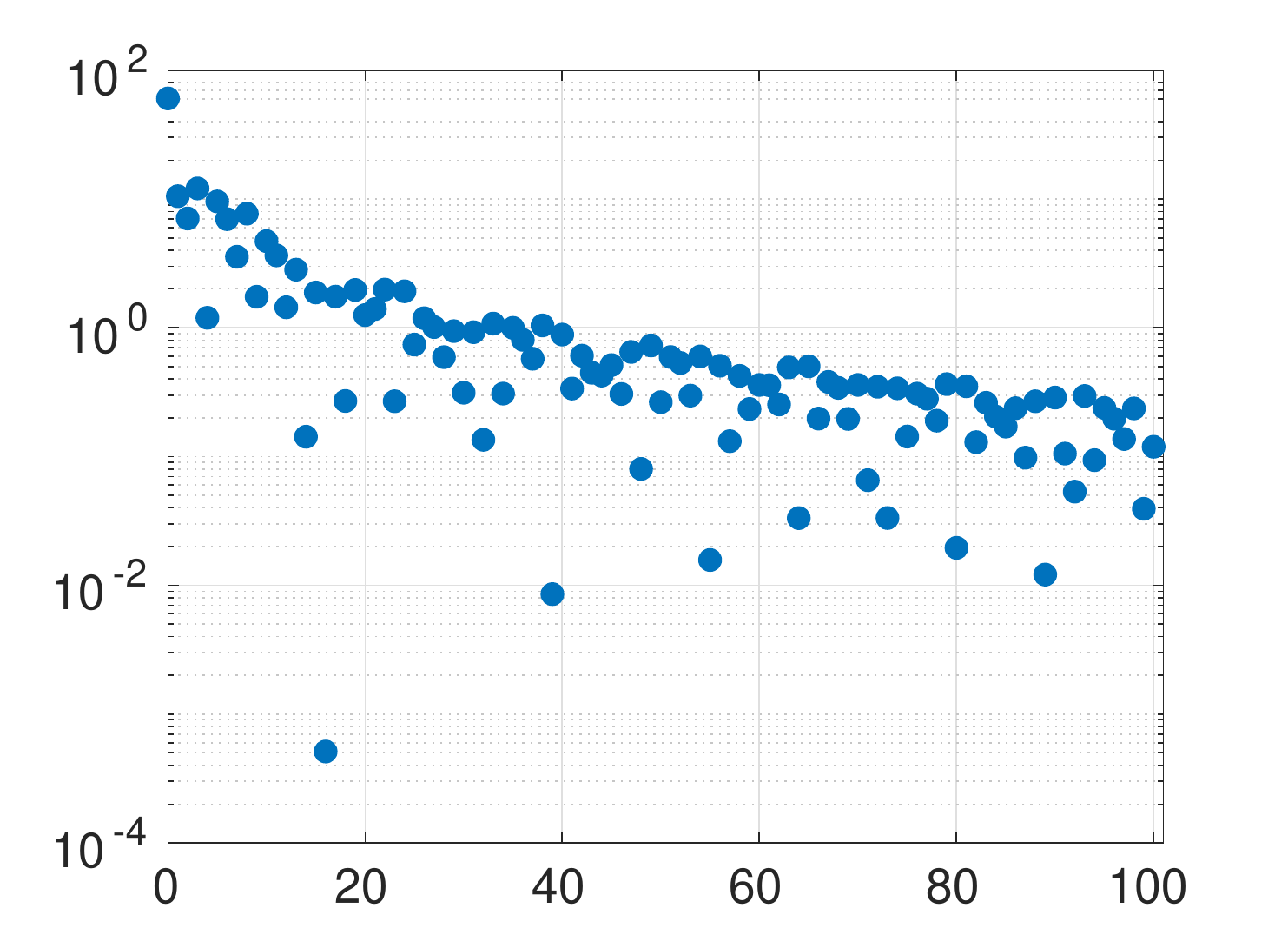}
        \put(23,-3){Degree of polynomial}
        \put(-3,7){\rotatebox{90}{Magnitude of coefficient}}
        \end{overpic}
        \vspace{0.2cm}
    \caption{Left: Neuron potential $\phi_n$ (blue) together with its
    polynomial approximation of degree $100$ (red) computed using the Chebfun
    software\protect\cite{driscoll2014chebfun}. Right: Magnitude of the Chebyshev
    coefficients of the polynomial approximation to $\phi_n$.}
    \label{fig:poly:approx}
\end{figure}

In Figure~\ref{fig:poly:approx} (left), we display the original function $\phi_n$ (in
blue), together with its polynomial approximation of degree $100$ (in red). We
see that the approximation has large oscillating errors near the wave front
located at $x=4$ mm, which is likely due to Wilbraham--Gibbs phenomenon\cite{trefethen2019approximation}. On the right panel of Figure~\ref{fig:poly:approx}, we report the
coefficients of the polynomial approximation to $\phi_n$ of degree $100$
constructed by Chebfun. The magnitude of the Chebyshev coefficients seems to
decay algebraically at a rate of $-1.2$, suggesting that a very large
polynomial degree (of order $10^4-10^5$) would be required to approximate the
solution to~\eqref{eq:alpha:r}
and~\eqref{eq_conservation_rR}--\eqref{eq_phi_rR} accurately. 

\subsection{Numerical performance for the zero flow limit}
\label{sec:performance}
So far, we have investigated the convergence of the different schemes via studying key functionals of the solution. To evaluate the performance
and scalability of the implementation of the schemes for the model in the zero
flow limit, we consider an additional set of experiments measuring the memory
usage and CPU timings\footnote{Timings were performed on a Lenovo ThinkPad
2.70GHz x 4 Intel Core i7-7500U CPU using FEniCS 2019.1.0 without
parallelization.}. We evaluate both the standard implementation and the
implementation allowing for higher order finite elements. For the standard
implementation, we consider simulations with BDF2 and ESDIRK4 time stepping for
the PDEs and the ODEs, respectively, and compare the performance of Godunov and
Strang splitting.  For the implementation allowing for higher order finite
elements we consider simulations with Strang splitting, BDF2 time PDE stepping, 
BE time ODE stepping, and compare performance of using $P_0$-$P_1^c$ and $P_1$-$P_2^c$
elements.

\begin{table}[ht]
    \textbf{A}
    \begin{center}
    \begin{tabular}{cccccccc}
    \toprule
        $N$ & Dofs & M (MiB) & T$_{\rm A}$ (s) &  T$_{\rm LU}$ (s) & T$_{\rm PDE}$ (s) & T$_{\rm ODE}$ (s) & T$_{\rm tot}$ (s)\\
    \midrule
        8000 & 72008  & 117  & 329 & 350 & 679  & 129 & 813 \\
        16000 & 144008 & 199 &  659 &  707 & 1366 & 258 & 1632 \\
        32000 & 288008 & 396 & 1438  & 1520 & 2958 & 559 & 3526 \\
    \bottomrule
    \end{tabular}
    \end{center}
    \textbf{B}
    \begin{center}
    \begin{tabular}{cccccccc}
    \toprule
        $N$ & Dofs & M (MiB) & T$_{\rm A}$ (s) &  T$_{\rm LU}$ (s) & T$_{\rm PDE}$ (s) & T$_{\rm ODE}$ (s) & T$_{\rm tot}$ (s)\\    
        \midrule
        8000  & 72008  & 118 & 333 & 352 & 685  & 259 & 954 \\
        16000 & 144008 & 199 & 677 & 719 & 1396 & 536 & 1932 \\
        32000 & 288008 & 396 & 1367 & 1451 & 2818 & 1086 & 3916 \\
    \bottomrule
    \end{tabular}
    \end{center}
    \textbf{C}
    \begin{center}
    \begin{tabular}{cccccccc}
    \toprule
        $N$ & Dofs & M (MiB) & T$_{\rm A}$ (s) & T$_{\rm LU}$ (s) & T$_{\rm
        PDE}$ (s) & T$_{\rm ODE}$ (s) &  T$_{\rm tot}$ (s)\\    
        \midrule
        8000   & 72008  & 130 & --  & -- & 696  & 253  & 953 \\
        16000  & 144008 & 223 & --  & -- & 1411 & 506  & 1922 \\
        32000  & 288008 & 446 &  -- & -- & 3076 & 1115 & 4198 \\
    \bottomrule
    \end{tabular}
    \end{center}
    \textbf{D} 
    \begin{center}
    \begin{tabular}{cccccccc}
    \toprule
        $N$ & Dofs & M (MiB) & T$_{\rm A}$ (s) &  T$_{\rm LU}$ (s) & T$_{\rm PDE}$ (s) & T$_{\rm ODE}$ (s) & T$_{\rm tot}$ (s)\\    
        \midrule
        4000  & 72008  & 165  & -- & -- & 1017 & 309  & 1331 \\
        8000  & 144008 & 285  & -- & -- & 1693 & 511  & 2209 \\
        16000 & 288008 & 567  & -- & -- & 3902 & 1182 & 5091 \\
    \bottomrule
    \end{tabular}
    \end{center}
\caption{CPU timings and memory usage for approximating solutions in the zero
    flow limit.  Dofs: number of degrees of freedom in the linear (PDE) system,
    M: Maximal memory usage of simulation relative to baseline.  T$_{\rm A}$:
    CPU time for finite element assembly, T$_{\rm LU}$: CPU time for LU solver,
    T$_{\rm ODE}$: CPU time for ODE stepping, T$_{\rm PDE}$: CPU time for PDE
    stepping (in \textbf{A} and \textbf{B} this equals the sum of T$_{\rm A}$
    and T$_{\rm LU}$), T$_{\rm tot}$: Total CPU time for simulation. T$_{\rm
    A}$ and  T$_{\rm LU}$ are not reported for \textbf{C} and \textbf{D}. All
    simulations have $\Delta t = 3.125$ ms and final time $T=5$ s (i.e.~$1600$
    timesteps). Results from simulations with the standard implementation with
    BDF2, ESDIRK4, $P_0$-$P_1^c$, and either \textbf{A}: Godunov splitting; or
    \textbf{B}: Strang splitting. Results from simulations with the
    implementation allowing for higher order elements with BDF2, BE, Strang
    splitting, and either \textbf{C}:  $P_0$-$P_1^c$ elements; or \textbf{D}: $P_1$-$P_2^c$
    elements.} 
    \label{tab:timings} 
\end{table}

For both the Godunov and Strang splitting schemes and the standard
implementation, we observe that the memory usage increases linearly with the
number of degrees of freedom: doubling the number of degrees of freedom leads to an increase in memory of a
factor $2$ (see Table~\ref{tab:timings}A and B). We observe that the CPU time for
the simulations grows linearly with the number of degrees of freedom: doubling the number of
degrees of freedom leads to an increase in total simulation time of a factor 2.
In the Strang splitting scheme, the ODEs are solved twice for each PDE step and
we thus expect the total ODE stepping time to be greater than for the Godunov
splitting scheme (where the ODEs are only solved once per PDE step). Indeed,
the total ODE stepping time is about twice as large for Strang splitting as for
Godunov splitting (Table~\ref{tab:timings}). Conversely, the time required for
finite element assembly and LU solves is comparable for Godunov splitting and
Strang splitting. In total, the simulation time is higher for Strang splitting
than Godunov ($11 - 18 \%$). Moreover, the total simulation time is
dominated by the cost of finite element assembly and LU solves ($84 \%$ for
Godunov splitting, $72 \%$ for Strang splitting). Finally, the time required
for finite element assembly is comparable to that of the LU solves for both
splitting schemes.

The simulation time of the higher order $P_1$-$P_2^c$ discretization (with the
implementation allowing for higher order schemes) is reported in
Table~\ref{tab:timings}D. For a given number of mesh elements
$N\in\{4000,8000,16000\}$, we compare the timings with the $P_0$-$P_1^c$
discretization with a mesh of $2N$ elements (Table~\ref{tab:timings}C) so that
the underlying systems have the same number of degrees of freedom as well as
approximately the same accuracy, following the linear convergence of the
$P_1$-$P_2^c$ discretization as observed in Section~\ref{sec_choice_disc}. We then
see that the total time needed to run the full simulation is approximately
$21 - 40 \%$ higher for the $P_1$-$P_2^c$ discretization. This difference is
explained by the higher density of the finite element matrices resulting from
the $P_1$-$P_2^c$ discretization compared to the choice of $P_0$-$P_1^c$.

In conclusion, Strang and Godunov splitting yield comparable accuracy and
memory usage. We observe that Strang splitting yields higher total CPU time
than Godunov. When varying (ODE and PDE) time stepping schemes we observe
minor variations in terms of accuracy. The higher order element scheme
($P_1$-$P_2^c$) has both a higher total CPU time and memory usage than the lower
order spatial scheme ($P_0$-$P_1$), whereas the accuracy is comparable.  We find
that the accurate computation of CSD wave characteristics (wave speed and wave
width) requires a very fine spatial and fine temporal resolution for all
schemes tested.

\section{Numerical solution of model including fluid dynamics}
\label{sec:full:scheme}

The previous schemes and experiments considered only the zero flow
limit.  Here, we present a numerical scheme for the full mathematical
model. The variational formulation~\eqref{eq:full:summary} is obtained
by multiplying~\eqref{eq_alpha_rR}--\eqref{eq_phi_rR}
and~\eqref{eq:p:R} by suitable test functions, integrating over the
domain $\Omega$, integrating terms with higher order derivatives by
parts, and inserting the boundary
conditions~\eqref{bc_ionflux_u}. Further, we use a backward Euler
scheme for the PDE time discretization. Let $S_r \subset H^1(\Omega)$,
$V_r^k \subset H^1(\Omega)$, $V_R^k \subset H^1(\Omega)$, $T_r \subset
H^1(\Omega)$, $T_R \subset H^1(\Omega)$, and $Q \subset H^1(\Omega)$
be spaces of functions for $r = 1, \dots, R-1$ and $k = 1, \dots,
|K|$.  Given $\alpha_r^n$, $[k]_r^n$, and $[k]_R^n$ at time level $n$,
at each time level $n+1$ find the volume fractions $\alpha_r \in S_r$,
the ion concentrations $[k]_r \in V_r^k$, $[k]_R \in V_R^k$, the
potentials $\phi_r \in T_r$, $\phi_R \in T_R$, and the extracellular
mechanical pressure $p_R \in Q$ such that:
\begin{subequations} \label{eq:full:summary}
\begin{align}
    \label{eq:full:summary:alpha:r}
    \frac{1}{\Delta t} \inner{\alpha_r - \alpha_r^n}{s_r}
    - \inner{\alpha_r u_r}{\nabla s_r}
    + \gamma_{rR} \inner{w_{rR}}{s_r}
    &= 0, \\
    \label{eq:full:summary:k:r}
     \frac{1}{\Delta t} \inner{ \alpha_r [k]_r - \alpha_r^n [k]_r^n}{v_r^k}
    - \inner{J_r^k}{\nabla v_r^k}
    +\gamma_{rR}\inner{{J}_{rR}^k}{v_r^k}
    &= 0 ,\\
    \label{eq:full:summary:k:R}
    \frac{1}{\Delta t} \inner{\alpha_R [k]_R - \alpha_R^n [k]_R^n}{v_R^k}
    - \inner{J_R^k}{\nabla v_R^k}
    - \sum_r^{R-1} \gamma_{rR}\inner{{J}_{rR}^k}{v_R^k}
    &= 0, \\
    \label{eq:full:summary:phi:r}
    \inner{\gamma_{rR} C_{rR} \phi_{rR}}{ t_r}- \inner{z^0 F a_r}{
        t_r} - \inner{F\alpha_r \sum_k z^k [k]_r}{ t_r} &= 0, \\
    \label{eq:full:summary:phi:R}
    - \sum_r^{R-1} \inner{\gamma_{rR} C_{rR}\phi_{rR}}{t_R}
    - \inner{z^0 F a_R}{t_R} - \inner{ F\alpha_R \sum_k z^k [k]_R}{
        t_R} &= 0,  \\
    \label{eq:full:summary:p:R}
     -\inner{\sum_r^R \alpha_r u_r}{\nabla q}
    &= 0,
\end{align}
\end{subequations}
for all $s_r \in S_r$ $v_r^k \in V_r^k$, $v_R^k \in V_R^k$, $t_r \in
T_r$, $t_R \in T_R$, $q \in Q$. The compartmental ion flux $J_{r}^{k}$
is given by~\eqref{eq:fluxJ:r}, the compartmental fluid velocity $u_r$
is given by~\eqref{eq:u:r}, the transmembrane water flux $w_{rR}$ is
given by~\eqref{eq:w:r}, while the transmembrane ion fluxes $J^k_{rR}$
are subject to modelling. As before, the potentials $\phi_r$ for $r =
1, \dots, R$ and the extracellular mechanical pressure $p_R$ are only
determined up to a constant. We employ continuous piecewise linear
elements for all variables. Note that $\alpha_r \in L^2(\Omega)$ is
sufficient for the weak formulation of the zero flow limit model to be
well-defined, whereas in the full model, the gradient of $\alpha_r$
appears via the expression for the compartmental fluid velocities
$u_r$ in~\eqref{eq:u:r}, thus suggesting $\alpha_r \in
H^1(\Omega)$. For the ODE time stepping, we apply an ESDIRK4 scheme and
first order Godunov splitting.

 
\section{Numerical convergence study with fluid dynamics: smooth analytical solution}
\label{sec:full:MMS}
To evaluate the numerical accuracy of the scheme presented above
in Section~\ref{sec:full:scheme}, we construct an analytical solution using the
method of manufactured solutions\cite{roache1998verification}. We consider the
full model with three compartments, namely a neuronal ($n$, $r=1$), a glial ($g$, $r=2$)
and an extracellular compartment ($e$, $r=3$).
We use $1,2,3$ and $n,g,e$ interchangeably for subscripts of our variables and
model parameters. In each compartment, we model the movement of potassium
(K$^+$), sodium (Na$^+$) and chloride (Cl$^-$). The transmembrane ion fluxes
$J^{k}_{rR}$ are taken to be passive leak fluxes,~i.e.:
\begin{align}
    J^{k}_{rR} = \frac{1}{F z^k}I_{r,\T{leak}}^k,
\end{align}
where $I_{r,\T{leak}}^k$ is defined by~\eqref{eq:mem:leak} for
compartment $r$ and ion species $k$. The neuronal leak conductances
$g_{n,\T{leak}}^k$ are given in Table~\ref{tab:params:mem:n}, and the glial leak
conductances $g_{g,\T{leak}}^k$ for $\Na^+$ and $\Cl^-$ are given in
Table~\ref{tab:params:full}, whereas the glial $\K^+$ conductance is given by~\eqref{eq:mem:KIR:g}.  We consider a one dimensional domain $\Omega =
[0,1]$ uniformly meshed with $N \in \{8,16,32,64,128\}$ elements. We initially
set $\Delta t = 10^{-3}$ s, and then halve the timestep with each spatial
refinement. The errors are evaluated at $t = 2 \times 10^{-3}$ s. The
analytical solutions are given by~\eqref{eq:MMS:compartmental}--\eqref{eq:MMS:compartmental_2} for the
neuronal and extracellular tissue variables, and by the following for the glial
tissue variables and the extracellular mechanical pressure:
\begin{equation}
\begin{aligned}
    \label{eq:MMS:compartmental:full}
    \alpha_g &= 0.2 - 0.1\sin(2\pi x)\exp(-t), & [\Na]_g &= 0.5 + 0.6\sin(\pi x)\exp(-t), \\
    [\K]_g &= 0.5 + 0.2\sin(\pi x)\exp(-t), & [\Cl]_g &= 1.0 + 0.8\sin(\pi x)\exp(-t), \\
    \phi_g &= \sin(2\pi x)\exp(-t), & p_R &= \sin(2\pi x)\exp(-t).
\end{aligned}
\end{equation}
Parameters values are given in Table~\ref{tab:params:n} and
Table~\ref{tab:params:full}. Initial and boundary
conditions are governed by the exact solutions~\eqref{eq:MMS:compartmental}
and~\eqref{eq:MMS:compartmental:full}.

\begin{table}[ht]
\begin{center}
    \begin{tabular}{lcccc}
        \toprule
        Parameter & Symbol & Value & Unit & Ref. \\
        \midrule
        $\Na^+$ leak conductance glial &$g_{g,\T{leak}}^{\Na}$ & $0.072$ & S/m$^{2}$ & \cite{ostby2009astrocytic}\\
        $\Cl^-$ leak conductance glial &$g_{g,\T{leak}}^{\Cl}$ & $0.5$ & S/m$^{2}$ &\cite{ostby2009astrocytic}\\
        KIR resting conductance glial &$ g_{\KIR}^0$ & $1.3$ & S/m$^{2}$ & \cite{steinberg2005}\\
        Maximum NaKCl rate glial &$g_{\Na\K\Cl}$ & $8.13 \times 10^{-4}$ &
        A/m$^{2}$ & \cite{o2016effects}\\
        Maximum pump rate glial &$\hat{I}_g$ & $0.0372$ & A/m$^{2}$
        &\cite{o2016effects} \\
        Membr. area-to-volume glial  &$\gamma_{ge}$ & $6.3849\times 10^{5}$ & 1/m &
        \cite{kager2000simulated} \\
        Membr. capacitance glial &$C_{ge}$ & $7.5 \times 10^{-3}$ & F/m$^2$  &
        \cite{kager2000simulated} \\
        Membr. water permeability glial &$\eta_{ge}$ &  $5.4\times 10^{-10}$ &
        m$^4$/(mol s) & \cite{o2016effects} \\
        Membrane stiffness neuron & $S_{ne}$ & $2.85 \times 10^{3}$ & Pa/m$^3$ & -- \\
        Membrane stiffness glial & $S_{ge}$ & $2.85 \times 10^{3}$ & Pa/m$^3$ & -- \\
        Gap junction connectivity glial & $\chi_g$ & $0.05 $ &  & -- \\
        Neuronal water permeability & $\kappa_{n}$ & $0$ & m$^4$/(N s) & * \\
        Glial water permeability & $\kappa_{g}$ & $5.0 \times 10^{-16}$ & m$^4$/(N s) & -- \\
        ECS water permeability & $\kappa_{e}$ & $5.0 \times 10^{-16}$ & m$^4$/(N s) & -- \\
        \bottomrule
\end{tabular}
\end{center}
  \caption{Physical parameters for the glial membrane and mechanical pressure.
    We use SI base units, that is, meter (m), mole (mol), Siemens (S) and
    ampere (A). The values are collected from O'Connell et
    al.\protect\cite{o2016effects}, \O{}stby et
    al.\protect\cite{ostby2009astrocytic}, Steinberg et
    al.\protect\cite{steinberg2005}, and Yao et al.\protect\cite{yao2011}. The symbol --
    indicates that the value is chosen by the authors, as we could not find any relevant values in
    the literature. * The neuronal water permeability is set to zero as we
    assume no gap junctions connecting the neurons.}
    \label{tab:params:full}
\end{table}

Based on properties of the approximation spaces and the time discretization,
we expect the optimal rate of convergence to be $1$ in the $H^1$-norm and $2$ in
the $L^2$-norm for the volume fractions $\alpha_n$, $\alpha_g$, $\alpha_e$, the
ion concentrations $[k]_n$, $[k]_g$, $[k]_e$, the potentials $\phi_n$,
$\phi_g$, $\phi_e$ and the mechanical pressure $p_e$. Our numerical observations
are in agreement with the theoretically optimal rates
(Table~\ref{tab:MMS:full}).

\begin{table}[ht]
\begin{center}
    \begin{tabular}{ c  c  c  c  c}
        \toprule
        $N$
        & $\norm{[{\rm Na}]_e - {[{\rm Na}]_e}_h}_{L^2}$
        & $\norm{\phi_n - {\phi_n}_h}_{L^2}$
        & $\norm{\alpha_n - {\alpha_n}_h}_{L^2}$
        & $\norm{p_e - {p_e}_h}_{L^2}$ \\
        \midrule
        8 & 2.47E-03(2.05) & 7.19E-02(1.94) & 1.73E-03(2.14) & 1.31E+01(3.40)\\
        16 & 6.11E-04(2.02) & 1.81E-02(1.99) & 4.13E-04(2.06) & 2.26E+00(2.53)\\
        32 & 1.52E-04(2.00) & 4.54E-03(2.00) & 1.02E-04(2.02) & 4.79E-01(2.24)\\
        64 & 3.80E-05(2.00) & 1.14E-03(2.00) & 2.54E-05(2.00) & 1.15E-01(2.06)\\
        128 & 9.51E-06(2.00) & 3.15E-04(1.85) & 6.35E-06(2.00) & 2.41E-02(2.25)\\
        \midrule
        $n$
        & $\norm{[{\rm Na}]_e - {[{\rm Na}]_e}_h}_{H^1}$
        & $\norm{\phi_n - {\phi_n}_h}_{H^1}$ \\
        \cmidrule{1-3}
        8 & 1.51E-01(1.00) & 1.02E+00(1.01)\\
        16 & 7.55E-02(1.00) & 5.05E-01(1.02)\\
        32 & 3.77E-02(1.00) & 2.52E-01(1.01)\\
        64 & 1.89E-02(1.00) & 1.26E-01(1.00)\\
        128 & 9.42E-03(1.00) & 6.28E-02(1.00)\\
        \cmidrule{1-3}
    \end{tabular}
    \end{center}
 \caption{Selected $L^2$-errors (upper panel) and $H^1$-errors (lower panel)
    and convergence rates (in parenthesis) for the full scheme at time
    $t=0.002$ s.  The test was run on the unit interval, and we initially let
    $\Delta t = 0.001$ s, and then half the timestep in each series. The
    spatial discretization consists of $N$ intervals.}
    \label{tab:MMS:full}
\end{table}

\section{Numerical convergence study: physiological CSD model
with microscopic fluid mechanics}
\label{sec:full:CSD}

\subsection{Problem description} \label{sec:sub:full:model}
Finally, we consider the full model with a neuronal ($n$, $r=1$), glial ($g$, $r=2$)
and extracellular compartment ($e$, $r=3$) and three ion species, namely
potassium (K$^+$), sodium (Na$^+$) and chloride (Cl$^-$). Cortical spreading
depression is triggered by applying excitatory fluxes to the neurons as
described in~\eqref{eq:trigger} in the one-dimensional domain of length $10$ mm. The
physiological parameters values are given in Table~\ref{tab:params:n} and
Table~\ref{tab:params:full}, whereas  the initial conditions are given in
Table~\ref{tab:init}B (Supplementary Tables).

The neuronal membrane mechanisms are as described in
Section~\ref{sec:sub:convergence:model}.  For the glial membrane mechanisms we
follow O'Connell\cite{o2016effects} and consider leak currents modelled as
in~\eqref{eq:mem:leak} (with $r=g$) for sodium (Na$^+$) and chloride (Cl$^-$),
and a potassium inward rectifier current ($I_\T{KIR}$, A/m$^2$).
Following Steinberg et al.\cite{steinberg2005}, the KIR conductance $g_{\KIR}$ (S/m$^2$) is
modelled as:
\begin{align}
    \label{eq:mem:KIR:g}
    g_{\KIR} &= g_{\KIR}^0 \sqrt{\frac{[\K^+]_R}{3}} \frac{1 +
    \exp(\frac{18.5}{42.5})}{1 + \exp(\frac{\phi_{ge} - E^g_{\K} + 18.5}{42.5})}
     \frac{1 + \exp(\frac{-118.6 - 85.2}{44.1})}
    {1 + \exp(\frac{-118.6 + \phi_{ge}}{44.1})},
\end{align}
where $g_{\KIR}^0$ (S/m$^{2}$) is the resting membrane conductance,
and corresponds to the conductance when $\phi_{ge} = E^g_{\K}$ and
$[{\rm K}^+]_e = [{\rm K}^+]_e^0$.
The Na/K/ATPase pump (ATP) occurs in both the neuron and glial membrane, and is
modelled as in~\eqref{eq:mem:ATP}. Finally, the current through the NaKCl
cotransporter $I_{\Na\K\Cl}$ (A/m$^{2}$) is modelled as:
\begin{equation}
    I_{\Na\K\Cl} = g_{\Na\K\Cl} \ln{\left(\frac{[\Na^+]_g [\K^+]_g
    [\Cl^-]_g^2}{[\Na^+]_e [\K^+]_e [\Cl^-]_e^2}\right)}.
\end{equation}

In summary, the total currents over the glial membrane are
modelled by~\eqref{eq:totmem:fluxes:g} (with the currents
(A/m$^{2}$) are converted to ion fluxes (mol/(m$^{2}$s)) by dividing by
Faraday's constant $F$ times the valence $z^k$):
\begin{subequations} \label{eq:totmem:fluxes:g}
\begin{align}
    J_g^{\Na} &= \frac{1}{F z^{\Na}}\left(I_{g,\T{leak}}^{\Na} + 3 I_{g,\ATP} + I_{\Na\K\Cl} + I^{\Na}_{\rm ex} \right),\\
    J_g^{\K} &= \frac{1}{F z^{\K}}\left(I_{\KIR} - 2 I_{g,\ATP} + I_{\Na\K\Cl} + I^{\K}_{\rm ex} \right), \\
    J_g^{\Cl} &= \frac{1}{F z^{\Cl}}\left(I_{g,\T{leak}}^{\Cl} + 2I_{\Na\K\Cl}
    + I^{\Cl}_{\rm ex} \right).
\end{align}
\end{subequations}
Here, $I^{\Na}_{\rm ex}$, $I^{\K}_{\rm ex}$, and $I^{\Cl}_{\rm ex}$
are excitatory fluxes used to trigger a cortical spreading depression
wave, see~\eqref{eq:trigger} in Section~\ref{sec_init_CSD_wave}.

\subsection{CSD wave characteristics}
\label{sec:full:CSD:characteristics}

\begin{figure}
\centering
\includegraphics[width=1.0\linewidth]{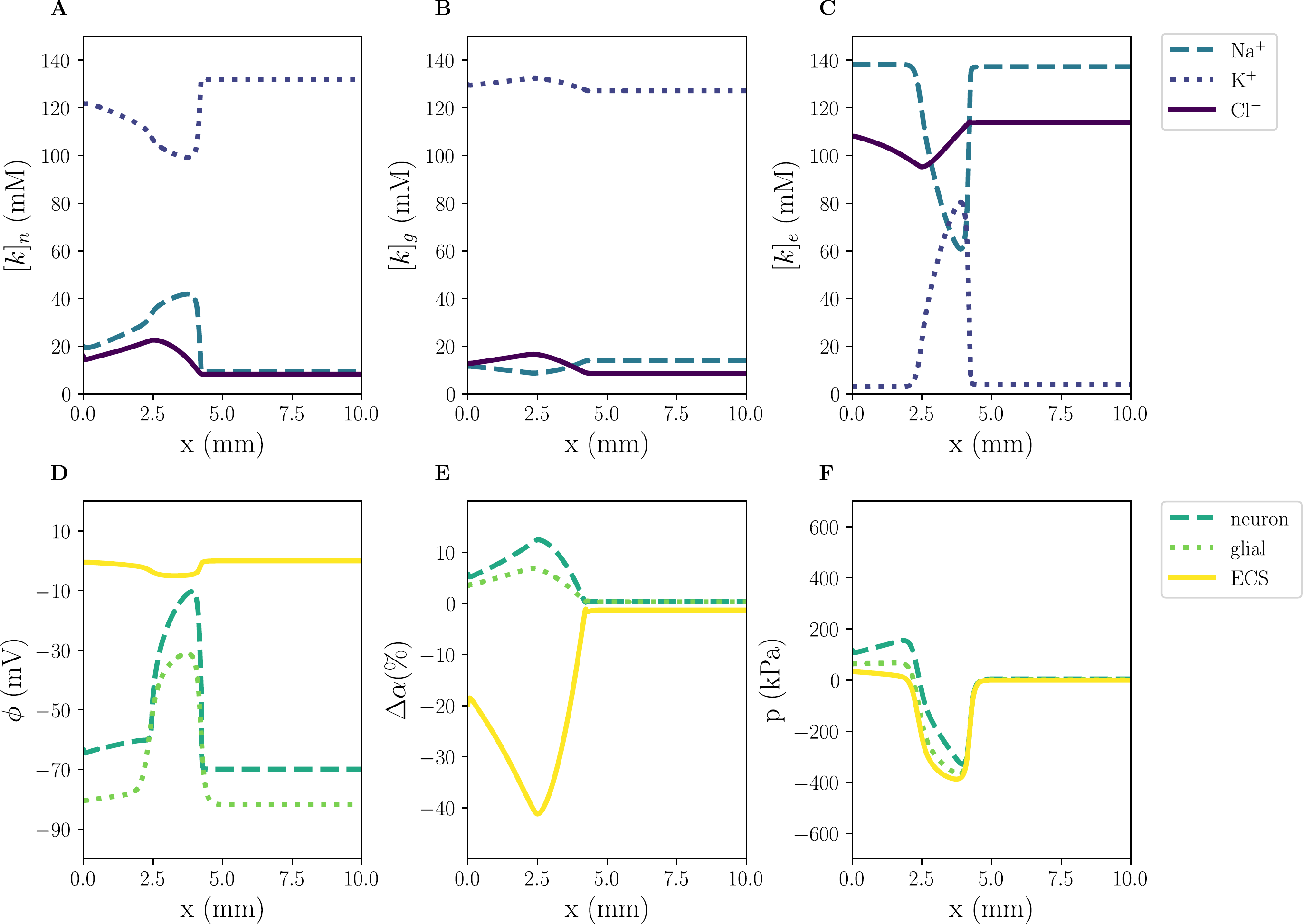}
\caption{Full model simulation of a CSD wave. Snapshots at $t=50$ s of neuronal, glial
  and extracellular ion concentrations (\textbf{A}, \textbf{B}, \textbf{C}),
  electrical potentials (\textbf{D}), change in volume fractions (\textbf{E})
  and mechanical pressure (\textbf{F}).} \label{fig:full:main}
\centering
\includegraphics[width=1.0\linewidth]{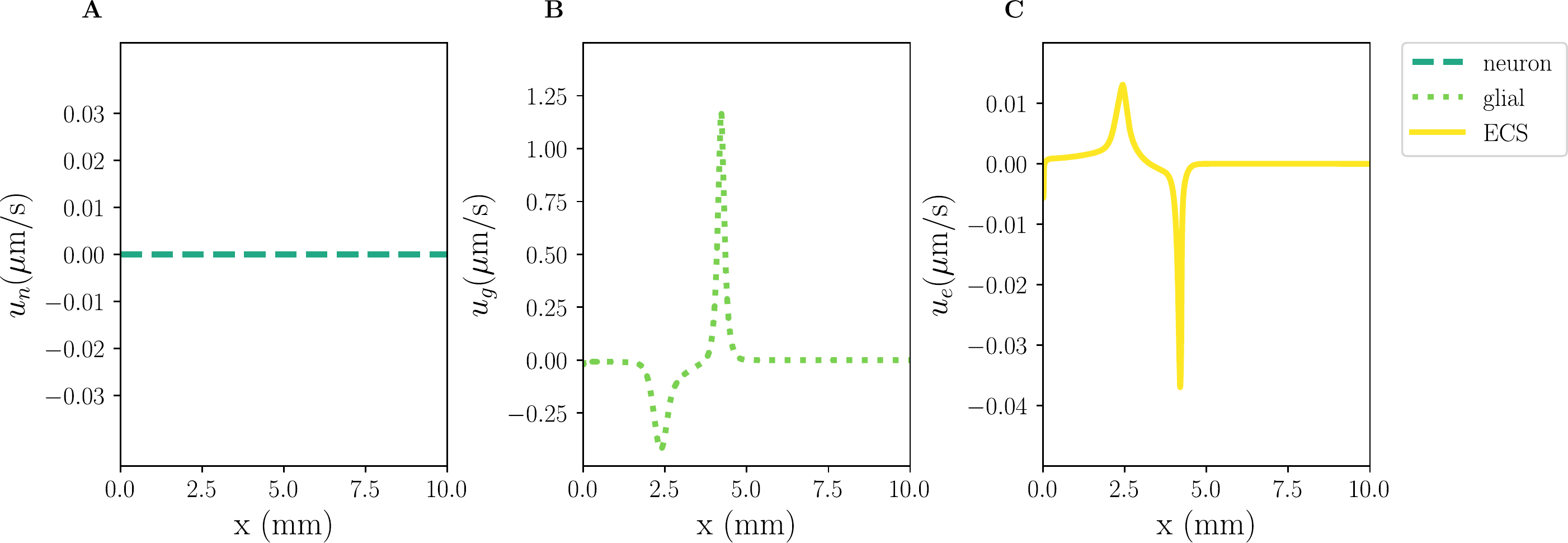}
    \caption{Full model simulation of a CSD wave: fluid velocities. Snapshots at $t = 50$ s of compartmental fluid
    velocities in the neurons (\textbf{A}), the glial cells (\textbf{B}), and
    the ECS (\textbf{C}).}
    \label{fig:full:velocity}
\end{figure}

As in the zero flow limit scenario, excitatory flux stimulation leads
to a CSD wave traveling through the tissue: we observe neuronal
depolarization, neuronal and ECS ionic concentration changes, and
neuronal swelling (Figure~\ref{fig:full:main}). Moreover, we observe
that the glial potential depolarizes: from $-81$ to $-31$ mV,
accompanied by a small drop in the extracellular potential from $0$ to
$-5$ mV. Substantial alterations in the intra- and extracellular ion
compositions follow the depolarization wave: we observe an
increase in the concentrations of extracellular potassium, and
decreases in extracellular sodium and chloride. In response to the
ionic shifts, the neurons and glial cells swell with an increase in
volume fractions of respectively $12.5\%$ and $6.8\%$, while the
extracellular space shrinks correspondingly. We note that the neural
and extracellular dynamics are qualitatively similar to those in the
case of two compartments (neurons, ECS) in the zero flow limit
(c.f.~Figure~\ref{fig:zerolimit:wave}), which is in accordance with
the (numerical) findings reported by O'Connell et
al\cite{o2016effects}.

The CSD wave is accompanied by a decrease in the mechanical pressures,
from the baseline of 0 kPa down to $-288$, $-334$ and $-389$ kPa in
the neuronal, glial and extracellular compartments, respectively
(Figure~\ref{fig:full:main} F). The pressure gradients (mechanical and
osmotic) drive microscopic fluid flow within the glial and the
extracellular compartments. We observe flow rates of up to $1.1$ and
$-0.03$ $\mu$m/s in the glial and extracellular compartments,
respectively; i.e. fluid flows in opposite directions
(Figure~\ref{fig:full:velocity} B,C). During glial swelling, water
moves across the glial membrane from the ECS into the glial cells. In
response, water within the glial cell network will be pushed away from
this area. Indeed, we observe a positive flow rate to the right of the
swelling, and a negative flow rate left of the swelling in the glial
cells and vice versa in the ECS. We observe no flow in the neuronal
compartment, which is expected as the neuronal water permeability is
set to zero (Figure~\ref{fig:full:velocity}A).

Next, we (quantitatively) evaluate the numerical convergence of the full model
scheme presented in Section~\ref{sec:full:scheme}. To this end, we consider the
quantities of interest defined in Section~\ref{sec_CSD_wave_charac} in addition
to the (spatial) width $w_{\rm CSD, p}$ of the extracellular mechanical pressure wave at
$t = 50$ s, given by
    \begin{align*}
        \mathbb{X} &= \{x \ | \ p_R(x, 50) > p_{\rm thres}\} \\
        w_{\rm CSD,p} &= \max\mathbb{X} - \min\mathbb{X}
    \end{align*}
    with $p_{\rm thres} = -10$ kPa.

\subsection{Convergence of CSD wave characteristics during refinement}

\begin{figure}
        \textbf{A}
        \vspace{0.25cm}
        \begin{center}
        \begin{overpic}[width=0.13\linewidth]{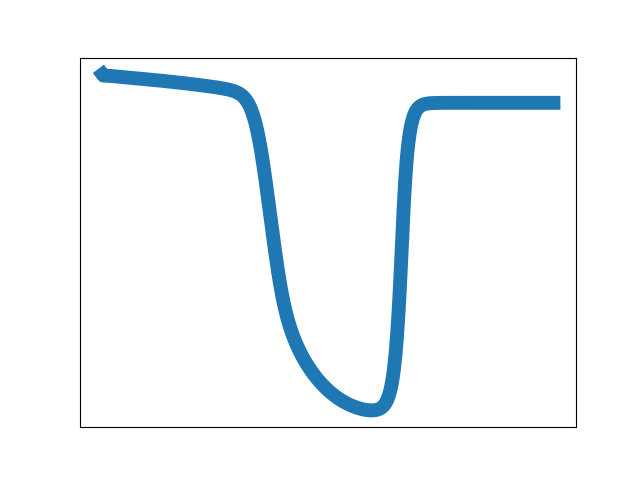}
            \put(-3,82){\color{black}\vector(1,0){50}}
            \put(55,81){Temporal refinement}
            \put(-18,-140){\rotatebox{90}{Spatial refinement}}
            \put(-3,82){\color{black}\vector(0,-1){50}}
        \end{overpic}
        \includegraphics[width=0.13\linewidth]{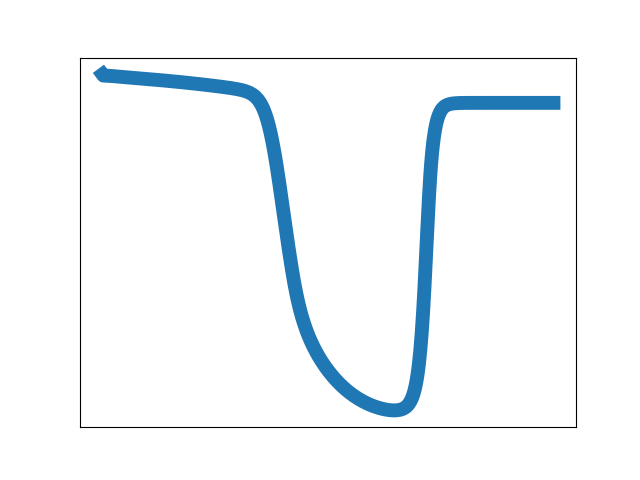}
        \includegraphics[width=0.13\linewidth]{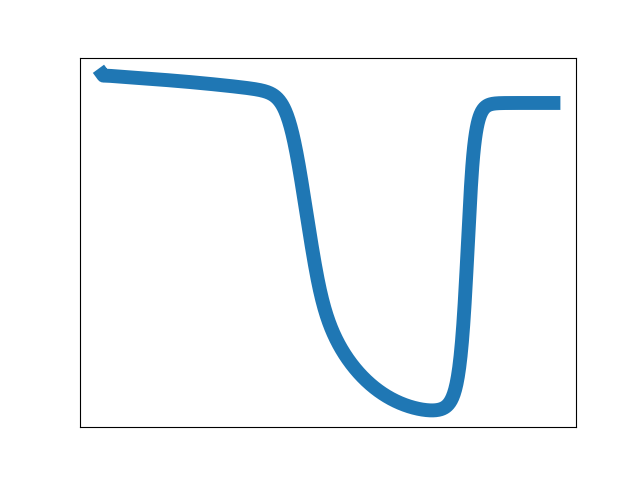}
        \includegraphics[width=0.13\linewidth]{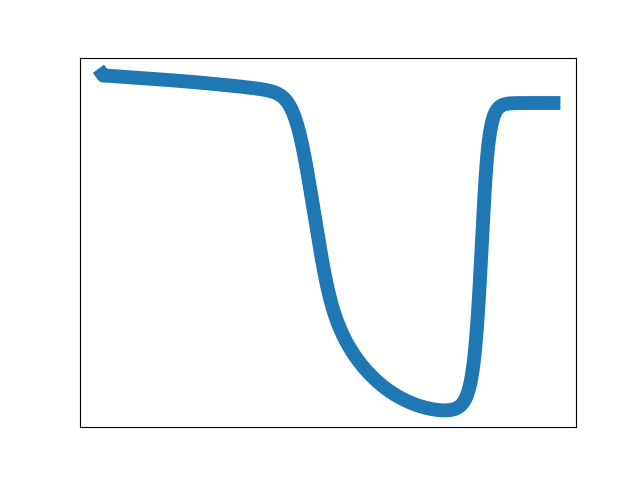}
        \includegraphics[width=0.13\linewidth]{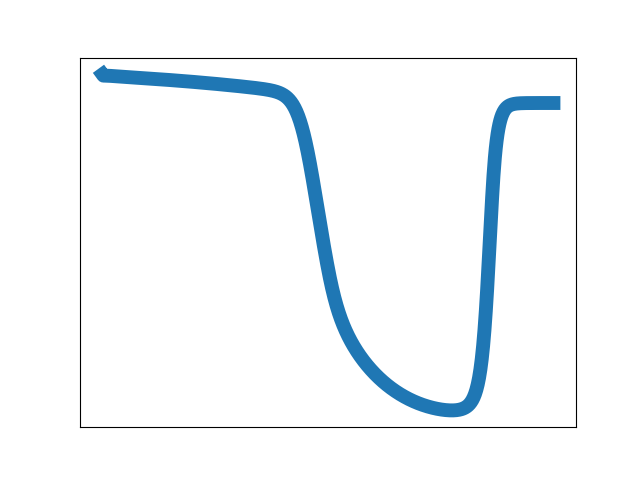}
        \includegraphics[width=0.13\linewidth]{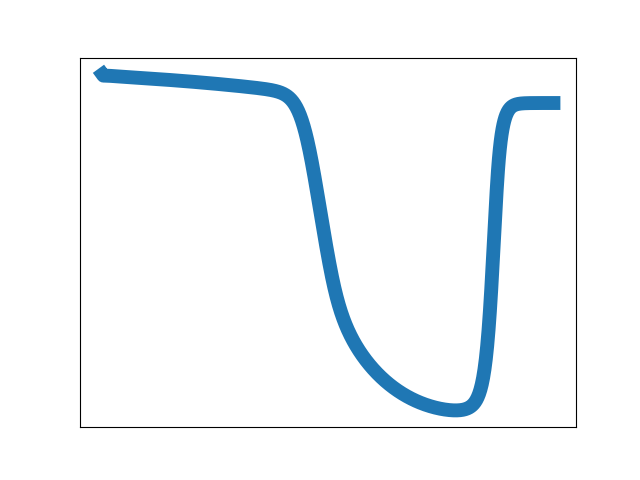}
        \includegraphics[width=0.13\linewidth]{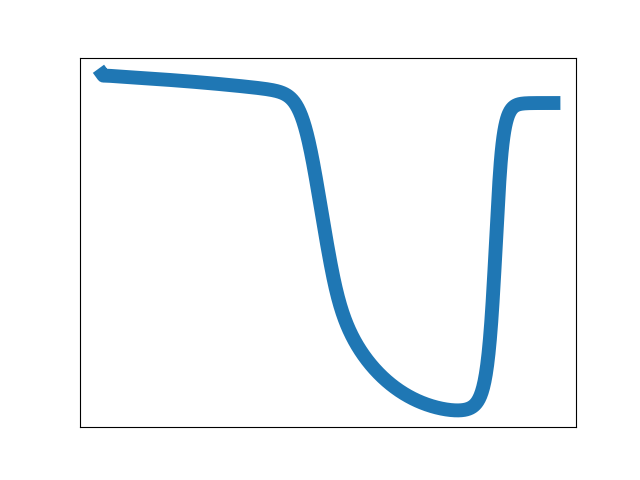} \\
        \includegraphics[width=0.13\linewidth]{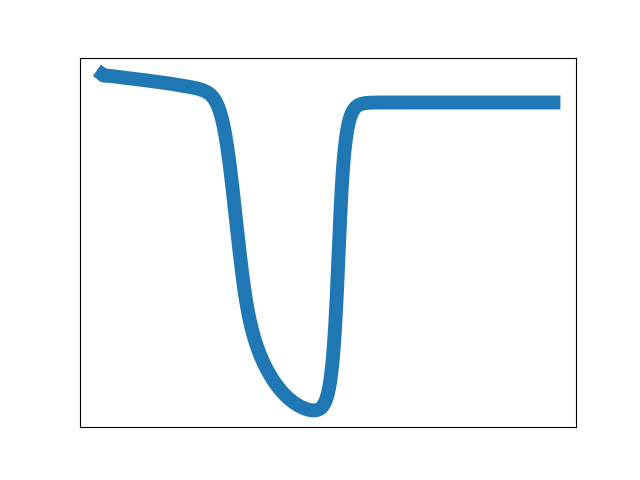}
        \includegraphics[width=0.13\linewidth]{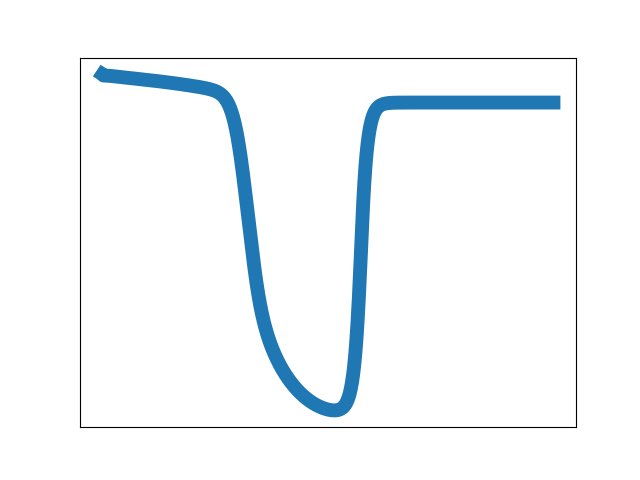}
        \includegraphics[width=0.13\linewidth]{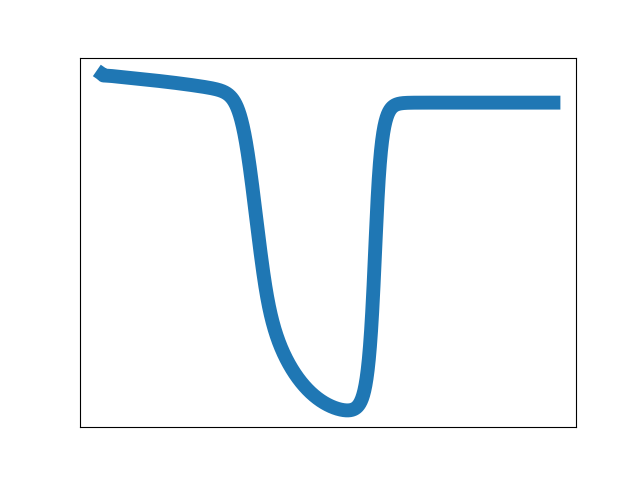}
        \includegraphics[width=0.13\linewidth]{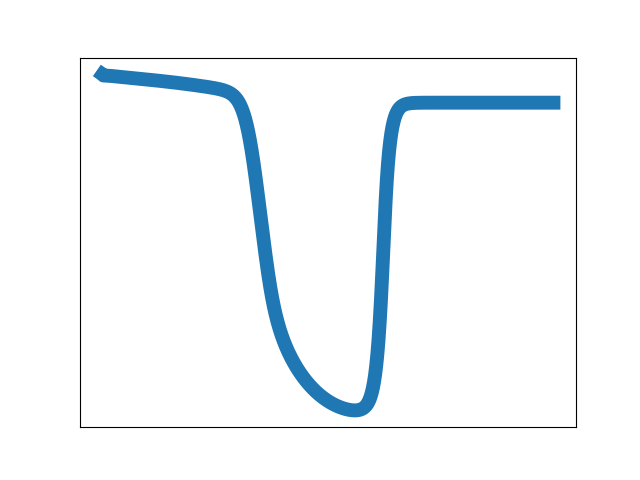}
        \includegraphics[width=0.13\linewidth]{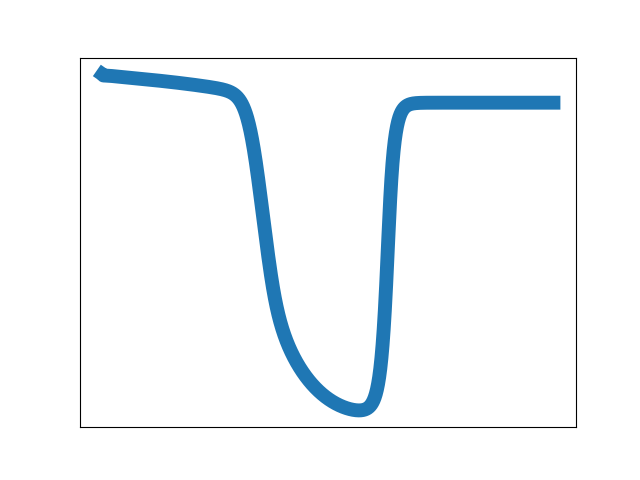}
        \includegraphics[width=0.13\linewidth]{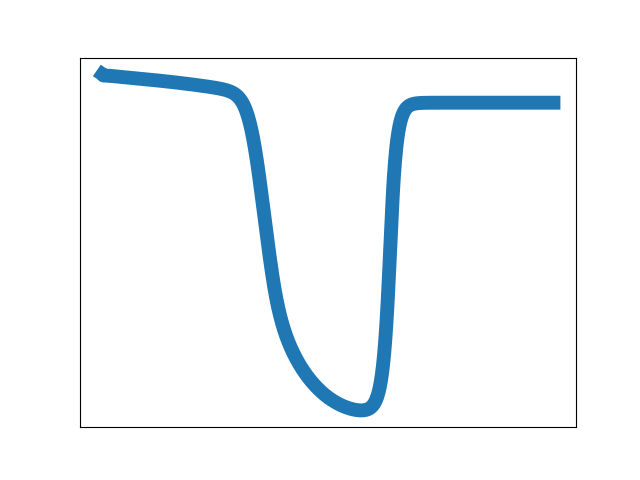}
        \includegraphics[width=0.13\linewidth]{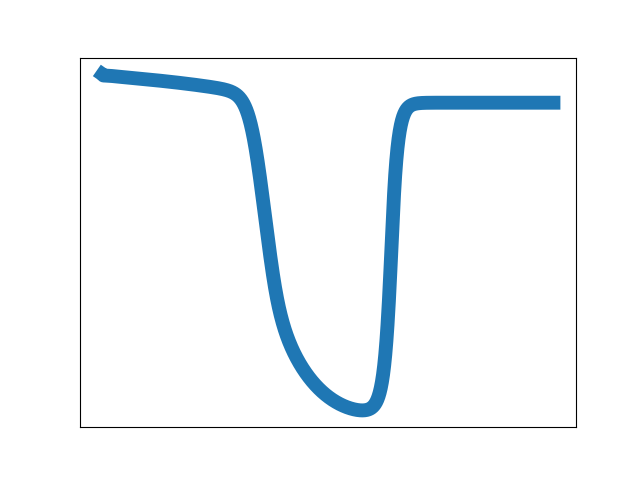} \\
        \includegraphics[width=0.13\linewidth]{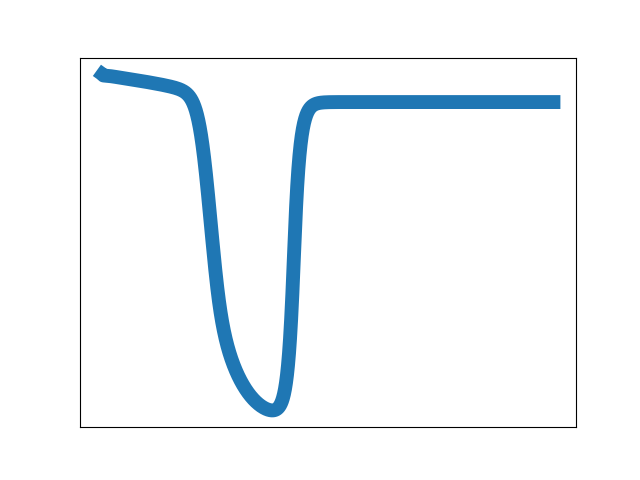}
        \includegraphics[width=0.13\linewidth]{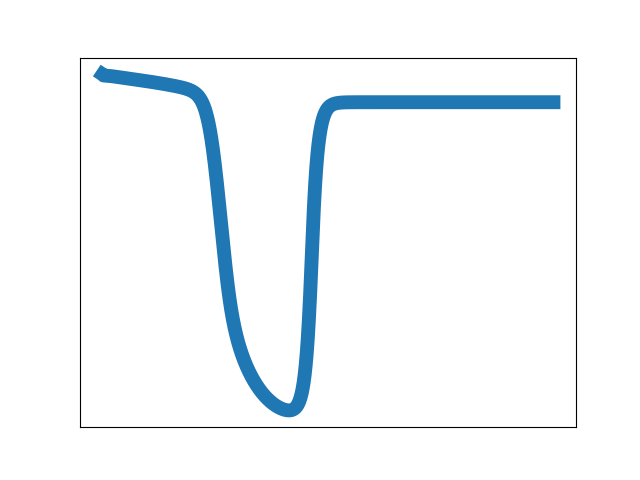}
        \includegraphics[width=0.13\linewidth]{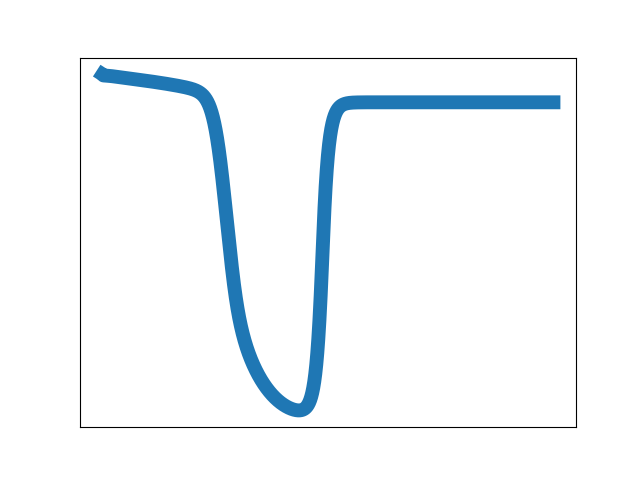}
        \includegraphics[width=0.13\linewidth]{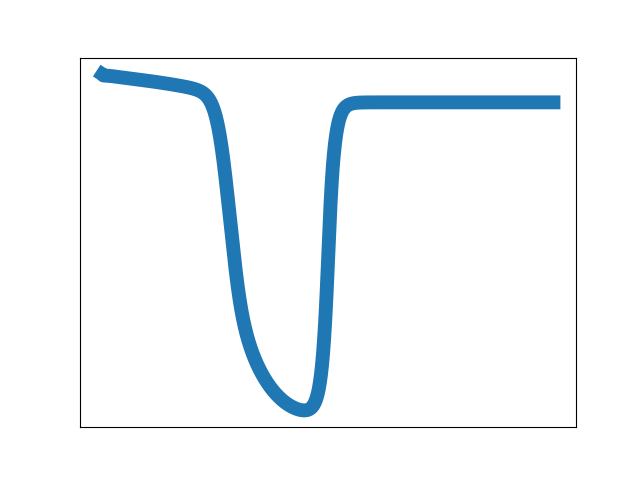}
        \includegraphics[width=0.13\linewidth]{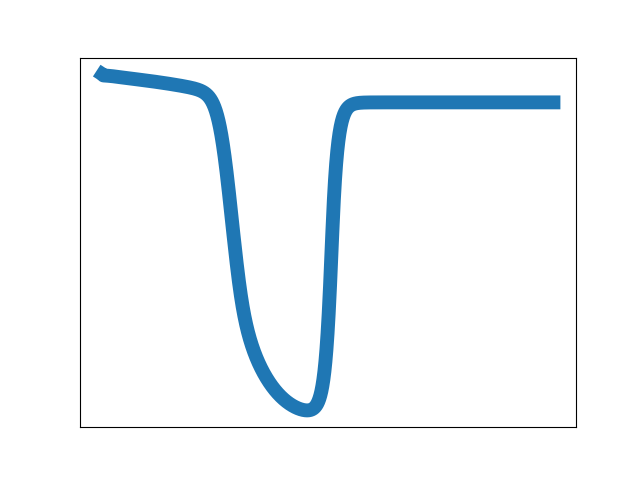}
        \includegraphics[width=0.13\linewidth]{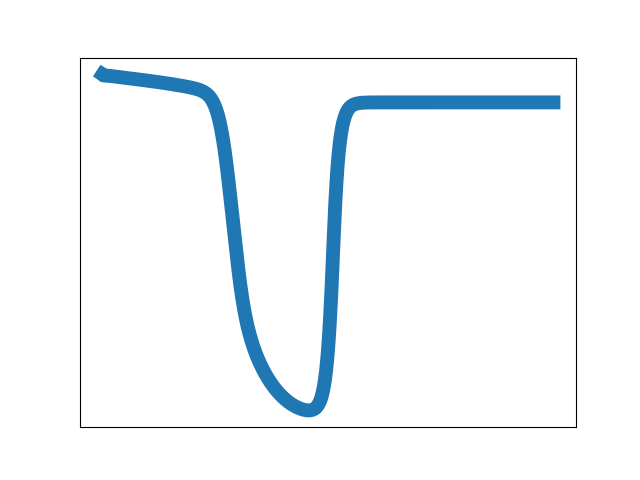}
        \includegraphics[width=0.13\linewidth]{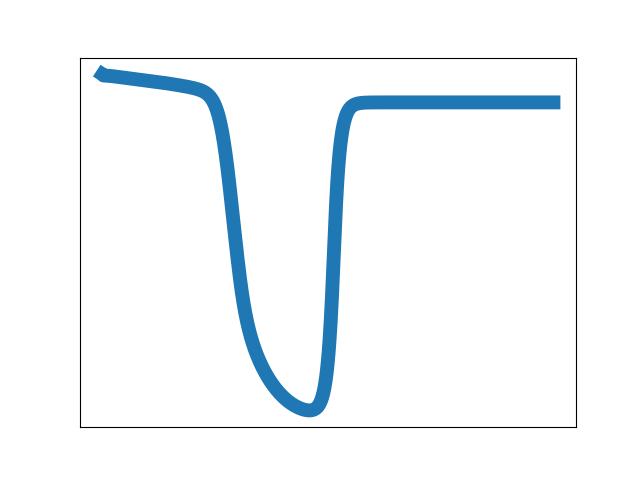} \\
        \includegraphics[width=0.13\linewidth]{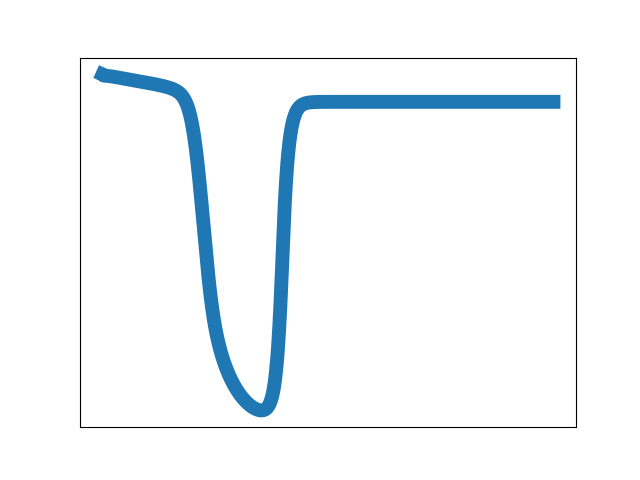}
        \includegraphics[width=0.13\linewidth]{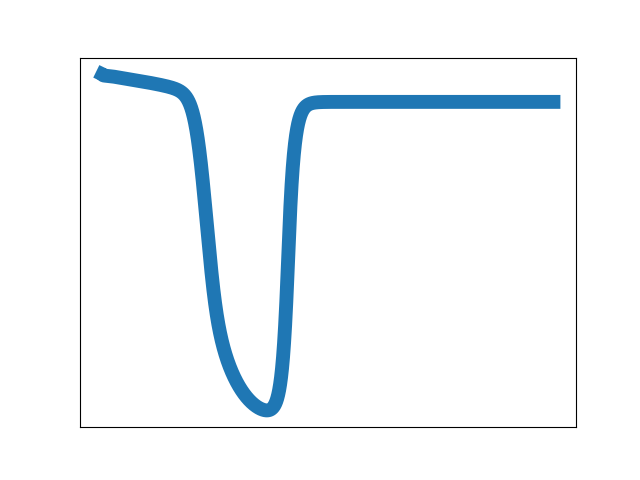}
        \includegraphics[width=0.13\linewidth]{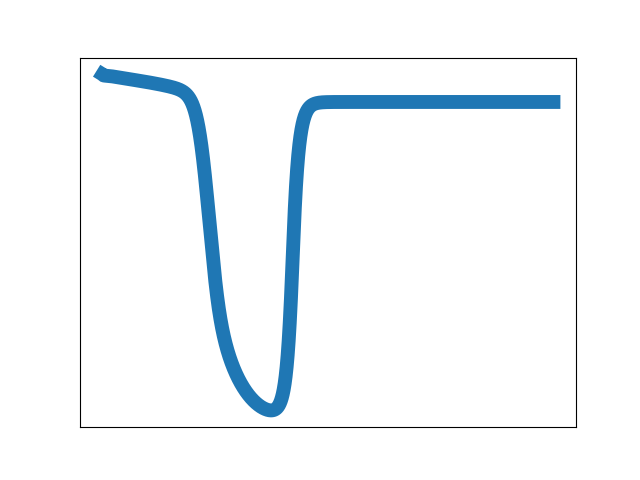}
        \includegraphics[width=0.13\linewidth]{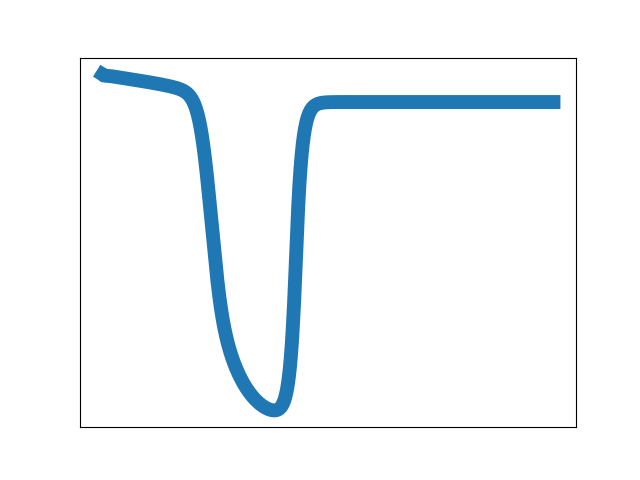}
        \includegraphics[width=0.13\linewidth]{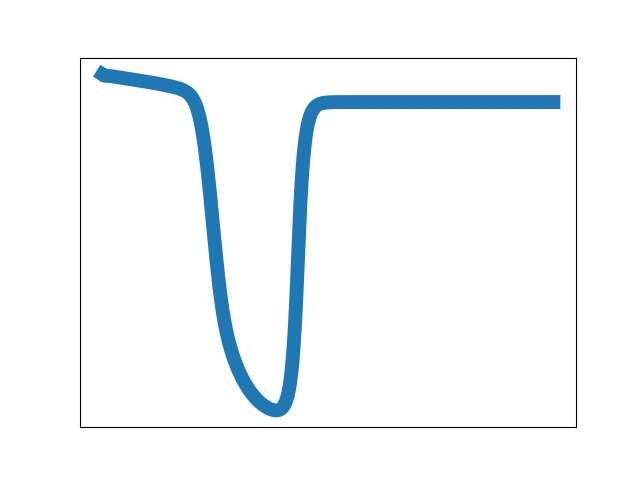}
        \includegraphics[width=0.13\linewidth]{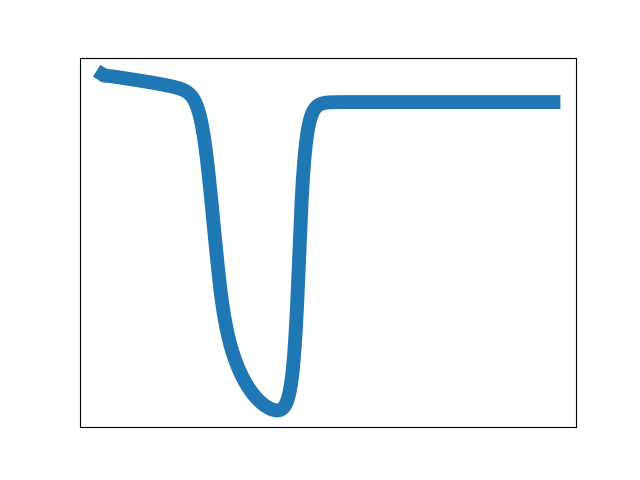}
        \includegraphics[width=0.13\linewidth]{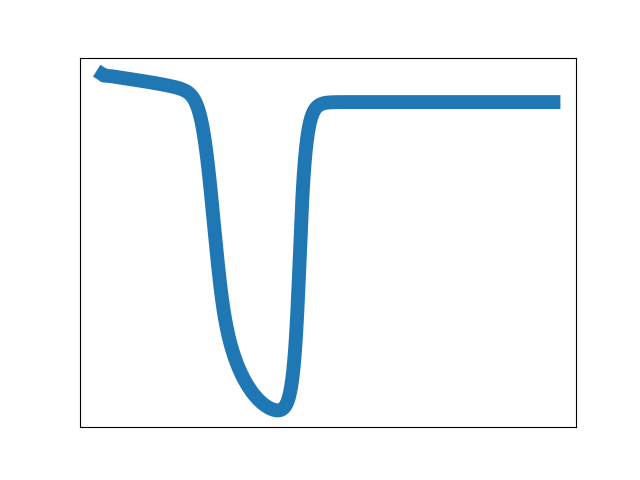}
    \end{center}
    \textbf{B} 
    \begin{center}
    \begin{tabular}{ c | c c c c c c c | c}
        \toprule
        \diagbox{$N$}{$\Delta t$} & 12.5 & 6.25 & 3.125 & 1.563 & 0.781 & 0.391
        & 0.195 & $\Delta w_{\rm CSD,p}$ \\
        \midrule 
        1000 &3.640 &3.910 &4.350 &4.520 &4.600 &4.650 &4.670 &\\
        2000 &2.960 &3.200 &3.350 &3.440 &3.485 &3.510 &3.525 &1.145\\
        4000 &2.502 &2.687 &2.795 &2.857 &2.887 &2.907 &2.915 &0.610\\
        8000 &2.409 &2.464 &2.505 &2.536 &2.555 &2.566 &2.571 &0.344\\
        \midrule 
        $\Delta w_{\rm CSD,p}$ & &0.055 &0.041 &0.031 &0.019 &0.011 &0.005 &\\
        \bottomrule
    \end{tabular}
    \end{center}
    \textbf{C}
    \begin{center}
    \begin{tabular}{ c | c c c c c c c | c}
        \toprule
        \diagbox{$N$}{$\Delta t$} & 12.5 & 6.25 & 3.125 & 1.563 & 0.781 & 0.391 & 0.195 & $\Delta \bar{v}_{\rm CSD}$ \\
        \midrule
        1000 &7.873 &8.528 &9.429 &10.000 &10.250 &10.250 &10.500 &\\
        2000 &6.256 &6.667 &7.000 &7.305 &7.516 &7.516 &7.516 &2.984\\
        4000 &4.977 &5.435 &5.654 &5.821 &5.987 &5.987 &5.987 &1.529\\
        8000 &4.583 &4.667 &4.977 &5.091 &5.091 &5.049 &5.162 &0.825\\
        \midrule
        $\Delta \bar{v}_{\rm CSD}$ & &0.084 &0.310 &0.114 &0.000 &0.042 &0.113 &\\
        \bottomrule
    \end{tabular}
    \end{center}
    \caption{Full model simulation: CSD wave properties during refinement
      in space (N) and time ($\Delta$t, ms) in a 1D domain of length
      $10$ mm at $t=50$ s. \textbf{A}: Extracellular mechanical
      pressure $p_R(x, 50)$ (kPa) versus $x \in \Omega$ (mm).
      \textbf{B}: Extracellular mechanical pressure wave width (mm)
      and difference $\Delta w_{\rm CSD,p}$ between consecutive
      refinements.  \textbf{C}: CSD mean wave speed $\bar{v}_{\rm
        CSD}$ (mm/min) and difference $\Delta \bar{v}_{\rm CSD}$
      between consecutive refinements.}
    \label{fig:convergence:CSD:full}
\end{figure}

Drawing on the findings in Section~\ref{sec:full:scheme}, here, we
consider the implicit lower-order scheme based on Godunov splitting, a
BE method for the PDE time-discretization, and ESDIRK4 for the ODE
time-stepping -- to compute the mean speed (cf.
Section~\ref{sec_CSD_wave_charac}) and the (spatial) width of the
extracellular mechanical pressure wave
(cf.~Section~\ref{sec:full:CSD:characteristics}) for different mesh
resolutions and time steps: $\Delta x_N = 10/N$ mm for $N = 1000,
2000, 4000, 8000, 16000, 32000$ and $\Delta t_i = 12.5/i$ s for $i =
1, 2, 4, 8, 16, 32, 64$. The results are presented in
Figure~\ref{fig:convergence:CSD:full}. Wave speeds are converted from
the native m/s to mm/min for interpretability.

As in the zero flow limit, the computed mean wave speed and the extracellular
mechanical pressure wave width increases with decreasing $\Delta t$, and
decreases with decreasing $\Delta x$: the smaller the time step, the faster and
wider the wave, while the smaller the mesh size, the slower and narrower the
wave (Figure~\ref{fig:convergence:CSD:full}A--C). The behavior of the mean
wave speed is qualitatively similar to that in the schemes for the zero flow
limit (cf.~\ref{fig:intro}). The computed extracellular mechanical pressure
wave width vary substantially, ranging from $2.409$ to $4.67$ mm.

The differences $\Delta w_{\rm CSD,p}$ in the extracellular mechanical
pressure wave width between the coarsest mesh sizes $N = 1000$ and
$2000$ are in the range $2.96-4.67$ mm, while $\Delta w_{\rm CSD,p}$
between the coarsest time steps $\Delta t = 12.5$ and $6.25$ are in
the range $2.409-3.910$ mm. For the finest time and mesh resolutions,
we observe that $\Delta w_{\rm CSD,p} = 0.005$ and $0.344$ mm,
respectively; thus the spatial error continues to dominate. Finally,
we observe that $\Delta w_{\rm CSD,p}$ decreases as we refine the
discretizations in time and space. There is however no clear
convergence rate (Figure~\ref{fig:convergence:CSD:full}A, B).

\section{Discussion and concluding remarks}
\label{sec:conclusions}
We have presented and analyzed finite element-based splitting schemes
for a mathematical framework modelling ionic electrodiffusion and
water movement in biological tissue in general, and brain
tissue in particular. We have evaluated the schemes in terms of their
numerical properties, including accuracy, convergence, and
computational efficiency, for idealized scenarios as well as for
challenging, physiologically realistic problem settings.

The schemes display optimal convergence rates in space for problems with smooth
manufactured solutions. However, the physiological CSD setting is
challenging: we find that the accurate computation of CSD wave
characteristics (wave speed and wave width) requires a very fine spatial
and fine temporal resolution for all schemes tested. Indeed, different
splitting and time stepping schemes and lower and higher order finite element
schemes give comparable results in terms of accuracy. Overall, the error
associated with the spatial discretization dominates. Explicit PDE and/or
ODE time stepping schemes easily fail to converge even for only moderately
coarse time steps, but yield accurate results for very fine timesteps. In light
of the long time scale associated with CSD (seconds to minutes), the small time
steps imposed by the explicit schemes (less than a millisecond) represent a
severe restriction.

The mathematical framework studied here was presented by Mori in
2015\cite{mori2015multidomain}, and has been used to simulate cortical
spreading depression in a three-compartment setting (including neurons, glial
cells and extracellular space)\cite{o2016effects}, and in multiple spatial
dimensions\cite{tuttle2019computational}. However, little has been reported on
numerical properties of discretizations of this model. The aforementioned
studies\cite{o2016effects,tuttle2019computational} have used time steps of
the order 10 ms and mesh sizes of the order 0.156--0.02 mm (corresponding to $N
= 64-500$ cf.~Figure~\ref{fig:intro}). Our findings indicate that high resolution is
required to accurately compute CSD wave properties and that low-to-moderate
resolutions can substantially overestimate (or, but more rarely, underestimate) the
CSD wave speed. We expect our finite element findings to extend also to
comparable finite difference or finite volume discretizations.

In terms of limitations, we have here compared different numerical
schemes in terms of accuracy, with less emphasis on computational
complexity or cost. We consider these numerical investigatons as a
starting point and guide for future theoretical studies. Another
research direction would be the extension of this study to the
two-dimensional CSD model studied by
O'Connell\cite{o2016computational} and
Tuttle\cite{tuttle2019modeling}, where the development of multigrid
solvers seems crucial to reach the high spatial resolution needed to
obtain accurate solutions.

This paper focuses on numerical challenges related to approximating
systems for ionic electrodiffusion and microscopic water movement. We
remark that the full model simulation, including extracellular
mechanical pressure, yields pressure differences far greater than what
one might expect in this setting ($\sim 6$ times atmospheric
pressure). It seems natural to reevaluate whether the current
compartmental fluid velocity model best represents the physiology, in
particular the fluid velocity component driven by electrostatic
forces. On the other hand, it is well-established that large osmotic
pressure gradients indeed are present in the brain environment. 

In conclusion, our findings show that numerical simulation of ionic
electrodiffusion and water movement in brain tissue is feasible, but
requires care numerically and substantial computational
resources. Numerical schemes or solution approaches that retain
accuracy at a lower computational expense would enable the study of a
wide array of phenomena in brain physiology, including in the context
of pathological conditions.

\appendix
\section{Supplementary Tables}
\label{sec:supp:tab}
\begin{table}[ht]
  \begin{center}
    \begin{tabular}{ c | c c c c c c c | c}
        \toprule
        \diagbox{$N$}{$\Delta t$} & 12.5 & 6.25 & 3.125 & 1.563 & 0.781 & 0.391 & 0.195 & $\Delta \bar{v}_{\rm CSD}$ \\
        \midrule
        1000 &8.000 &8.631 &9.385 &9.738 &9.938 &10.046 &10.092 & -- \\
        2000 &6.262 &6.862 &7.223 &7.431 &7.538 &7.600 &7.638 &2.454\\
        4000 &5.138 &5.636 &5.938 &6.096 &6.181 &6.227 &6.246 &1.392\\
        8000 &4.798 &4.978 &5.128 &5.242 &5.311 &5.349 &5.366 &0.880\\
        \midrule
        $\Delta \bar{v}_{\rm CSD}$ & -- &0.180 &0.150 &0.114 &0.069 &0.038 &0.017 &\\
        \bottomrule
    \end{tabular}
  \end{center}
    \caption{CSD mean wave speed $\bar{v}_{\rm CSD}$ (mm/min) and difference in
        CSD mean wave speed $\Delta \bar{v}_{\rm CSD}$ between consecutive
        refinements in space (rows) and time (columns). Numerical scheme: Godunov
        splitting, BDF2, ESDIRK4.}
  \label{tab:godunov}
\end{table}
\begin{table}
    \textbf{A} 
    \begin{center}
    \begin{tabular}{ c | c c c c c c c | c}
        \toprule
        \diagbox{$N$}{$\Delta t$} & 12.5 & 6.25 & 3.125 & 1.563 & 0.781 & 0.391 & 0.195 & $\Delta \bar{v}_{\rm CSD}$ \\
        \midrule
        1000 & $*$ & $*$ & $*$ & $*$ & 9.938 &10.046 &10.092 & -- \\
        2000 & $*$ & $*$ & $*$ & $*$ & 7.538 &7.600 &7.638 &2.454\\
        4000 & $*$ & $*$ & $*$ & $*$ & 6.181 &6.227 &6.246 &1.392\\
        8000 & $*$ & $*$ & $*$ & $*$ & 5.311 &5.349 &5.366 &0.880\\
        \midrule
        $\Delta \bar{v}_{\rm CSD}$ & -- & -- & -- & -- &  -- & 0.038 &0.017 &\\
        \bottomrule
    \end{tabular}
    \end{center}
    \textbf{B} 
    \begin{center}
    \begin{tabular}{ c | c c c c c c c | c}
        \toprule
        \diagbox{$N$}{$\Delta t$} & 12.5 & 6.25 & 3.125 & 1.563 & 0.781 & 0.391 & 0.195 & $\Delta \bar{v}_{\rm CSD}$ \\
        \midrule
        1000 &7.323 &8.338 &9.077 &9.554 &9.831 &9.985 &10.062 & -- \\
        2000 &5.938 &6.608 &7.054 &7.338 &7.492 &7.585 &7.623 &2.439\\
        4000 &4.970 &5.487 &5.835 &6.038 &6.150 &6.212 &6.238 &1.385\\
        8000 &4.727 &4.922 &5.080 &5.212 &5.293 &5.339 &5.361 &0.877\\
        \midrule
        $\Delta \bar{v}_{\rm CSD}$ &  -- &0.195 &0.158 &0.132 &0.081 &0.046 &0.022 &\\
        \bottomrule
    \end{tabular}
  \end{center}
    \caption{CSD mean wave speed $\bar{v}_{\rm CSD}$ (mm/min) and difference in
        CSD mean wave speed $\Delta \bar{v}_{\rm CSD}$ between consecutive
        refinements in space (rows) and time (columns). Numerical scheme: Strang,
        BDF2, and RK4 (\textbf{A}) or BE (\textbf{B}). $*$ indicates that the solver
        failed to converge.}
    \label{tab:ode}
\end{table}
\begin{table}
    \textbf{A}
    \begin{center}
    \begin{tabular}{  l  c  c  c }
        \toprule
        Parameter & Symbol & Value & Unit \\
        \midrule
        neuron volume fraction & $\alpha_{\rm n}^0$ & $0.8$ &  \\
        $\Na^+$ concentration neuron& $[\Na]_{\rm n}^0$ & $9.3$&  mol/m$^{3}$\\
        $\K^+$ concentration neuron&  $[\K]_{\rm n}^0$  & $132$& mol/m$^{3}$   \\
        $\Cl^-$ concentration neuron& $[\Cl]_{\rm n}^0$ & $8.0$& mol/m$^{3}$ \\
        $\Na^+$ concentration ECS&    $[\Na]_{\rm e}^0$ & $137$& mol/m$^{3}$   \\
        $\K^+$ concentration ECS&     $[\K]_{\rm e}^0$  & $4$& mol/m$^{3}$     \\
        $\Cl^-$ concentration ECS&    $[\Cl]_{\rm e}^0$ & $114 $& mol/m$^{3}$  \\
        potential neuron &            $\phi_{\rm n}^0$ & $-0.070 $& V \\
        potential ECS &               $\phi_{\rm e}^0$ & $0.0 $& V \\
        \bottomrule
    \end{tabular}
    \end{center}
    \textbf{B} 
    \begin{center}
    \begin{tabular}{  l  c  c  c }
        \toprule
        Parameter & Symbol & Value & Unit \\
        \midrule
        neuron volume fraction & $\alpha_{n}^0$ & $0.5$ & \\
        glial volume fraction& $\alpha_{g}^0$ & $0.3$ &  \\
        $\Na^+$ concentration neuron& $[\Na]_{n}^0$ & $9.3$&  mol/m$^{3}$ \\
        $\K^+$ concentration neuron&  $[\K]_{n}^0$  & $132$& mol/m$^{3}$    \\
        $\Cl^-$ concentration neuron& $[\Cl]_{n}^0$ & $8.0$& mol/m$^{3}$  \\
        $\Na^+$ concentration glial& $[\Na]_{g}^0$ & $13$&  mol/m$^{3}$   \\
        $\K^+$ concentration glial&  $[\K]_{g}^0$  & $128$& mol/m$^{3}$    \\
        $\Cl^-$ concentration glial& $[\Cl]_{g}^0$ & $8.0$& mol/m$^{3}$   \\
        $\Na^+$ concentration ECS&   $[\Na]_{e}^0$ & $137$& mol/m$^{3}$    \\
        $\K^+$ concentration ECS&    $[\K]_{e}^0$  & $4$& mol/m$^{3}$      \\
        $\Cl^-$ concentration ECS&   $[\Cl]_{e}^0$ & $114 $& mol/m$^{3}$   \\
        potential neuron &   $\phi_{n}^0$ & $-0.070 $& V \\
        potential glial &   $\phi_{g}^0$ & $-0.082 $& V \\
        potential ECS &   $\phi_{e}^0$ & $0.0 $& V \\
        mechanical pressure ECS &   $p_{e}^0$ & $0.0 $& V \\
        \bottomrule
    \end{tabular}
\end{center}
      \caption{Initial values for state variables in the zero flow limit
      (\textbf{A}) and in the full model (\textbf{B}).
      We use SI base units; that is, meter (m), and mole (mol).}
    \label{tab:init}
\end{table}

\clearpage

\section*{Acknowledgments}

The authors thank Didrik Bakke Dukefoss, Rune Enger, Geir Halnes,
Erlend A. Nagelhus, and Klas Pettersen for valuable and constructive
discussions on cortical spreading depression and brain
electrophysiology.

AJE and MER have received support the European Research Council (ERC)
under the European Union's Horizon 2020 research and innovation
programme under grant agreement 714892. PEF acknowledges support from
the Engineering and Physical Sciences Research Council (grants
EP/R029423/1 and EP/V001493/1). NB acknowledges support from the
Engineering and Physical Sciences Research Council Centre for Doctoral
Training in Industrially Focused Mathematical Modelling (grant
EP/L015803/1) in collaboration with Simula Research Laboratory.

\bibliographystyle{ws-m3as}
\bibliography{references}

\end{document}